\documentclass[12pt]{article}
\usepackage{geometry}
\usepackage{graphicx}
\usepackage{tikz-cd}
\usepackage{kotex}
\usepackage{amssymb,amsmath,amsthm,mathtools}
\usepackage{color}
\usepackage{enumerate}
\usepackage[T1]{fontenc}
\usepackage[utf8]{inputenc} 
\usepackage{bm}
\usepackage{array}
\usepackage{hyperref}
\usepackage{titlesec}
\usepackage{mathscinet}
\usepackage[toc,page]{appendix}
\usepackage{bbold}
\usepackage{wasysym} 
\usepackage{easyReview}
\usepackage{adjustbox}
\usetikzlibrary{shapes.geometric}
\titleformat{\section}[block]{\sffamily\scshape\filcenter\Large}{}{0pt}{}
\renewcommand{\thesubsection}{(\thesection\alph{subsection})}
\titleformat{\subsection}[runin]{\normalfont\bfseries}{\indent\thesubsection\ }{0pt}{}
\renewcommand{\emph}[1]{{\sffamily\textbf{#1}}}

\DeclareMathOperator{\Tot}{Tot}

\DeclareMathOperator{\id}{\mathsf{Id}}

\DeclareMathOperator{\Fun}{\mathsf{Fun}}
\DeclareMathOperator{\diff}{diff}
\DeclareMathOperator{\QCoh}{\mathrm{QCoh}}
\DeclareMathOperator{\Spec}{\mathrm{Spec}}
\DeclareMathOperator{\Mod}{\mathrm{Mod}}
\DeclareMathOperator*{\hocolim}{hocolim}
\DeclareMathOperator*{\holim}{holim}
\DeclareMathOperator{\Res}{\mathrm{Res}}
\DeclareMathOperator{\Coh}{\mathsf{Coh}}
\DeclareMathOperator*{\can}{\mathrm{can}}
\DeclareMathOperator{\Perf}{\mathsf{Perf}}
\DeclareMathOperator{\IndCoh}{IndCoh}
\DeclareMathOperator{\Ch}{\mathsf{Ch}}
\DeclareMathOperator{\Aut}{Aut}
\DeclareMathOperator{\Obj}{Obj}
\DeclareMathOperator{\Ob}{Ob}
\DeclareMathOperator{\field}{\mathbb{k}}
\DeclareMathOperator{\Hom}{\operatorname{Hom}}

\DeclareMathOperator{\dg}{\operatorname{dg}}

\DeclareMathOperator{\VB}{\mathsf{VB}}
\DeclareMathOperator{\vb}{\mathsf{vb}}
\DeclareMathOperator{\str}{\operatorname{str}}

\DeclareMathOperator{\pr}{\operatorname{pr}}
\DeclareMathOperator{\Gm}{\mathbb{G}_{\operatorname{m}}}
\DeclareMathOperator{\cha}{\lambda}
\DeclareMathOperator{\chaa}{\mu}
\DeclareMathOperator{\comp}{\vcenter{\hbox{\scalebox{0.65}{$\bullet$}}}}
\DeclareMathOperator{\Nbd}{Nbd}
\DeclareMathOperator{\Bl}{Bl}
\DeclareMathOperator{\Symp}{Symp}
\DeclareMathOperator{\rF}{\reflectbox{F}}
\DeclareMathOperator{\objectdegree}{{\text{$i$}}}
\DeclareMathOperator{\homdegree}{{\text{$k$}}}
\DeclareMathOperator{\yo}{\text{よ}}
\DeclareMathOperator{\Rep}{Rep}
\DeclareMathOperator{\Ku}{\textrm{Kü}}
\DeclareMathOperator{\Sh}{Sh}
\newtheoremstyle{them}
{4pt}
{4pt}
{\normalfont }
{}
{\sffamily\scshape }
{.}
{.5em}
{}
\newtheoremstyle{exam}
{4pt}
{4pt}
{\normalfont }
{}
{\sffamily\scshape }
{.}
{.5em}
{}
\swapnumbers

\theoremstyle{them}
\newtheorem{thm}{Theorem}[section]
\newtheorem{dfn}[thm]{Definition}
\newtheorem{defn}[thm]{Definition}
\newtheorem{prp}[thm]{Proposition}
\newtheorem{prop}[thm]{Proposition}
\newtheorem{crl}[thm]{Corollary}
\newtheorem{cor}[thm]{Corollary}

\newtheorem{ques}[thm]{Question}

\theoremstyle{exam}
\newtheorem{exmp}[thm]{Example}

\newtheorem{lmm}[thm]{Lemma}
\newtheorem{lemma}[thm]{Lemma}
\newtheorem{rmk}[thm]{Remark}

\title{Local Homological Mirror Symmetry}
\author{Si-Yang Liu\thanks{Department of Mathematics, Imperial College London, London, SW7 2AZ.
  Email: \texttt{s.liu1@imperial.ac.uk}}
  \and Kyungmin Rho\thanks{Max-Planck Institute for Mathematics in Bonn.
  Email: \texttt{rho@mpim-bonn.mpg.de}}
  }
\date{\today}

\begin{document}

\maketitle
\begin{abstract}
    Local mirror symmetry, originated from the work of Chiang-Klemm-Yau-Zaslow \cite{CKYZ99}, indicated a potential relation between mirrors to a vector bundle and mirrors to its zero section. In this paper, we constructed a homological mirror to these local varieties and proved the corresponding Homological Mirror Symmetry statement. As an application, we constructed infinite families of exotic Weinstein sectors as mirrors to Hirzebruch surfaces.
\end{abstract}
\tableofcontents
\section{Introduction}
The Mirror Symmetry phenomena was originally discovered by physicists and later formulated mathematically in various different forms. It was firstly understood as two Calabi-Yau manifolds having mirror symmetric Hodge diamonds, and then as relations between periodic integrals and the Gromov-Witten invariants. In 1994, Kontsevich \cite{Kontsevich95} proposed a framework to understand the mirror symmetry phenomena systematically as an equivalence between triangulated categories, one from coherent complexes of sheaves on an algebraic variety and another from Lagrangian submanifolds of a symplectic manifold. This was known as the \emph{Homological Mirror Symmetry} proposal. More precisely in modern language, given an expectedly mirror symplectic manifold $(X,\omega)$ and an algebraic variety $(Y,J)$, we expect a quasi-equivalence of $A_{\infty}$-categories \[
    \mathcal{F} (X,\omega)\simeq\Coh (Y,J)
\]between a version of the \emph{Fukaya category} of the symplectic manifold $(X,\omega)$ and the dg category of complexes of coherent sheaves on $Y$. Knowing Homological Mirror Symmetry, we are then expected to derive enumerative version of mirror symmetry via purely algebraic constructions \cite{KKP08}.\par
However, Homological Mirror Symmetry proposal didn't address the question of how the mirror pairs are constructed, and the corresponding geometric construction was proposed by Strominger-Yau-Zaslow \cite{SYZ96} in 1996. Roughly speaking, suppose a symplectic manifold $(X,\omega)$ admits a fibration $X\to B$ whose fibres are Lagrangian tori, then we can construct a mirror space $\check{X}\to B$ by dualizing the fibre torus. There will be a natural complex structure $J$ on the dual space in good cases so we get a candidate mirror space $(\check{X},J)$. This construction was then generalized by the famous work of Gross and Siebert \cite{GS06,GS10,GS26} to the situation when $X\to B$ is a Lagrangian torus fibration with possibly elliptic singularities.\par 
In \cite{Fukaya05}, Fukaya proposed an approach to construct a mirror functor $\mathcal{F} (X,\omega)\to\Coh (\check{X},J)$ when $\check{X}$ is obtained from $X$ via SYZ-like construction. Roughly speaking, instead of naively taking dual Lagrangian fibrations, we look at the moduli space of Lagrangian tori (and their ``closures'') objects in the Fukaya category $\mathcal{F} (X,\omega)$, so that the existence of mirror functor is almost by definition. Significant partial progress under this proposal aligning the SYZ construction with Homological Mirror Symmetry was made, by Abouzaid \cite{Abouzaid21} when $X\to B$ is a smooth torus fibration and $\pi_2(B)=0$, and very recently Bai-Haney \cite{BH26} verified a situation when $X\to B$ fits into the Gross-Siebert program with some extra technical assumptions via intrinsic mirror symmetry.\par
In this paper, we look at the mirror pair $(X,\check{X})$ with $A$- and $B$-models swapped. In other words, we want to understand the following question:
\begin{ques}\label{Q:HMS-toric-degeneration}
    Suppose that $(X_t)_{t\in\mathbb{C}}$ is a toric degeneration family and $(Y,\omega_{p(t)})$ is the mirror symplectic manifold to $(X_t, J_t)$. Is there a quasi-equivalence of $A_{\infty}$-categories \[
        \Coh (X_t,J_t)\simeq\mathcal{F} (Y,\omega_{p(t)})?
    \]
\end{ques}
Note that answers to question \ref{Q:HMS-toric-degeneration} when $X_t\equiv X_{\Sigma}$ is a toric variety with moment polytope $\Sigma$ have been answered by Fang-Liu-Treumann-Zaslow \cite{FLTZ12} following an earlier idea of Bondal \cite{Bondal06} and an unpublished work of Bondal-Ruan, where the mirror $(Y,\omega)=\left((\mathbb{C}^{\ast})^n,\Lambda\right)$ is a Landau-Ginzburg model-like symplectic manifold, or more rigorously a \emph{Liouville sector} where the Fukaya category we talk about is the \emph{partially wrapped Fukaya category}, using the language of \cite{GPS20,GPS24M}. Sectorial Descent, established in \cite{GPS24S} and is expected to be mirror to proper descent of coherent complexes, should quickly imply Homological Mirror Symmetry results for the central fibre $\Coh (X_0,J_0)\simeq\mathcal{W} (Y_0,\omega_0)$ of the toric degeneration family. Seidel's proposal \cite{Seidel02} then suggests that $(X_t,J_t)$ should be mirror to a compactification $(Y,\omega_{p(t)})$ of $Y_0$, where $\omega_{p(t)}$ keeps track of the volume near the compatification data.\par
Before proceeding on this direction, there's a sanity problem to check. Given a toric degeneration family $(X_t,J_t)_{t\in\mathbb{C}}$, we can always replace the central fibre $X_0$ by its \emph{log blow-up} $\Bl_Z (X_0)$ along a divisor $Z$ with non-trivial log structure. This operation is local to the fibre at $0$ and should not change its smoothings, and therefore we wonder the following question:
\begin{ques}\label{Q:log-blow-up-mirror}
    What is the mirror space to log blow-ups? Does mirror to log blow-up affect the compactification?
\end{ques}
Answer to the second part of \ref{Q:log-blow-up-mirror} should be clear if we believe in question \ref{Q:HMS-toric-degeneration}: it should not affect the compactification. However, to answer the second part we have to first understand the first question, i.e. to construct the mirror space to log blow-ups. In the simplest situation, suppose that $X_0=Z_1\cup_D Z_2$ is a union of smooth algebraic varieties along a divisor $D$, the log blow-up $\Bl_Z X_0=Z_1\cup\mathbb{P}\left(N_{Z_1} D\oplus N_{Z_2} D\right)\cup Z_2$ is obtained by replacing $D$ by the projectivization of direct sums of normal bundles of $D$ to $Z_1$ and $Z_2$, respectively. It then suffices for us to understand the mirror to the projective bundle $\mathbb{P}\left(N_{Z_1} D\oplus N_{Z_2} D\right)$, and then apply sectorial and proper descent to both models. Therefore the question we need to address is 
\begin{ques}\label{Q:mirror-projective-bundle}
    Given a smooth projective algebraic variety $X$, what is the mirror space to the projectivized bundle $\mathbb{P} (\mathcal{O}_X\oplus\mathcal{L})$ for some line bundle $\mathcal{L}$ on $X$?
\end{ques}
This is in some sense a local question, and the name \emph{local mirror symmetry}, originated from the work of Chiang-Klem-Yau-Zaslow \cite{CKYZ99}, referred to the study of Gromov-Witten curve counts on total space of anti-canonical line bundles over Calabi-Yau surfaces. Observe that for any line bundles $\mathcal{L}_1,\mathcal{L}_2$ over a scheme $X$ we have $\mathbb{P}_X (\mathcal{L}_1\oplus\mathcal{L}_2)=\mathbb{P}_X (\mathcal{L}_1\otimes\mathcal{L}_2^{-1}\oplus\mathcal{O}_X)$, the spaces considered in \cite{CKYZ99} are also close to the cases we are interested in.\par
Before we go to the statement of the main result, let's gi ve intuitions on what the mirror to $\mathbb{P} (\mathcal{O}_X\oplus\mathcal{L})$ should look like. Suppose that $F\to P\to B$ is a general holomorphic Kähler bundle (so fibres and bases are all Kähler manifolds and they produce global Kähler forms on the total space) such that both the fibre $F$ and the base $B$ admits a Lagrangian torus fibration onto some smooth manifolds $\Sigma (F)$ and $\Sigma (B)$, which induces a global torus fibration on $P$. From the SYZ perspective, the mirror should be obtained by dualizing the Lagrangian torus firations, hence resulting in a fibre bundle \[
    \check{B}\to\check{P}\to\check{F}
\]with arrows reversed. Therefore, we expect the mirror $\mathcal{Y}$ to $\mathbb{P} (\mathcal{O}_X\oplus\mathcal{L})$ to be a symplectic fibre bundle over mirror to $\mathbb{P}^1$, which is known by \cite{FLTZ12} to be the Weinstein sector $(\mathbb{C}^{\ast}, z+1/z)$, with fibre $Y$ mirror to the base manifold $X$. Symplectic fibre bundles over $\mathbb{C}^{\ast}$ with given monodromy data $\phi$ was previously studied by Kartal \cite{Kartal21}, who provided a beautiful formula to obtain the Fukaya category $\mathcal{W} (\mathcal{Y})$ from the Fukaya category $\mathcal{W} (Y)$ and the induced auto-equivalence $\phi\colon\mathcal{W} (Y)\to\mathcal{W} (Y)$. We therefore expect that the ``mapping cylinder category'', with the stops equipped, is quasi-equivalent to the coherent complexes of sheaves over the projective bundle $\mathbb{P}\left(\mathcal{O}_X\oplus\mathcal{L}\right)$, once $\phi$ is mirror to tensoring by $\mathcal{L}$. This, however, is different from the mirror symmetry expectation mentioned in \cite{Kartal21}.
\subsection{Statement of the result.} The main result of this paper is a complete answer to question \ref{Q:mirror-projective-bundle} and several implications. Let $Y$ be a smooth quasi-projective variety over $\mathbb{C}$, and let $L(Y,\xi)^{\times}\to Y$ be a principal bundle over $Y$ with fiber the algebraic torus $\mathbb{G}_m$ and with characteristic class $\xi\in CH^2 (Y)$. Let $L(Y,\xi)\to Y$ and $P(Y,\xi)\to Y$ be the associated line and projective line bundle over $Y$. Let $X$ be a symplectic manifold mirror to $Y$ and $\phi\colon X\to X$ a symplectomorphism such that $\phi\colon\mathcal{W} (X)\to\mathcal{W} (X)$ is mirror to the functor $-\otimes L(Y,\xi)^{-1}$ obtained by tensoring with the line bundle $L(Y,\xi)^{-1}$. The mapping cylinder $M(X,\phi)\to\mathbb{C}^{\ast}$ of $\phi$ is a fibre bundle over $\mathbb{C}^{\ast}$, and we write $z$ and $1/z$ to be the composed function $M(X,\phi)\to\mathbb{C}^{\ast}\to\mathbb{C}$.
\begin{thm}\label{thm:HMS-projective-bundle}
    For HMS pairs $(X,Y)$ and mirror automorphisms $(\phi,\xi)$ with $X$ Weinstein sector and $Y$ smooth quasi-projective, there exist a canonical $A_{\infty}$-functor\[
        \Coh\left(P(Y,\xi)\right)_{dg}\xrightarrow{\simeq}\Perf\mathcal{W}\left(M(X,\phi),z+\dfrac{1}{z}\right)^{op}
    \]
    which induces equivalences of triangulated categories after taking homotopy categories.
\end{thm}
It's well expected that inclusions $Y\xrightarrow[i_0]{i_{\infty}} P(Y,\xi)$ of $Y$ into $P(Y,\xi)$ as the zero-section or section at infinity are mirror to the Orlov functors $\mathcal{W} (X)\xrightarrow[\cup_0]{\cup_{\infty}}\mathcal{W} (M(X,\phi),z+1/z)$ sending $X$ to the fibre near $0$ or $\infty$, respectively. With this expectation, we can quickly obtain that 
\begin{crl}\label{crl:mirror-symmetry-line-principal-bundles}
    Under the same convention as Theorem \ref{thm:HMS-projective-bundle}, there exist quasi-equivalences of $A_{\infty}$-categories \[
        \Coh (L(Y,\xi))_{dg}\xrightarrow{\simeq}\Perf\mathcal{W}\left(M(X,\phi),z\right)^{op}
    \]and\[
        \Coh (L(Y,\xi)^{\times})_{dg}\xrightarrow{\simeq}\Perf\mathcal{W} (M(X,\phi))^{op}
        \]between coherent complexes of $\mathbb{G}_m$-principal bundles and line bundles over $Y$ and perfect modules over wrapped Fukaya categories of mapping cylinders built on the pair $(X,\phi)$ possibly with a stop.
\end{crl}
We want to note that, although the theorem is pretty general, we don't know many examples where both the Homological Mirror Symmetry $\Coh (Y)_{dg}\simeq\Perf\mathcal{W} (X)^{op}$ and the mirrors between (graded) symplectomorphisms $\phi\in\Symp (X)$ and derived auto-equivalences $\check{\phi}\colon\Coh (Y)\to\Coh (Y)$ are known. A few known cases are for toric Fano manifolds \cite{FLTZ12,Abouzaid09,Hanlon19,GPS24M} and tensor by toric vector bundles and for log Calabi-Yau surfaces \cite{HK22,HK25} and tensor by any line bundles. Recently the HMS result in Hacking-Keating was further generalized by Varolgunes \cite{Varolgunes26} but it's not known whether we have similar statements about mirror auto-equivalence groups.\par
A very recent work of Mamaev \cite{Mamaev26} identified the category $\Coh (P(Y,\xi))_{dg}$ with Fukaya categories with dg coefficients $\Perf\mathcal{W} ((\mathbb{C}^{\ast},z+z^{-1});\Coh (Y)_{dg})^{op}$, so combining the two results suppposedly gives the quasi-equivalence \[
    \Perf\mathcal{W} (M(X,\phi);z+z^{-1})^{op}\simeq\Perf\mathcal{W} (\mathbb{C}^{\ast},z+z^{-1};\Perf\mathcal{W} (X)^{op})^{op}
\]which verifies in special situations that Fukaya categories of fibre bundles over $(\mathbb{C}^{\ast}, z+z^{-1})$ should correspond to perverse Schobers of Fukaya categories.
\subsection{Multiple monodromies.} In addition to the main mirror symmetry statement, we make some important further quick implications without much effort.\par
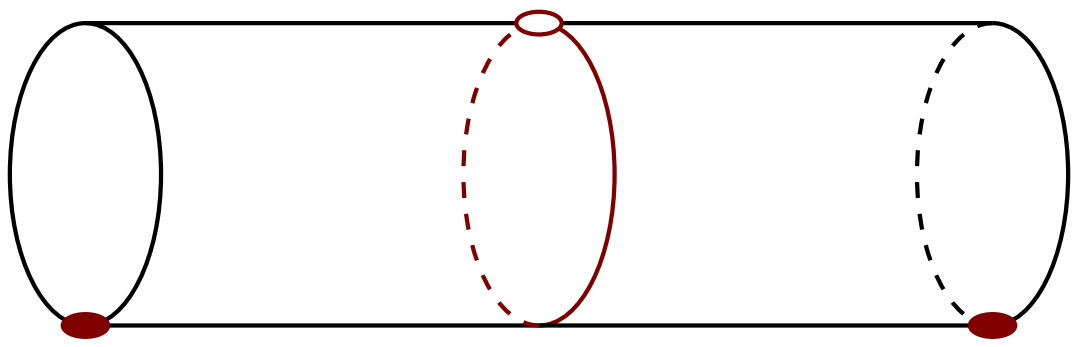
\begin{figure}[ht]
    \centering
    \begin{tikzpicture}
        \begin{scope}
        \draw [thick] (-2,0) arc[start angle=0, end angle=360, x radius=0.5, y radius=1];
        \draw [thick] (-2.5,1) -- (0.5,1);
        \draw [thick] (-2.5,-1) -- (0.5,-1);
        \draw [thick,red!50!black] (0.5,1) arc[start angle=90, end angle=-90, x radius=0.5, y radius=1];
        \draw [thick,red!50!black, dashed] (0.5,-1) arc[start angle=270, end angle=90, x radius=0.5, y radius=1];
        \end{scope}
        \node [ellipse,
    minimum width=3mm,
    minimum height=1.5mm,
    inner sep=0pt,
    fill=red!50!black,
    draw=red!50!black,
    thick] at (-2.5,-1) {};
    \begin{scope}[xshift=3cm]
        \draw [thick] (-2.5,1) -- (0.5,1);
        \draw [thick] (-2.5,-1) -- (0.5,-1);
        \draw [thick] (0.5,1) arc[start angle=90, end angle=-90, x radius=0.5, y radius=1];
        \draw [thick, dashed] (0.5,-1) arc[start angle=270, end angle=90, x radius=0.5, y radius=1];
        \node [ellipse,
    minimum width=3mm,
    minimum height=1.5mm,
    inner sep=0pt,
    fill=red!50!black,
    draw=red!50!black,
    thick] at (0.5,-1) {};
    \end{scope}
    \node [ellipse,
    minimum width=3mm,
    minimum height=1.5mm,
    inner sep=0pt,
    fill=white,
    draw=red!50!black,
    thick] at (0.5,1) {};
    \end{tikzpicture}
    \caption{Gluing $A$-models to get $\mathcal{X} (\phi_1,\phi_2)$.}\label{fig:gluing-A-model-bundle}
\end{figure}

Given two characteristic classes $\xi_1, \xi_2$ on $Y$ with mirror symplectomorphisms $\phi_1, \phi_2\in\Symp (X)$, we can form the glued space \[
    \mathcal{Y} (\xi_1, \xi_2)\coloneqq P(Y,\xi_1)\cup_Y P(Y,\xi_2)
\]from the projective bundles $P(Y,\xi_i)$ along the zero-section and section-at-infinity. By construction $\mathcal{Y} (\xi_1,\xi_2)$ is still a fibre bundle over $Y$ with fibres nodal curves $\mathbb{P}^1\cup\mathbb{P}^1$, and since each $P(Y,\xi_i)$ is an irreducible component of $\mathcal{Y}$, we have the proper descent diagram \[
\begin{tikzcd}
    \Coh (Y)_{\infty}\arrow[r,"i_{\infty}"]\arrow[d,"i_0"]& \Coh\left(P(Y,\xi_2)\right)_{\infty}\arrow[d]\\ 
    \Coh\left(P(Y,\xi_1)\right)_{\infty}\arrow[r]&\Coh\left(\mathcal{Y} (\xi_1,\xi_2)\right)_{\infty}
\end{tikzcd}
\]so that $\Coh (\mathcal{Y} (\xi_1, \xi_2))_{\infty}$ is completely determined by the homotopy pushforward diagram. This pushforward diagram is mirror to the homotopy pushforward diagram \[
    \begin{tikzcd}
        \mathcal{W} (X)\arrow[r,"\cup_0"]\arrow[d,"\cup_{\infty}"]&\mathcal{W} (M(X,\phi_2),z_2+z_2^{-1})\arrow[d]\\ 
        \mathcal{W} (M(X,\phi_1),z_1+z_1^{-1})\arrow[r]&\mathcal{W} (\mathcal{X} (\phi_1,\phi_2))
    \end{tikzcd}
\]where $\mathcal{X} (\phi_1,\phi_2)$ is the Weinstein sector obtained from gluing $\left(M(X,\phi_1),z_1+z_1^{-1}\right)$ and $\left(M(X,\phi_2),z_2+z_2^{-1}\right)$ along the identification of fibres-at-infinity $\{z_1=\infty\}\cong\{z_2=0\}$. What we obtain will be a sector which after removing the stops admits a smooth fibration over the twice-punctured disk $\mathbb{D}^2\setminus\{\pm 1/2\}$, with stops on one of the punctures and at infinity. Figure \ref{fig:gluing-A-model-bundle} shows how the two mapping cylinders are glued via the gluing on the base. In general, we can prove the following mirror symmetry statement.
\begin{crl}\label{crl:gluing-mirror-symmetry}
    Given HMS pairs $(X,Y)$ with $X$ Weinstein sector and $Y$ smooth quasi-projective and mirror automorphisms $(\phi_i,\xi_i)$, $1\leq i\leq k$, we have quasi-equivalences of $A_{\infty}$-categories \[
        \Perf\mathcal{W}\left[\mathcal{X} (\phi_1,\dotsb ,\phi_k)\right]^{op}\simeq\Coh\mathcal{Y} (\xi_1,\dotsb ,\xi_k)_{dg}
    \]where $\mathcal{X} (\phi_1,\dotsb,\phi_k)$ is obtained from gluing $k$ copies of mapping cylinders $\left(M(X,\phi_i),z_i+z_i^{-1}\right)$, $\mathcal{Y} (\xi_1,\dotsb ,\xi_k)$ from gluing $k$ copies of projective bundles $P(Y,\xi_i)$, under the rule provided by the graph\par
    
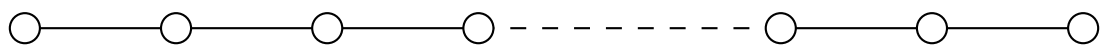
\begin{center}
        \begin{tikzpicture}
            \draw (0,0) -- (3,0);
            \draw [dashed] (3,0) -- (5,0);
            \draw (5,0) -- (7,0);
            \foreach \x in {0,1,2,3,5,6,7} {
            \node [circle, fill=white, draw=black, inner sep=2pt] at (\x,0) {};
            }
        \end{tikzpicture}
    \end{center}
    where vertices represent the spaces $\left(M(X,\phi_i),z_i+z_i^{-1}\right)$ and $P(Y,\xi_i)$, and two vertices are connected by an edge if their representing objects are glued along the fibre/divisor.
\end{crl}
Naturally, we can also connect the first and the last vertices in Corollary \ref{crl:gluing-mirror-symmetry} to obtain mirror symmetries for ``mapping punctured tori'':
\begin{crl}\label{crl:gluing-mirror-symmetry-circ}
    With the same convention as Corollary \ref{crl:gluing-mirror-symmetry}, let \[
        \mathcal{X}_{\circ} (\phi_1,\dotsb,\phi_k)\quad\text{and}\quad\mathcal{Y}_{\circ} (\xi_1,\dotsb,\xi_k)
    \]be Weinstein sectors and algebraic varieties, respectively, obtained by gluing together Weinstein sectors $(M(X,\phi_i),z_i+z_i^{-1})$, resp. $\mathbb{P} (\mathcal{O}_Y\oplus\mathcal{L}_{\xi_i})$ with the same rule as in Corollary \ref{crl:gluing-mirror-symmetry}, with respect to the graph
    
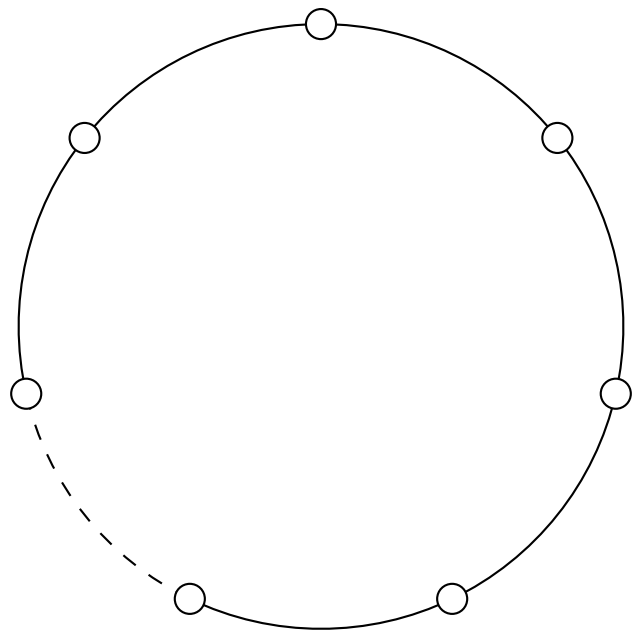
\begin{center}
    \begin{tikzpicture}
        \def\R{2}
        \coordinate (v1) at (90:\R);
        \coordinate (v2) at (38.571:\R);
        \coordinate (v3) at (-12.857:\R);
        \coordinate (v4) at (-64.286:\R);
        \coordinate (v5) at (-115.714:\R);
        \coordinate (v6) at (-167.143:\R);
        \coordinate (v7) at (-218.571:\R);
        \draw (v1) arc (90:38.571:\R);
        \draw (v2) arc (38.571:-12.857:\R);
        \draw (v3) arc (-12.857:-64.286:\R);
        \draw (v4) arc (-64.286:-115.714:\R);
        \draw[dashed] (v5) arc (-115.714:-167.143:\R);
        \draw (v6) arc (-167.143:-218.571:\R);
        \draw (v7) arc (-218.571:-270:\R);
        \foreach \vertex in {v1,v2,v3,v4,v5,v6,v7} {
            \node[circle, fill=white, draw=black, inner sep=2pt]
                at (\vertex) {};
        }
    \end{tikzpicture}
\end{center}
Then we have quasi-equivalences of $A_{\infty}$-categories \[
    \Perf\mathcal{W}\left[\mathcal{X}_{\circ} (\phi_1,\dotsb,\phi_k)\right]^{op}\simeq\Coh\mathcal{Y}_{\circ} (\xi_1,\dotsb,\xi_k)_{dg}.
\]
\end{crl}
Together with some deformation argument modified from \cite{LK17}, we should be able to achieve mirror symmetry for fibre bundles with fibres elliptic curves.\par
Aside from gluing pieces together, one can also take direct sums of line bundles on $B$-side and then fibrewise compactify to a toric bundle, while on the $A$-side taking products of $T^{\ast}\mathbb{S}^1$s and putting related stops on the base, using the rules described in \cite{FLTZ12,Kuwagaki20}.\par
More precisely, given mirror pairs $(X;\phi_1,\dotsb ,\phi_k)$ on the $A$-side and $(Y;\xi_1,\dotsb, \xi_k)$ on the $B$-side where $(\phi_1,\dotsb ,\phi_k)$ are mutually commutative symplectomorphisms, together with an integral convex polytope $\Sigma\subseteq\mathbb{R}^k$, one can construct \[
    \mathcal{X}_{\Sigma} (X;\phi_1,\dotsb ,\phi_k )\coloneqq\left(T^{\ast}\mathbb{T}^n,\Lambda_{\Sigma}\right)\times_{\phi_1,\dotsb ,\phi_n} X\to (T^{\ast}\mathbb{T}^n,\Lambda_{\Sigma})
\]as a fibre bundle with fibre $X$ and monodromy around each generator of $\pi_1(\mathbb{T}^n)$ $\phi_i$, and $\mathcal{Y}_{\Sigma} (Y;\xi_1,\dotsb ,\xi_k)$ a fibre bundle over $Y$ with fibres $X_{\Sigma}$, where $(\xi_1,\dotsb ,\xi_k)$ characterizes its structure principal bundle.
\begin{crl}\label{crl:toric-bundle-HMS}
    With the same convention as Corollary \ref{crl:gluing-mirror-symmetry} and let $\Sigma$ be an integral convex polytope in $\mathbb{R}^k$, there is a quasi-equivalence of $A_{\infty}$-categories \[
    \Perf\mathcal{W}\left[\mathcal{X}_{\Sigma} (X;\phi_1,\dotsb ,\phi_k)\right]^{op}\simeq\Coh\left[\mathcal{Y}_{\Sigma} (Y;\xi_1,\dotsb ,\xi_k)\right]_{dg}.
    \]
\end{crl}
\subsection{Hirzebruch surfaces.} As another application, we looked at a particular example of projective toric surfaces that are $\mathbb{P}^1$-bundles over $\mathbb{P}^1$, which are known as \emph{Hirzebruch surfaces}. More precisely, given a non-negative integer $n$, the Hirzebruch surface $\mathbb{F}_n$ is the projectization \[
    \mathbb{F}_n =\mathbb{P}\left(\mathcal{O}_{\mathbb{P}^1}\oplus\mathcal{O}_{\mathbb{P}^1}(n)\right).
\]of the rank $2$ vector bundle $\mathcal{O}\oplus\mathcal{O}_{\mathbb{P}^1} (n)$ over $\mathbb{P}^1$. Since $\mathbb{P}^1$ is toric, $\mathbb{F}_n$ are toric for all $n$.\par
Note that although coherent-constructible correspondence for all toric varieties are known, this result is not compatible with the mirror construction of Hori-Vafa \cite{HV00} when the toric variety is not Fano. For instance, the number of critical points of the Hori-Vafa superpotential $W_n$ is not the same as the rank of the quantum cohomology of the mirror $\mathbb{F}_n$ \cite[Example 8.2]{FOOO10}. The SYZ mirror construction, as explained by Auroux \cite{Auroux09} and later studied by Jing \cite{Jing26}, is highly non-trivial and involves possibly gluing infinitely many charts together.\par
Using theorem \ref{thm:HMS-projective-bundle}, we provide a rather direct description of the homological mirror to $\mathbb{F}_n$.
\begin{crl}\label{crl:mirror-to-Hirzebruch}
    Write \[
        \rF_n =\left(M((\mathbb{C}^{\ast}, x+1/x),\phi^n),z+1/z\right)
    \]to be a Landau-Ginzburg model, where $\phi\colon\mathbb{C}^{\ast}=T^{\ast}\mathbb{S}^1\to T^{\ast}\mathbb{S}^1$ is the Dehn twist along the zero section. Then we have diffeomorphisms as manifolds with boundaries and isomorphisms as holomorphic manifolds with superpotentials\[
        \rF_n\cong\left(\mathbb{C}^{\ast}\times\mathbb{C}^{\ast}, x+1/x+y+1/y\right)
    \]and there is a quasi-equivalence of $A_{\infty}$-categories \[
        \Perf\mathcal{W} (\rF_n)^{op}\simeq\Coh\left(\mathbb{F}_n\right)_{dg}.
    \]
\end{crl}
In particular, from the work of Bondal-Orlov \cite{BO01} and Bridgeland-Maciocia \cite{BM01}, we know that $\Coh (\mathbb{F}_n)\cong\Coh (\mathbb{F}_m)$ iff $m=n$, and therefore we immediately get the following possibly interesting Corollary:
\begin{crl}\label{crl:exotic-Weinstein-structure}
    $\rF_n$ is exact symplectomorphic to $\rF_m$ iff $m=n$.
\end{crl}
Therefore although $\rF_n$s are all diffeomorphic and isomorphic as Landau-Ginzburg pairs, they're not exact symplectomorphic and are therefore a sequence of ``exotic'' Weinstein sectors, since exact symplectomorphisms preserve Fukaya categories \cite{HS26}. However, we should note that Fukaya categories cannot distinguish their convex completions and therefore we have no knowledge of whether this creates exotic Weinstein structures on $\left(\mathbb{C}^{\ast}\right)^2$. Another mystery is about the relation between local homological mirror and the Hori-Vafa construction, i.e.
\begin{ques}\label{Q:HV-vs-LHMS}
    What is the relation between the local homological mirror $\rF_n$ and the Hori-Vafa mirror $\left((\mathbb{C}^{\ast})^2,W_n\right)$?
\end{ques}
Answering this question would hopefully provide the desired SYZ mirror for Hirzebruch surfaces via the family Floer program.
\subsection{Structure of the paper.} In section \ref{sec:proof-of-theorems} we provide all the proofs to the statements in the main theorem with technical details postponed to later sections. In section \ref{sec:algebraic-preparations} we give basic definitions of the higher algebras and corresponding operations that we will use throughout the paper. Section \ref{sec:Fukaya-category-discrete-quotient} is devoted to some of the new structural results we need for the $A$-side computation, while sections \ref{sec:G-equiv-sheaves} and \ref{sec:kunneth-formula} establish necessary materials for the $B$-side computation. More explicitly, section \ref{sec:G-equiv-sheaves} provides basic properties of the category of $G$-equivariant sheaves, and section \ref{sec:kunneth-formula} establishes the key formula: the K\"unneth formula for equivariant sheaves. In the appendix \ref{sec:Eilenberg-Zilber-constraction}, we recall the classical Eilenberg-Zilber construction and adapted the construction to our situation, which will be used in the proof of the K\"unneth formula \ref{thm:Kunneth-embedding}.
\subsection*{Acknowledgement.} The main idea of this work arose during the authors' discussion wth Tim Hu and Maximilian Schimpf during their stay in Bonn/Max-Planck-Insititut in late 2025, and we thank all the useful discussions with them.
We also want to thank Yanki Lekili, Cheuk Yu Mak, Longke Tang and Bingyu Zhang for many useful discussions and suggestions.
S.-Y. Liu is funded by the Deutsche Forschungsgemeinschaft (DFG, German Research Foundation) under Germany's Excellence Strategy - GZ 2047/1, Projekt-ID 390685813.
Kyungmin Rho is grateful to Max Planck Institute for Mathematics in Bonn for its hospitality and financial support.
\subsection*{AI disclosure.} The result \cite{BM01} was found with the aid of LLMs when the authors were trying to find the proof that coherent complexes over non-isomorphic Hirzebruch surfaces are not quasi-isomorphic. The authors believed that there should be a better reference than the overly strong result of \cite{BM01}, while the LLM claimed it to be essential.\par
The authors used LLMs to explore whether the definition of ``small'' Fukaya category in \ref{defn:wrapped-Fukaya-cat-inf-sectors} is the same as the category of compact objects in $\mathcal{W}^{big} (X)$, without a success. After discussion with Bingyu Zhang and Longke Tang, the authors were convinced that some necessary geometric argument was needed and decided to leave this as a further question to investigate.\par
The reference \cite{Christ26} relating homotopy quotients of an $A_{\infty}$-category $\mathcal{A}$ by a group $G$ to the $A_{\infty}$-orbit category was also found with the aid of LLMs.\par
In the proof of Corollary \ref{crl:gluing-mirror-symmetry}, the authors asked LLM to separate the first huge equation into two rows. The original code was input into LLM and the output was checked and copied directly to the file.\par
The authors tried to find further refernces for the proof of Proposition \ref{prp:dg-enhancement}, and the LLM didn't provide answers outside the authors' knowledge.\par
The authors used LLM to turn the right module action of Seidel \cite{Seidel08} into a left module action, checked and copied the generated codes to the file.
\section{Proof of the Theorems\label{sec:proof-of-theorems}}
In this section, with technical details postponed to later sections, we will provide proofs of the claimed results in the Introduction.
\subsection{The preliminary and the blueprint.} Throughout this paper, on the $A$-model we will regard all symplectic manifolds $X$ as \emph{Liouville sectors}, which will be defined in Definition \ref{defn:Liouville-sector-inf}. For Liouville sectors of finite type their wrapped Fukaya categories are defined in \cite{GPS20}, which share a lot of good functorial properties explained in \cite{GPS24S}.\par
Given an exact symplectomorphism $\phi\colon X\to X$ of the Liouville sector $X$, we can define the \emph{mapping cylinder} \[
    M (X,\phi )\coloneqq X\times_{\phi}\mathbb{C}^{\ast}\coloneqq\frac{T^{\ast}\mathbb{R}\times X}{(z,x)\sim (z+1,\phi (x))}
\]which can be regarded as the quotient of the product space $X\times\mathbb{C}$ by a $\mathbb{Z}$-action which acts by $\phi$ on $X$ and by translation by $1$ on $\mathbb{C}$. Its Weinstein structure is induced from $X\times\mathbb{C}$ following definition \ref{defn:G-invariant-Weinstein-sector}, and more explicitly by performing a modification to the product Weinstein structure, using an idea of Kartal \cite{Kartal21}. This modification is done in Proposition \ref{prp:Kartal-modification-Weinstein}. The wrapped Fukaya category of $T^{\ast}\mathbb{R}$ with stops given by preimages of $e^z+e^{-z}$ near infinity was computed in Example \ref{exmp:mirror-to-Coh-G-P1}, and we have
\begin{prp}\label{prp:Fukaya-cat-mapping-cylinder}
    $\mathcal{W} (M(X,\phi),z+1/z)\simeq\left(\mathcal{W} (T^{\ast}\mathbb{R}, e^z+e^{-z})\otimes\mathcal{W} (X)\right)/\mathbb{Z}$.
\end{prp}
\begin{proof}
    By Künneth formula \ref{crl:Kunneth-Fukaya-infinite} combined with the definition of tensor products of $A_{\infty}$-categories in \cite{Ornaghi25}, we otain a quasi-equivalence of $A_{\infty}$-categories \[
        \mathcal{W} (T^{\ast}\mathbb{R},e^z+e^{-z})\otimes\mathcal{W} (X)\hookrightarrow\mathcal{W} (T^{\ast}\mathbb{R}\times X, e^z+e^{-z}).
    \]By theorem \ref{thm:wrapped-Fukaya-quotient}, as we have modified the Weinstein structure on $T^{\ast}\mathbb{R}\times X$ so that it's $\mathbb{Z}$-invariant, we know that \[
        \mathcal{W} (M(X,\phi),z+z^{-1})=\mathcal{W}\left((T^{\ast}\mathbb{R}\times X,e^z+e^{-z})/\mathbb{Z}\right)\simeq\mathcal{W} (T^{\ast}\mathbb{R}\times X,e^z+e^{-z})/\mathbb{Z}
    \]which concludes the Proposition.
\end{proof}
On the mirror side, we start with an algebraic variety $Y$ together with a choice of characteristic classes $\xi\in CH^2 (Y)$ that defines a holomorphic line bundle $L_{\xi, Y}$ on $Y$. The subbundle $L_{\xi,Y}^{\times}$ removing the zero section $Y\hookrightarrow L_{\xi, Y}$ is then a principal $\mathbb{G}_m$-bundle over $Y$, and hence carries a free algebraic $\mathbb{G}_m$-action on $Y$. Then we have
\begin{thm}\label{thm:twisted-tensor-P1-bundle}
    Suppose that $Y$ is smooth and quasi-projective, then there is a quasi-equivalence of dg categories\footnote{Here the definition of $\Perf$ is different from the usual definition. See Definition \ref{dfn:dg-enhancement-Coh}.}
$$
\Perf_{\mathbb{G}_m}(\mathbb{P}^1)_{\dg}
\underset{((1),\mathcal{L}(\xi)\otimes-)}{\otimes}
\Perf(Y)_{\dg}
\rightarrow
\Perf(P(Y,\xi))_{\dg}.
$$
\end{thm}

\begin{proof}
    Note that \[
        Y=\left[\frac{L(Y,\xi)^{\times}}{\mathbb{G}_m}\right],
    \]and hence we have quasi-equivalences $\Perf (Y)_{dg}\simeq\Perf_{\mathbb{G}_m} (L(Y,\xi)^{\times})_{dg}$. Note that the action $((1),\mathcal{L} (\xi)\otimes -)$ on $\Perf_{\mathbb{G}_m} (\mathbb{P}^1)_{dg}\otimes\Perf_{\mathbb{G}_m} (L(Y,\xi)^{\times})_{dg}$ is diagonal, by twisted K\"unneth formula we have quasi-embedding \[
        \Perf_{\mathbb{G}_m} (\mathbb{P}^1)\otimes_{\mathbb{Z}=\widehat{\mathbb{G}}_m}\Perf_{\mathbb{G}_m} (L(Y,\xi)^{\times})_{dg}\hookrightarrow\Perf_{\mathbb{G}_m} (\mathbb{P}^1\times L(Y,\xi)^{\times} )_{dg}\cong\Perf (P(Y,\xi))_{dg},
    \]and it suffices to show that this quasi-embedding is essentially surjective. The work of Orlov \cite{Orlov92} showed that in this situation there is a semi-orthogonal decomposition of triangulated categories \[
    D\Perf (P(Y,\xi))_{dg}\simeq\langle D\Perf(Y)_{-1}, D\Perf(Y)_0\rangle,
    \]which implies that objects of the form $\pi^{\ast}\mathcal{F}$ and $\pi^{\ast}\mathcal{F}\otimes\mathcal{O}_{P(Y,\xi)} (1)$, where $\mathcal{O}_{P(Y,\xi)} (1)$ is the invertible sheaf whose restriction to each fibre is $\mathcal{O}_{\mathbb{P}^1} (1)$ and $\mathcal{F}\in\Perf (Y)$, generate $\Perf(P(Y,\xi))_{dg}$. Note that these objects are images of the objects $\mathcal{O}_{\mathbb{P}^1} (i)\widetilde{\otimes}\mathcal{F}$ in\\$\Perf_{\mathbb{G}_m} (\mathbb{P}^1)_{dg}\widetilde{\otimes}\Perf (Y)_{dg}$, and therefore the image of the given quasi-embedding generates the target dg category.
\end{proof} 
It's clear from \cite{FLTZ12} and Example \ref{exmp:mirror-to-Coh-G-P1} that $\Perf_{\mathbb{G}_m} (\mathbb{P}^1)_{dg}\simeq\Perf\mathcal{W} (T^{\ast}\mathbb{R},e^z+e^{-z})$, and therefore we have collected enough materials to conclude the proof.
\begin{proof}[Proof of Theorem \ref{thm:HMS-projective-bundle}]
    Since by assumption $\Perf\mathcal{W} (X)^{op}\simeq\Coh (Y)_{dg}=\Perf (Y)_{dg}$ and from the previous discussions $\Perf\mathcal{W} (T^{\ast}\mathbb{R},e^z+e^{-z})^{op}\simeq\Perf_{\mathbb{G}_m} (\mathbb{P}^1)_{dg}$, we have quasi-equivalences between their tensor products \[
        \Perf\mathcal{W} (T^{\ast}\mathbb{R},e^z+e^{-z})^{op}\otimes\Perf\mathcal{W} (X)^{op}\simeq\Perf_{\mathbb{G}_m} (\mathbb{P}^1)_{dg}\otimes\Perf (Y)_{dg}.
    \]On the $A$-model the $\mathbb{Z}$-action on the product is provided by shifting by $1$ on $T^{\ast}\mathbb{R}$ and by the symplectomorphism $\phi$ on $\mathcal{W} (X)$, which coincides with the categorical $\mathbb{Z}$-action on the $B$-side coming from the $\mathbb{G}_m$-action. Therefore after taking $\mathbb{Z}$-orbit categories and combining with Proposition \ref{prp:Fukaya-cat-mapping-cylinder} and Theorem \ref{thm:twisted-tensor-P1-bundle} we conclude the proof.
\end{proof}
\begin{proof}[Proof of Corollary \ref{crl:mirror-symmetry-line-principal-bundles}]
    It suffices to show that the inclusion functor $i_0\colon\Coh (Y)\to\Coh (P(Y,\xi))$ is mirror to the Orlov functor $\cup_0\colon\mathcal{W} (X)\to\mathcal{W} (M(X,\phi),z+z^{-1})$. Note that by construction, $i_0$ can be factored as \[
        \Perf (Y)_{dg}\xrightarrow{i_0}\Perf (Y)_{dg}\otimes\Perf_{\mathbb{G}_m}(\mathbb{P}^1)_{dg}\to\Perf (P(Y,\xi))_{dg},
    \]and similarly $\cup_0$ can be factored as \[
        \Perf\mathcal{W} (X)^{op}\xrightarrow{\cup_0}\Perf\mathcal{W} (X\times (T^{\ast}\mathbb{R},e^z+e^{-z}))^{op}\to\Perf\mathcal{W} (M(X,\phi),z+z^{-1})^{op},
    \]where the second map are all quotient by the $\mathbb{Z}$-action, which we know to be mirror to each other. On the $B$-side the functor $i_0$ is given by sending any $\mathcal{F}\in\Coh (Y)$ to the sheaf $\mathcal{F}\otimes\mathcal{O}_0 [0]$, while on the $A$-side $\cup_0$ sends $L\in\mathcal{W} (X)$ to $L\otimes D[0]$, where $D[0]\in\Perf\mathcal{W} (T^{\ast}\mathbb{R},e^z+e^{-z})^{op}\simeq\Perf_{\mathbb{G}_m} (\mathbb{P}^1)_{dg}$ is the linking disk to the stop labelled as $0$, which up to a shift is identified with $\mathcal{O}_0 [0]$. Therefore $i_0$ and $\cup_0$ are mirror to each other and by taking localizations we conclude the proof. The proof for principal bundles is basically the same, with $0$ replaced by $\infty$.
\end{proof}
\subsection{Sectorial and proper descents.} Now we move forward to the proof of the theorems \ref{crl:gluing-mirror-symmetry} and \ref{crl:gluing-mirror-symmetry-circ}. Although stated in a pretty general way, it suffices for us to understand the case when $k=2$, which will follow from the discussion of Example \ref{exmp:wrapped-Fukaya-cat-twice-punctured-equivariant}.\par
Given $B$-model data $(Y,\xi_1,\xi_2)$, we construct the space \[
    \mathcal{Y} (\xi_1,\xi_2)=P(Y,\xi_1)\cup_Y P(Y,\xi_2)
\]by gluing the two projective bundles by identifying the section at $0$ and the section at $\infty$. Note that there is a $\mathbb{G}_m^2$-action on $\mathcal{Y} (\xi_1,\xi_2)$, given by the $\mathbb{G}_m$-actions on each $P(Y,\xi_i)$ that extend trivially to the other component. From the proper pushforward diagram \[
    \begin{tikzcd}
        Y\arrow[r,"i_0"]\arrow[d,"i_{\infty}"]&P(Y,\xi_2)\arrow[d]\\ 
        P(Y,\xi_1)\arrow[r]&\mathcal{Y} (\xi_1,\xi_2)
    \end{tikzcd}
\]we get directly the proper pushforward diagram of quotient stacks\footnote{This follows directly from the definition of quotient stacks and valuative criterion. See e.g. \cite{Alper26}.} \[
    \begin{tikzcd}
        {[Y/\mathbb{G}_m^2]}\arrow[r,"i_0"]\arrow[d,"i_{\infty}"]&{\big[P(Y,\xi_2)/\mathbb{G}_m^2\big]}\arrow[d]\\ 
        \big[P(Y,\xi_1)/\mathbb{G}_m^2\big]\arrow[r]&\big[\mathcal{Y}(\xi_1,\xi_2)/\mathbb{G}_m^2\big]
    \end{tikzcd}
\]The homotopy pushforward formula for quasi-coherent sheaves for pushforward diagrams of closed immersions, which follows from Gaitsgory-Rozenblyum \cite[Theorem A.1.2]{GR17II} together with descent property for quasi-coherent sheaves, implies the homotopy pushforward diagram \[
    \begin{tikzcd}
        {\IndCoh [Y/\mathbb{G}_m^2]_{\infty}}\arrow[r,"i_0"]\arrow[d,"i_{\infty}"]&{\IndCoh\big[P(Y,\xi_2)/\mathbb{G}_m^2\big]_{\infty}}\arrow[d]\\ 
        \IndCoh\big[P(Y,\xi_1)/\mathbb{G}_m^2\big]_{\infty}\arrow[r]&\IndCoh\big[\mathcal{Y}(\xi_1,\xi_2)/\mathbb{G}_m^2\big]_{\infty}
    \end{tikzcd}
\]and since coherent sheaves are pushed to coherent sheaves via closed immersions, it follows that the same diagram for coherent sheaves is also a homotopy pushforward diagram. Using the quasi-equivalence between $dg$- and $\infty$-categories proved in \cite{Cohn16}, we have
\begin{prp}\label{prp:pushforward-coherent-stack}
    There is a quasi-equivalence of $A_{\infty}$-categories \[
        \Coh_{\mathbb{G}_m^2}\mathcal{Y} (\xi_1,\xi_2)_{dg}\simeq\hocolim\left(\Perf_{\mathbb{G}_m^2} (P(Y,\xi_1))_{dg}\leftarrow\Perf_{\mathbb{G}_m^2} (Y)_{dg}\rightarrow\Perf_{\mathbb{G}_m^2} (P(Y,\xi_2))_{dg}\right).
    \]
\end{prp}
Now we turn to the mirror. Since the extra $\mathbb{G}_m$-action is also identity on $Y$, we have that $[P(Y,\xi_i)/\mathbb{G}_m^2]=[\mathbb{P}^1/\mathbb{G}_m]\times Y\times [\mathrm{pt}/\mathbb{G}_m]$. As $\Perf_{\mathbb{G}_m} (\mathbb{P}^1)_{dg}\simeq\Perf\mathcal{W} (T^{\ast}\mathbb{R},e^z+e^{-z})$, it suffices to understand Homological Mirror Symmetry for $[\mathrm{pt}/\mathbb{G}_m]$.
\begin{lmm}
    There is a quasi-equivalence of $A_{\infty}$-categories \[
        \Perf\mathcal{W}\left(\bigsqcup_{i\in\mathbb{Z}}T^{\ast}[0,1]\right)^{op}\simeq\Perf_{\mathbb{G}_m} (\mathrm{pt})_{dg}.
    \]
\end{lmm}
\begin{proof}
    Note that the left-hand side is the homotopy coproduct \[
       \Perf \mathcal{W}\left(\bigsqcup_{i\in\mathbb{Z}}T^{\ast} [0,1]\right)^{op}\simeq\coprod_{i\in\mathbb{Z}}\Perf\mathcal{W} (T^{\ast} [0,1])^{op}\simeq \coprod_{i\in\mathbb{Z}}\Perf (\mathrm{pt})_{dg}\simeq\coprod_{i\in\mathbb{Z}}\Perf\field\text{-}\Mod_{dg},
    \]while on the other hand, from \cite[Proposition 4.9]{TV18} we know that the $\mathbb{G}_m$-equivariant quasi-coherent category $\QCoh_{\mathbb{G}_m} (\mathrm{pt})_{dg}\simeq\Rep_{\Bbbk} (\mathbb{G}_m)_{dg}$ is equivalent to the category of representations of $\mathbb{G}_m$, and since a linear $\mathbb{G}_m$-representation is the same as a $\mathbb{Z}$-grading, we see that $\Rep_{\Bbbk} (\mathbb{G}_m)_{dg}$ is equivalent to category of complexes of $\field$-vector spaces equipped with an extra $\mathbb{Z}$-grading. By taking compact objects, we know that \[
        \Perf_{\mathbb{G}_m} (\mathrm{pt})_{\infty}\simeq\left(\Rep_{\Bbbk} (\mathbb{G}_m)_{\infty}\right)^{\omega}\simeq\coprod_{i\in\mathbb{Z}}\Perf\Bbbk\text{-}\Mod_{\infty}
    \]using the fact that $\Perf_{\mathbb{G}_m} (\mathrm{pt})_{dg}$ is the compact subcategory of $\QCoh_{\mathbb{G}_m} (\mathrm{pt})_{dg}$ \cite{BZFN10}. We then use the equivalence between $dg$-categories and $\Bbbk$-linear stable $\infty$-categories \cite{Cohn16} to conclude the proof.
\end{proof}
From this lemma and the previous discussions, we get a quasi-equivalence of $A_{\infty}$-categories \begin{align}\label{eq:equivariant-mirror-Gm2}
    \begin{split}
    \Perf\mathcal{W}\left(\bigsqcup_{i\in\mathbb{Z}} (T^{\ast}\mathbb{R},e^z+e^{-z})\right)^{op}&\simeq\Perf\mathcal{W} (T^{\ast}\mathbb{R},e^z+e^{-z})^{op}\otimes\Perf_{\mathbb{G}_m} (\mathrm{pt})\\ 
    &\simeq\Perf_{\mathbb{G}_m} (\mathbb{P}^1)\otimes\Perf_{\mathbb{G}_m} (\mathrm{pt})\simeq\Perf_{\mathbb{G}_m^2} (\mathbb{P}^1).
    \end{split}
\end{align}
\begin{proof}[Proof of Corollary \ref{crl:gluing-mirror-symmetry}]
    As already discussed, it suffices for us to verify \[
        \Perf\mathcal{W}\left(\mathcal{X} (\phi_1,\phi_2)\right)^{op}\simeq\Coh\left(\mathcal{Y}(\xi_1,\xi_2)\right)
    \]for chosen mirror data $(X,\phi_1,\phi_2)$ and $(Y,\xi_1,\xi_2)$. Equivalence \eqref{eq:equivariant-mirror-Gm2} together with analyse of mirror functors provides a quasi-equivalence between the diagrams \[
\begin{aligned}
&\Perf\left(
\begin{tikzcd}
\mathcal{W}\left(\bigsqcup_{\mathbb{Z}^2}T^{\ast}[0,1]\right)
  \arrow[r]\arrow[d]
&
\mathcal{W}\left(\bigsqcup_{\mathbb{Z}}(\mathbb{C},e^{z_1}+e^{-z_1})\right)
\\
\mathcal{W}\left(\bigsqcup_{\mathbb{Z}}(\mathbb{C},e^{z_2}+e^{-z_2})\right)
&
\end{tikzcd}
\right)^{op}
\\[1em]
&\simeq
\left(
\begin{tikzcd}
\Coh[\mathrm{pt}/\mathbb{G}_m^2]_{dg}\arrow[r]\arrow[d]
&
\Coh[\mathbb{P}^1/\mathbb{G}_m^2]_{dg}
\\
\Coh[\mathbb{P}^1/\mathbb{G}_m^2]_{dg}
&
\end{tikzcd}
\right).
\end{aligned}
\]and therefore we have $\Perf\mathcal{W} (\mathcal{X})^{op}\simeq\Coh_{\mathbb{G}_m^2} (\mathbb{P}^1\cup_{0\sim\infty}\mathbb{P}^1)$, which provides the desired categorical equivalence by tensoring with $\Perf\mathcal{W} (X)^{op}\simeq\Coh (Y)_{dg}$ and then quotient by the canonical $\mathbb{Z}^2$-action.
\end{proof}
\begin{proof}[Proof of Corollary \ref{crl:gluing-mirror-symmetry-circ}]
    The proof is completely parallel to the proof of Corollary \ref{crl:gluing-mirror-symmetry}.
\end{proof}
\subsection{Toric bundles.} Note from equivariant coherent-constructible correspondence \cite{FLTZ11} that we have a quasi-equivalence of dg categories \[
    \Sh_{\widetilde{\Lambda}_{\Sigma}} (\mathbb{R}^n)^{cc}\simeq\Coh_{\mathbb{G}_m^n} (X_{\Sigma})
\]when $X_{\Sigma}$ is a smooth toric variety.
Here $\Sh_{\widetilde{\Lambda}_{\Sigma}} (\mathbb{R}^n)^{cc}$ is the dg category of constructible complexes of sheaves with compact supports, and by \cite[Lemma 4.18]{GPS24M} we have \[
    \Sh_{\widetilde{\Lambda}_{\Sigma}} (\mathbb{R}^n)^{cc}\simeq\Sh_{\widetilde{\Lambda}_{\Sigma}} (\mathbb{R}^n)^{\omega},
\]which we know from sheaf-Floer correspondence \ref{prp:inf-type-sheaf-Floer-Rn} to be quasi-equivalent to the wrapped Fukaya category $\Perf\mathcal{W} (T^{\ast}\mathbb{R}^n,\widetilde{\Lambda}_{\Sigma})^{op}$.
Therefore we have a quasi-equivalence of $A_{\infty}$-categories \[
    \Perf\mathcal{W} (T^{\ast}\mathbb{R}^n,\widetilde{\Lambda}_{\Sigma})^{op}\simeq\Coh_{\mathbb{G}_m^n} (X_{\Sigma})
\]whenever $X_{\Sigma}$ is smooth, or equivalently $\Sigma$ is an integral convex polytope.
\begin{proof}[Proof of Corollary \ref{crl:toric-bundle-HMS}]
    From $A$-side we have \[
        \mathcal{W}\left[\mathcal{X}_{\Sigma} (X;\phi_1,\dotsb ,\phi_k)\right]\simeq\mathcal{W} (T^{\ast}\mathbb{R}^k,\widetilde{\Lambda}_{\Sigma})\otimes_{\mathbb{Z}^k}\mathcal{W} (X),
    \]while on the $B$-side we can write $\mathcal{Y}_{\Sigma} (Y;\xi_1,\dotsb ,\xi_k)=\left(X_{\Sigma}\times \prod_{i=1}^k L(Y,\xi_i)^{\times}\right)/\mathbb{G}_m^k$ and so \[
        \Coh\left[\mathcal{Y} (Y;\xi_1,\dotsb ,\xi_k)\right]\simeq\left[\Coh_{\mathbb{G}_m^k} (X_{\Sigma})\otimes\Coh_{\mathbb{G}_m^k} \left(\prod_{i=1}^kL(Y,\xi_i)^{\times}\right)\right]\Big/\mathbb{G}_m^k.
    \]Pairing these two equivalences concludes.
\end{proof}
\subsection{Hirzebruch surfaces.} The last part of the proof is contributed to Homological Mirror Symmetry for Hirzebruch surfaces. Note that the HMS part of Corollary \ref{crl:mirror-to-Hirzebruch} is a direct Corollary of Theorem \ref{thm:HMS-projective-bundle}, and the only non-trivial part is the geometric identification of the constructed space $\rF_n$.
\begin{proof}[Proof of Corollary \ref{crl:mirror-to-Hirzebruch}]
     Consider the affine hypersurface \[
            Z\coloneqq\{(x,y,z)\in\mathbb{C}^3\vert x^{n+1}+yz=0\}\subseteq\mathbb{C}^3
        \]which admits a holomorphic map \[
            \pi\colon Z\to\mathbb{C},\ (x,y,z)\mapsto x
        \]so that $\pi^{-1} (0)=\{(y,z)\vert yz=0\}$ is singular nodal while $\pi^{-1} (t)=\{(y,z)\vert yz =t^{n+1}\}$ is smooth affine, which from its expression (and Riemann-Hurwitz formula) is clearly isomorphic to $\mathbb{C}^{\ast}$. We can also verify that the monodromy action around the singular fibre $x=0$ is the same as the Dehn twist $\phi^n$, and therefore we have $\rF_n\cong_{alg} Z\setminus\{x=0\}$. Note that $\rF_n$ admits a projection $\rF_n\to\mathbb{C}_x^{\ast}\times\mathbb{C}_y^{\ast}$ where the fibre at each pair $(x,y)$ consists exactly of a single point, and therefore we have $\mathcal{X}\cong_{alg}\mathbb{C}_x^{\ast}\times\mathbb{C}_y^{\ast}$. This isomorphism does not change the coordinates of $x$ and $y$, so the superpotential under this algebraic isomorphism is still $x+x^{-1} +y+y^{-1}$.
\end{proof}
\section{Algebraic Preparations\label{sec:algebraic-preparations}}
Before filling in the details, in this section we recall some of the basic definitions of $A_{\infty}$-categories. Throughout the whole paper, except otherwise mentioned, $\field$ would mean a field of characteristic $0$.
\subsection{$A_{\infty}$-categories and dg-categories.} We start with a recollment of the definition of $A_{\infty}$-categories and functors between them.
\begin{dfn}
A
\emph{$\mathbb{Z}$-graded $A_{\infty}$-category $\mathcal{A}$
over $\field$}
is a category enriched in $A_{\infty}$-operads over cochain complexes of $\field$-vector spaces.
Explicitly, it consists of the following:
\begin{itemize}
\item
a class of \emph{objects} $\Ob(\mathcal{A})$.
\item for each pair of objects $L,L'\in\Ob (\mathcal{A})$, a $\mathbb{Z}$-graded $\field$-vector space of morphisms
$$
\hom_{\mathcal{A}}(L,L')\coloneqq \bigoplus_{\bullet\in\mathbb{Z}} \hom_{\mathcal{A}}^{\bullet}(L,L').
$$
\item An \emph{$A_{\infty}$-structure} \[
    \mu_{\mathcal{A}}^{\bullet} =\left(\mu_{\mathcal{A}}^k\colon\hom_{\mathcal{A}}^{\bullet} (L_{k-1},L_k)\otimes\dotsb\otimes\hom_{\mathcal{A}}^{\bullet} (L_0, L_1)\to\hom_{\mathcal{A}}^{\bullet +2-k} (L_0, L_k)\right)_{k\in\mathbb{Z}_{\geq 0}}
\]that satisfies the \emph{Maurer-Cartan equation} $(\mu_{\mathcal{A}}^{\bullet})^2=0$, i.e. \[
    \sum_{m,n}(-1)^{\maltese_n}\mu_{\mathcal{A}}^{d-m+1} (a_d,\dotsb, a_{n-m+1},\mu_{\mathcal{A}}^m (a_{n+m},\dotsb, a_{n+1}),a_n,\dotsb, a_1)=0,
\]where $\maltese_n =\vert a_1\vert+\dotsb+\vert a_n\vert -n$. In particular,
$\mu_{\mathcal{A}}^1$ is a \emph{differential},
i.e.,
$\mu_{\mathcal{A}}^1(\mu_{\mathcal{A}}^1(f))=0$
for any $f\in\hom_{\mathcal{A}}(L,L')$. \[
\mu_{\mathcal{A}}^2\colon\hom_{\mathcal{A}}^{\bullet} (L_1,L_2)\otimes\hom_{\mathcal{A}}^{\bullet} (L_0,L_1)\to\hom_{\mathcal{A}}^{\bullet} (L_0, L_2)
\]serves as the \emph{composition} of two morphisms, where the Maurer-Cartan equation turns into the Liebnitz rule $\mu_{\mathcal{A}}^2 (\mu_{\mathcal{A}}^1 (f),g)\pm\mu_{\mathcal{A}}^2 (f,\mu_{\mathcal{A}}^1 (g))\pm\mu_{\mathcal{A}}^1\mu_{\mathcal{A}}^2 (f,g)=0$ and the homotopy associativity \begin{equation}\label{eq:homotopy-associativity}
\mu^1_{\mathcal{A}}\mu_{\mathcal{A}}^3 (f,g,h)=\pm\mu_{\mathcal{A}}^2 (\mu_{\mathcal{A}}^2 (f,g),h)\pm\mu_{\mathcal{A}}^2 (f,\mu_{\mathcal{A}}^2 (g,h)).
\end{equation}
For simplicity, we also write $(-1)^{-\vert f\vert}\mu_{\mathcal{A}}^2 (f',f)$ as
$$
\comp
:
\hom_{\mathcal{A}}^{\bullet}(L_0,L_1)\underset{\field}{\otimes}\hom_{\mathcal{A}}^{\bullet}(L_1,L_2)
\rightarrow
\hom_{\mathcal{A}}^{\bullet}(L_0,L_2),
\quad
f\otimes f'
\mapsto
f\comp f'
$$
for any objects $L_0$, $L_1$, $L_2\in\Ob(\mathcal{A})$.
\end{itemize}
When $\mu_{\mathcal{A}}^k=0$ for $k\geq 2$, we say $\mathcal{A}$ is a \emph{dg-category}, and the homotopy associativity \eqref{eq:homotopy-associativity} reduces to strict associativity, so that $\mathcal{A}$ is an honest $1$-category. In this situation, we sometimes simply write $d_{\mathcal{A}}\coloneqq\mu_{\mathcal{A}}^1$.
\end{dfn}

\begin{defn}
An \emph{$A_{\infty}$-functor}
$\phi:\mathcal{A}\rightarrow\mathcal{B}$
between two $A_{\infty}$-categories is a collection of maps \[
    \left(\phi^d\colon\hom_{\mathcal{A}} (L_{d-1},L_d)\otimes\dotsb\otimes\hom_{\mathcal{A}} (L_0,L_1)\to\hom_{\mathcal{A}}^{\bullet +1-d} (\phi(L_0),\phi(L_d))\right)_{d\geq 1}
\]and a map of objects $\phi\colon\Ob (\mathcal{A})\to\Ob (\mathcal{B})$ such that the compatibility equation
\begin{align*}
    \sum_r&\sum_{s_1+\dotsb+s_r=d}\mu_{\mathcal{B}}^r (\phi^{s_r}(a_d,\dotsb,a_{d-s_r+1}),\dotsb,\phi^{s_1} (a_{s_1},\dotsb, a_1))\\ 
    &=\sum_{m,n}(-1)^{\maltese_n}\phi^{d-m+1}(a_d,\dotsb, a_{n+m+1},\mu_{\mathcal{A}}^m (a_{n+m},\dotsb, a_{n+1}),a_n,\dotsb, a_1).
\end{align*}
When $\mathcal{A}$ and $\mathcal{B}$ are dg categories, we say $(\phi^d)$ is a \emph{dg-functor} if $\phi^d\equiv 0$ for $d\geq 2$. In other words, we can regard a dg functor as a cochain map
$
\phi:
\hom_{\mathcal{A}}(L,L')
\rightarrow
\hom_{\mathcal{B}}(\phi(L),\phi(L'))
$
for any objects $L,L'\in\Ob(\mathcal{A})$
that preserves the composition,
i.e.,
$
\phi(f\comp f')
=
\phi(f)\comp \phi(f'),
$
and the units,
i.e.,
$
\phi(e_L)
=
e_{\phi(L)}.
$

Compositions of $A_{\infty}$/dg functors can then be defined in an obvious way.
\end{defn}

\begin{defn}
For an $A_{\infty}$-category $\mathcal{A}$,
its \emph{cohomology category}
$H^{\bullet}(\mathcal{A})$
is a category whose objects are the same as $\mathcal{A}$
and the morphism space between two objects $L_0, L_1$ is given by
$$
\Hom_{H^{\bullet}(\mathcal{A})}(L,L')
:=
H^{\bullet}(\hom_{\mathcal{A}}(L,L'),\mu_{\mathcal{A}}^1).
$$
An $A_{\infty}$-functor $(\phi^{\bullet}):\mathcal{A}\rightarrow\mathcal{B}$
induces an ordinary functor
$H^{\bullet}(\phi):H^{\bullet}(\mathcal{A})\rightarrow H^{\bullet}(\mathcal{B})$,
whose mapping on objects is the same as $\phi$ and
the mapping on morphisms is given by $[f]\mapsto [\phi^1(f)]$.

A dg-functor $\phi:\mathcal{A}\rightarrow\mathcal{B}$ is called a
\emph{quasi-embedding}
(resp. \emph{quasi-equivalence}, \emph{quasi-isomorphism})
if the induced functor
$H^{\bullet}(\phi):H^{\bullet}(\mathcal{A})\rightarrow H^{\bullet}(\mathcal{B})$
is an embedding
(resp. equivalence, isomorphism)
of categories.
\end{defn}
As an example, we describe a particular dg category constructed from a $\field$-linear category.
\begin{defn}
Let $\mathcal{C}$ be a $\field$-linear category.
The \emph{dg-category $\Ch(\mathcal{C})_{\dg}$ of cochain complexes in $\mathcal{C}$} is defined as follows:
\begin{itemize}
\item
\emph{objects} are given by cochain complexes in $\mathcal{C}$
$$
M
=
\left(\left(M^{\bullet}\right)_{\bullet\in\mathbb{Z}},
d_M:=\left(d_M^\bullet:M^{\bullet}\rightarrow M^{\bullet+1}\right)_{\bullet\in\mathbb{Z}}\right),
$$
where $M^{\bullet}$ is an object in $\mathcal{C}$ and
$d_M^{\bullet}\in\Hom_{\mathcal{C}}(M^{\bullet},M^{\bullet+1})$
is a morphism in $\mathcal{C}$
satisfying
$d_M^{\bullet}\comp d_M^{\bullet+1}=0$
for each $\bullet\in\mathbb{Z}$.

\item
a \emph{dg-morphism space}
is a cochain complex of $\field$-modules
whose degree $k$-component is given by
$$
\hom_{\Ch(\mathcal{C})_{\dg}}^k
\left(
M,M'
\right)
:=
\prod_{\bullet\in\mathbb{Z}}
\Hom_{\mathcal{C}}\left(M^{\bullet},M'^{\bullet+k}\right),
$$
and the differential is given by
$$
df:=d_M\comp f - (-1)^k f\comp d_{M'}
$$
for $f\in \hom_{\Ch(\mathcal{C})_{\dg}}^k(M,M')$.
More precisely,
the degree $\bullet$-component of $df$
is given by
$$
(df)^{\bullet}
:=
d_M^{\bullet}\circ f^{\bullet+1}
-
(-1)^k
f^{\bullet}\circ d_{M'}^{\bullet+k}.
$$
\item
the composition and the units are defined in an obvious way.

\end{itemize}

We denote by
$\Ch^b(\mathcal{C})_{\dg}$
the \emph{dg-category of bounded cochain complexes in $\mathcal{C}$},
whose objects are cochain complexes
$M=\left(\left(M^{\bullet}\right)_{\bullet\in\mathbb{Z}},d_M\right)$
in $\mathcal{C}$
such that
$M^{\bullet}=0$ except for finitely many $\bullet\in\mathbb{Z}$.
The morphism spaces are defined in the same way as in $\Ch(\mathcal{C})_{\dg}$.
\end{defn}
\begin{rmk}\label{rmk:EmbeddingOfCochainComplexes}
A functor
$\mathcal{F}:\mathcal{C}\hookrightarrow\mathcal{D}$
between $\field$-linear categories
naturally induces a dg-functor between dg-categories
$
\mathcal{F}_{\dg}:
\Ch(\mathcal{C})_{\dg}
\hookrightarrow
\Ch(\mathcal{D})_{\dg}
$
and
$
\Ch^b(\mathcal{C})_{\dg}
\hookrightarrow
\Ch^b(\mathcal{D})_{\dg},
$
in an obvious way.
If
$\mathcal{F}$ is an embedding
(resp., equivalence, isomorphism),
then
$\mathcal{F}_{\dg}$
is a quasi-embedding
(resp., quasi-equivalence, quasi-isomorphism).
\end{rmk}
\subsection[modules and perfect modules]{$A_{\infty}$-modules and perfect modules.} For an $A_{\infty}$-category $\mathcal{A}$, recall that an \emph{$\mathcal{A}$-module} is an $A_{\infty}$-functor $\mathcal{M}\colon\mathcal{A}\to\Ch_{\field}$ from $\mathcal{A}$ to the dg category of chain complexes over $\field$.
\begin{dfn}\label{dfn:A-inf-module}
     The $A_{\infty}$-module structure on $\mathcal{M}\in\Fun (\mathcal{A},\Ch_{\field})\eqqcolon\Mod (\mathcal{A})$ is a collection of $\field$-linear maps
    \begin{align*}
    \mu^{d\vert 1}&\colon\hom_{\mathcal{A}} (L_{d-1},L_d)\otimes\dotsb\otimes\hom_{\mathcal{A}} (L_0, L_1)\otimes\mathcal{M} (L_0)\to\mathcal{M} (L_d) [1-d]
    \end{align*}
    such that for each $d\geq 0$, \begin{equation}\label{eq:A-inf-bimodule-equation}\begin{aligned}
        &\sum (-1)^{\maltese_{d-i}}\mu_{\mathcal{M}}^{d-i\vert 1}(v_d,\dotsb,v_{i+1},\mu_{\mathcal{M}}^{i\vert 1} (v_i,\dotsb, v_1,\mathbf{b}))\\ 
        &+\sum (-1)^{\maltese_{d-k-i}}\mu_{\mathcal{M}}^{d-i+1\vert 1}(v_d,\dotsb, v_{k+i+1},\mu_{\mathcal{A}}^i (v_{k+i},\dotsb,v_{k+1}),v_k,\dotsb, v_1,\mathbf{b})=0.
    \end{aligned}\end{equation}
    Let $\mathcal{M}$ and $\mathcal{N}$ be two $\mathcal{A}$-modules.
    A \emph{premorphism} $s\colon\mathcal{M}\to\mathcal{N}$ of degree $k$ is a family of morphisms
    $s^d\colon\hom_{\mathcal A}(L_{d-2},L_{d-1})\otimes\dotsb\otimes\hom_{\mathcal A}(L_0,L_1)\otimes\mathcal M(L_0)\longrightarrow\mathcal N(L_{d-1})[k-d+1]$
    of chain complexes associated to tuples $(L_0,\dotsb, L_{d-1})$.\par
    The set of premorphisms from $\mathcal{M}$ to $\mathcal{N}$ form a chain complex over $\Bbbk$, and the category of $\mathcal{A}$-modules $\Mod (\mathcal{A})$ is a dg category whose differentials are given explicitly in \cite[page 20]{Seidel08}. Pre-morphisms $s$ such that $\mu^1_{\Mod_{\mathcal{A}}} (s)=0$ are called \emph{morphisms} of $A_{\infty}$-modules. Two $\mathcal{A}$-modules $\mathcal{M}$ and $\mathcal{N}$ are \emph{isomorphic} if there exists a quasi-isomorphism $f\colon\mathcal{M}\to\mathcal{N}$, i.e. there is a morphism $g\colon\mathcal{N}\to\mathcal{M}$ so that both $f\circ g$ and $g\circ f$ are identities after taking cohomologies.
\end{dfn}
For each object $X\in\Obj (\mathcal{A})$, the assignment $X\mapsto\hom_{\mathcal{A}} (X,-)$ defines an $A_{\infty}$-functor $\mathcal{A}^{op}\to\Mod_{\field} (\mathcal{A})$ in the way as described in \cite[section (1l)]{Seidel08}. This is the $A_{\infty}$-version of the \emph{Yoneda embedding}, denoted by $\yo_{-}\colon\mathcal{A}\to\Mod (\mathcal{A})$.
\begin{dfn}
    An $A_{\infty}$-module $\mathcal{M}\in\Mod_{\field} (\mathcal{A})$ is called \emph{representable} in $\mathcal{A}$ if there exists an object $X\in\Obj (\mathcal{A})$ and an isomorphism of $\mathcal{A}$-modules $\mathcal{M}\xrightarrow{\simeq}\yo_X$.\par
    We define the category of \emph{perfect $\mathcal{A}$-modules} $\Perf\mathcal{A}\subseteq\Mod (\mathcal{A})$ to be the smallest pre-triangulated dg subcategory containing images of the Yoneda functors.
    Equivalently, $\Perf\mathcal{A}$ is the idempotent completion of $\mathcal{A}^{op}$.
\end{dfn}
\subsection{Group actions on $A_{\infty}$-categories and orbit categories.} In this section we define orbit categories associated to a discrete group $G$ acting on an $A_{\infty}$-category $\mathcal{A}$. This would coincide with a special case of \cite{FKQ24} when $G=\mathbb{Z}$. We assume in this section that $G$ is either a finite group or a free abelian group of finite rank.
\begin{dfn}[Group Action on $A_{\infty}$-Categories]
    An \emph{action} of $G$ on $\mathcal{A}$ assigns to each $g\in G$ a quasi-equivalence $\{g^i\}_{i\geq 0}\colon\mathcal{A}\to\mathcal{A}$ of $A_{\infty}$-categories, where \[
    g^i\colon\hom_{\mathcal{A}} (L_{i-1}, L_i)\otimes\dotsb\otimes\hom_{\mathcal{A}} (L_0, L_1)\to\hom_{\mathcal{A}} (g(L_0),g(L_i))[i-1]
    \]forms an $A_{\infty}$-endofunctor of $\mathcal{A}$. We say an $G$-action is \emph{strict} if $g^k =0$ for $k\geq 2$.
\end{dfn}
\begin{dfn}\label{defn:Ainf-orbit-category}
    Given an $A_{\infty}$-category $\mathcal{A}$ with a strict $G$-action, the \emph{orbit category} of $\mathcal{A}$ under the action $G$, written as $\mathcal{A}/G$, is a category with objects $\Obj\mathcal{A}$ and morphism spaces\[
    \hom_{\mathcal{A}/G} (L,K)\coloneqq\bigoplus_{g\in G}\hom_{\mathcal{A}} (g(L),(K)).
    \]
    For given tuples of objects $L_0, L_1,\dotsb ,L_d$ and chosen morphisms $a_i\in\hom_{\mathcal{A}/G} (L_i, L_{i+1})$, the $A_{\infty}$-structure on $\mathcal{A}/G$ \begin{align*}
    \mu_{\mathcal{A}/G}^d &\colon\hom_{\mathcal{A}} ((g_1\dotsb g_d)(L_{d-1}),(g_1\dotsb g_{d-1})(L_d))\otimes\dotsb\otimes\hom_{\mathcal{A}} (g_1 (L_0),L_1)\\ 
    &\to\hom_{\mathcal{A}} ((g_1\dotsb g_d)(L_0),L_d)[2-d]
    \end{align*}
    is given explicitly by 
    \begin{align*}
        \mu_{\mathcal{A}/G}^d (a_{d-1}\otimes\dotsb\otimes a_0)&=\mu_{\mathcal{A}}^d ((g_1\dotsb g_d) (a_{d-1})\otimes\dotsb\otimes g_1 (a_1)\otimes a_0),
    \end{align*}
    where $a_i\in\hom_{\mathcal{A}} (L_i,g_i (L_{i+1}))$ is $G$-homogeneous, and extended $d$-linearly.
\end{dfn}
The corresponding $dg$- and ordinary category version is also analogous, and we have
\begin{lmm}
    $\mathcal{A}/G$ is an $A_{\infty}$-category.
\end{lmm}
\begin{proof}
    We verify directly that\begin{align*}
        \sum_{k,l,j}&\mu^{j+k}_{\mathcal{A}}(g_1\dotsb g_j h_1\dotsb h_l g_{j+1}\dotsb g_{j+k-1} (a_{j+k+l})\otimes\dotsb\otimes g_1\dotsb g_j h_1\dotsb h_l (p_{j+l+1})\\ 
        &\otimes g_1\dotsb g_j\mu^l_{\mathcal{A}}(h_1\dotsb h_{l-1} (a_{j+l})\otimes\dotsb\otimes a_{j+1})\otimes g_1\dotsb g_{j-1} (a_j)\otimes\dotsb\otimes a_1)\\ 
        =&\sum_{k,l,j}\mu_{\mathcal{A}}^{j+k}(g_1\dotsb g_j h_1\dotsb h_l g_{j+1}\dotsb g_{j+k-1} (a_{j+k+l})\otimes\dotsb\otimes g_1\dotsb g_j h_1\dotsb h_l (p_{j+l+1})\\ 
        &\otimes\mu^l_{\mathcal{A}}\left(\mathbf{g}(h_1\dotsb h_l(a_{j+l})\otimes\dotsb\otimes\mathbf{g} h_1 (a_{j+1}))\right)\otimes g_1\dotsb g_{j-1} (a_j)\otimes\dotsb\otimes a_1)=0
    \end{align*}
    where $\mathbf{g}=g_1\dotsb g_j$, using the Maurer-Cartan equation for $\mu_{\mathcal{A}}^{\bullet}$.
\end{proof}
Note that when $\mathcal{A}$ is a dg category, i.e. $\mu_{\mathcal{A}}^d=0$ for $d\geq 3$ and $G\curvearrowright\mathcal{A}$ is a strict action, then the orbit category $\mathcal{A}/\mathbb{Z}$ is exactly the dg enhanced orbit category in the sense of Fan-Keller-Qiu \cite{FKQ24}. In this case the functor induced from $1\in\mathbb{Z}$ is naturally invertible, and therefore the natural transform $q\colon Q_{\mathbb{N}}F\to Q_{\mathbb{N}}$ in \cite{FKQ24} is an equivalence of functors, so we have nothing to localize.\par
As in \cite{FKQ24}, there is a functor $\mathcal{A}\to\mathcal{A}/G$ which is identity on the set of objects, and on morphisms given by the embedding $\hom_{\mathcal{A}} (X,Y)\hookrightarrow\hom_{\mathcal{A}/G} (X,Y)=\displaystyle\bigoplus_{g\in G}\hom_{\mathcal{A}} (X,g^0(Y))$ into the component where $g=\id$. It's not hard to check that this defines a strict $A_{\infty}$-functor $Q_G\colon\mathcal{A}\to\mathcal{A}/G$.
\begin{exmp}
Let $\Lambda$ be a group,
$\mathcal{C}$ a $\field$-linear category,
and
$\phi:\Lambda\rightarrow\Aut_{\str}(\mathcal{C})$
a strict $\Lambda$-action on $\mathcal{C}$.
It induces a strict $\Lambda$-action
$\phi:\Lambda\rightarrow\Aut_{\str}(\Ch(\mathcal{C})_{\dg})$
on the dg-category $\Ch(\mathcal{C})_{\dg}$
in a straightforward way.
\end{exmp}
\subsection{Tensor products of dg-categories.} Given $A_{\infty}$-categories $\mathcal{A}$ and $\mathcal{B}$, their tensor product $\mathcal{A}\otimes\mathcal{B}$ was recently written down by Ornaghi \cite{Ornaghi25}, which is unfortunately very complicated. But when $\mathcal{A}$ and $\mathcal{B}$ are dg categories, we have a simpler formula that can be written down, and is a special situation of Ornaghi's construction, which we will recall in the following definition.
\begin{defn}
Let $\mathcal{A}$, $\mathcal{B}$ be dg-categories,
$\Lambda$ an abelian group,
and
$\phi:\Lambda\rightarrow\Aut_{\str}(\mathcal{A})$,
$\psi:\Lambda\rightarrow\Aut_{\str}(\mathcal{B})$
strict $\Lambda$-actions on $\mathcal{A}$, $\mathcal{B}$,
respectively.
The \emph{twisted tensor product}
$
\mathcal{A}\underset{\Lambda}{\otimes}\mathcal{B}
$
of $\mathcal{A}$ and $\mathcal{B}$
with respect to the actions $\phi$ and $\psi$ of $\Lambda$
is a dg-category
consisting of the following:

\begin{itemize}
\item
$
\Ob\left(\mathcal{A}\underset{\Lambda}{\otimes}\mathcal{B}\right)
:=
\Ob(\mathcal{A})\times \Ob(\mathcal{B}).
$

\item
a dg-morphism space
$$
\hom_{\mathcal{A}\underset{\Lambda}{\otimes}\mathcal{B}}
\left((L,R),(L',R')\right)
:=
\displaystyle\bigoplus_{\lambda\in\Lambda}
\left(
\hom_{\mathcal{A}}(\phi_{\lambda}(L),L')
\underset{\field}{\otimes}
\hom_{\mathcal{B}}(\psi_{\lambda}(R),R')
\right)
$$
as a cochain complex of $\field$-vector spaces.
In particular,
the differential is given by
$$
d_{\mathcal{A}\underset{\Lambda}{\otimes}\mathcal{B}}
(f\otimes g)
=
d_{\mathcal{A}}f \otimes g
+
(-1)^{|f|}
f \otimes d_{\mathcal{B}}g.
$$

\item
composition of
$\Lambda$-homogeneous
dg-morphisms
$$
f\otimes g
\in
\hom_{\mathcal{A}}(\phi_{\lambda}(L),L')
\underset{\field}{\otimes}
\hom_{\mathcal{B}}(\psi_{\lambda}(R),R')
\subseteq
\hom_{\mathcal{A}\underset{\Lambda}{\otimes}\mathcal{B}}
\left((L,R),(L',R')\right)
$$
and
$$
f'\otimes g'
\in
\hom_{\mathcal{A}}(\phi_{\lambda'}(L'),L'')
\underset{\field}{\otimes}
\hom_{\mathcal{B}}(\psi_{\lambda'}(R'),R'')
\subseteq
\hom_{\mathcal{A}\underset{\Lambda}{\otimes}\mathcal{B}}
\left((L',R'),(L'',R'')\right)
$$
is defined by
$$
(f\otimes g)\comp(f'\otimes g')
:=
(-1)^{|g||f'|}
(\phi_{\lambda'}(f)\comp f')
\otimes
(\psi_{\lambda'}(g)\comp g')
$$
$$
\in
\hom_{\mathcal{A}}(\phi_{\lambda\lambda'}(L),L'')
\underset{\field}{\otimes}
\hom_{\mathcal{B}}(\psi_{\lambda\lambda'}(R),R'')
\subseteq
\hom_{\mathcal{A}\underset{\Lambda}{\otimes}\mathcal{B}}
\left((L,R),(L'',R'')\right),
$$
and it is extended to whole dg-morphisms by $\field$-linearity.
\end{itemize}
When $\Lambda=\{1\}$ and $\phi=\id$, we simply write $\mathcal{A}\otimes_1\mathcal{B} =\mathcal{A}\otimes\mathcal{B}$ and call it the \emph{tensor product} of $\mathcal{A}$ and $\mathcal{B}$.
\end{defn}
\begin{defn}
Let $\mathcal{A}$, $\mathcal{B}$, $\mathcal{A}'$, $\mathcal{B}'$ be dg-categories,
and
$\phi:\mathcal{A}\rightarrow\mathcal{A}'$,
$\psi:\mathcal{B}\rightarrow\mathcal{B}'$
be
dg-functors.
The \emph{tensor product}
$\phi\otimes\psi
:
\mathcal{A}\otimes\mathcal{B}
\rightarrow
\mathcal{A}'\otimes\mathcal{B}'
$
of $\phi$ and $\psi$
is a dg-functor
consisting of
a mapping of objects
$
(L,R)\mapsto (\phi(L),\psi(R))
$
and
a mapping of morphisms
$
f\otimes g \mapsto \phi(f)\otimes\psi(g).
$
\end{defn}
The following proposition is straightforward.
\begin{prp}\label{prp:twisted-tensor-equal-quotient}
Let $\mathcal{A}$, $\mathcal{B}$ be dg-categories,
and
$\phi:\Lambda\rightarrow\Aut_{\str}(\mathcal{A})$,
$\psi:\Lambda\rightarrow\Aut_{\str}(\mathcal{B})$
strict $\Lambda$-actions on $\mathcal{A}$, $\mathcal{B}$,
respectively.
Consider a strict $\Lambda$-action
$\phi\otimes \psi:\Lambda\rightarrow\Aut_{\str}(\mathcal{A}\otimes\mathcal{B})$
on $\mathcal{A}\otimes\mathcal{B}$
defined by
$(\phi\otimes\psi)_{\lambda}:=\phi_{\lambda}\otimes \psi_{\lambda}
:
\mathcal{A}\otimes\mathcal{B}\rightarrow\mathcal{A}\otimes\mathcal{B}$
for each $\lambda\in\Lambda$.

Then the twisted tensor product of $\mathcal{A}$ and $\mathcal{B}$
with respect to $\phi$ and $\psi$
is isomorphic to
the orbit category of $\mathcal{A}\otimes\mathcal{B}$
with respect to $\phi\otimes\psi$:
$$
\mathcal{A}\underset{\Lambda}{\otimes}\mathcal{B}
\cong
(\mathcal{A}\otimes\mathcal{B})/\Lambda.
$$
\end{prp}

\begin{proof}
The class of objects for both categories is identically $\Ob(\mathcal{A})\times \Ob(\mathcal{B})$. 
For any two objects $(L,R)$ and $(L',R')$, the morphism space in the twisted tensor product is:
$$
\hom_{\mathcal{A}\underset{\Lambda}{\otimes}\mathcal{B}}((L,R),(L',R'))
=
\bigoplus_{\lambda\in\Lambda}
\left(
\hom_{\mathcal{A}}(\phi_{\lambda}(L),L')
\underset{\field}{\otimes}
\hom_{\mathcal{B}}(\psi_{\lambda}(R),R')
\right).
$$
On the other hand, by definition of the orbit category, the morphism space in $(\mathcal{A}\otimes\mathcal{B})/\Lambda$ is:
\begin{align*}
\hom_{(\mathcal{A}\otimes\mathcal{B})/\Lambda}((L,R),(L',R'))
&=
\bigoplus_{\lambda\in\Lambda}
\hom_{\mathcal{A}\otimes\mathcal{B}}((\phi\otimes\psi)_\lambda(L,R), (L',R'))
\\
&=
\bigoplus_{\lambda\in\Lambda}
\hom_{\mathcal{A}\otimes\mathcal{B}}((\phi_\lambda(L),\psi_\lambda(R)), (L',R'))
\\
&=
\bigoplus_{\lambda\in\Lambda}
\left(
\hom_{\mathcal{A}}(\phi_{\lambda}(L),L')
\underset{\field}{\otimes}
\hom_{\mathcal{B}}(\psi_{\lambda}(R),R')
\right).
\end{align*}

These spaces match exactly as cochain complexes of vector spaces. 
Furthermore, the composition in the orbit category $(\mathcal{A}\otimes\mathcal{B})/\Lambda$ 
coincides with the composition rule defined in the twisted tensor product. Thus, the identity map on objects and morphisms gives an isomorphism of dg-categories.
\end{proof}
\section{Fukaya Categories of Discrete Quotients\label{sec:Fukaya-category-discrete-quotient}}
Suppose $X$ is a Weinstein sector and $G\curvearrowright X$ be a discrete group acting exact symplectically, discrete continuously and freely on $X$. We further assume that the action is \emph{cylindrical}, i.e. in a cylindrical end $\partial_{\infty}X\times [0,+\infty)\hookrightarrow X$ of $X$, the $G$-action on $X$ preserves the cylindrical coordinate $[0,+\infty)$ and is a contactomorphism on each contact boundary at infinity $\partial_{\infty} X$ of $X$.\par
By assumption the quotient space $X/G$ is still a Weinstein sector on $X$ with ideal boundary at infinity $\partial_{\infty} X/G$. In this section, we are going to compute the wrapped Fukaya category $\mathcal{W} (X/G)$ of $X/G$ in terms of the wrapped Fukaya category $\mathcal{W} (X)$. Our main theorem for this section is:
\begin{thm}[Wrapped Fukaya Categories of Quotients]\label{thm:wrapped-Fukaya-quotient}
    There is an equivalence of $A_{\infty}$-categories
    \[
        \mathcal{W} (X/G)\simeq\mathcal{W} (X)/G \simeq\mathcal{W}(X)_G
    \]
    after twisted completion, where $\mathcal{W} (X)/G$ is the $A_{\infty}$-variant of dg orbit categories developed by Fan-Keller-Qiu \cite{FKQ24} and $\mathcal{W}(X)_G$ is the homotopy quotient of $\mathcal{W}(X)$ by the $G$-action, whose equivalence with $\mathcal{W}(X)/G$ was proved by Christ \cite{Christ26}.
\end{thm}
\subsection{Weinstein sectors with group action.} In this section we recall the definition of Weinstein sectors and construct Weinstein structures on Weinstein sectors that is invariant under a finitely generated discrete group action. Recall that
\begin{dfn}\label{defn:Liouville-sector-inf}
    A \emph{Liouville sector} is a triple $(X,\omega,Z)$ where $(X,\omega)$ is a symplectic manifold with boundary and $Z$ is a complete vector field on $X$ such that 
    \begin{enumerate}[\indent (i)]
        \item (Infinite-type) There exists an exhaustion $X=\displaystyle\bigcup_i X_i$ where $X_0\subseteq X_1\subseteq\dotsb$ is a sequence of compact symplectic manifolds with corners, such that $Z$ is outward-pointing along $\partial X_i\setminus\partial X$ for all $i$.
        \item We have $\mathcal{L}_Z\omega =\omega$, i.e. $Z$ scales $\omega$ exponentially.
        \item There exists a smooth function $I\colon\Nbd (\partial X)\to\mathbb{R}$ defined in an open neighbourhood of $\partial X$, a constant $0<\alpha\leq 1$ and a compact subset $K_i\subseteq\partial X\cap (X_i\setminus X_{i-1})$ for each $i$ such that $ZI=\alpha I$ away from $\displaystyle\bigcup_i K_i$ and the hamiltonian vector field $X_I$ is outward-pointing along $\partial X$.
    \end{enumerate}
    We say $X$ is \emph{Weinstein} if there exists a smooth function $\phi\colon X\to\mathbb{R}$ such that $Z$ is gradient-like for $\phi$, $\phi(\partial X_i)=\{r_i\}$ is a sequence of constants diverging to $\infty$, and the restriction of $\phi$ to $I^{-1}(0)$ defines a Weinstein structure for some chosen smooth function $I$.
\end{dfn}
This definition differs from the original one of Ganatra-Pardon-Shende \cite{GPS20} in that we allow Liouville sectors of infinite type, which when $\partial X=\emptyset$ coincides with \cite{CE12}.
\begin{dfn}\label{defn:G-invariant-Weinstein-sector}
    Let $G$ be a Lie group and $(X,Z_X,\omega_X)$ be a Liouville sector. A $G$-action $G\curvearrowright X$ is \emph{exact symplectic} if for all $g\in G$, there exists $f_g\in C^{\infty} (X;\mathbb{R})$, depending smoothly on $G$, such that $g^{\ast}\lambda_X =\lambda_X +df_g$ where the zeroes of $\lambda_X +df_g$ lie in a compact subset of $X$. It's \emph{compactly supported} if for all $g\in G$, there exists a compact subset $K_g\subseteq X$ such that $g$ is identity outside $K_g$.\par
    An exact symplectic $G$-action on $X$ is \emph{Liouville} if $g^{\ast}\lambda_X =\lambda_X$ for all $g\in G$, where $\lambda_X$ is the symplectic dual of $Z_X$. If $X$ is furthermore Weinstein with Weinstein structure $\phi_X$, then we say $G$-action is \emph{Weinstein} if furthermore $g^{\ast}\phi_X =\phi_X$ for all $g\in G$.
\end{dfn}
\begin{lmm}
    Suppose that $X$ is a Weinstein sector with an exact symplectic $G$-action with compact support and $G$ is a finite group, then there exists a Weinstein structure on $X$ such that $G$ is Weinstein.
\end{lmm}
\begin{proof}
    \textbf{Step 1.} In each cases we can assign to each $g\in G$ a smooth function $f_g\in C^{\infty} (X,\mathbb{R})$ so that $g^{\ast}\lambda_X =\lambda_X +df_g$. We look for Liouville $1$-form $\lambda\coloneqq\lambda_X +df$, for some $f\in C^{\infty}(X,\mathbb{R})$ so that $g^{\ast}\lambda =\lambda$ for all $g\in G$. Plugging in $g$ gives \[
        g^{\ast}\lambda -\lambda =d(f_g +f\circ g-f),
    \]and therefore $f$ should satisfy the identity \begin{equation}\label{eq:invariance-condition}
        f_g +f\circ g-f=0
    \end{equation}for all $g\in G$. Note that for $h\in G$, we have the composition \[
        h^{\ast}g^{\ast}\lambda_X =h^{\ast} (\lambda_X +df_g)=\lambda_X +d(f_h +f_g\circ h)=\lambda_X +df_{gh},
    \]the functions $\{f_g\}_{g\in G}$ must satisfy $f_{gh}=f_h +f_g\circ h$ up to a constant shift, which we can arrange so that the strict identity holds in both situations.\par
    \textbf{Step 2.} Suppose that $G$ is finite, then we consider the function \[
        f=\frac{1}{\vert G\vert}\sum_{g\in G} f_g.
    \]It's not difficult to check that equation \eqref{eq:invariance-condition} holds for any $g\in G$, and therefore the $1$-form $\lambda_G\coloneqq\lambda_X +df$ satisfies $d\lambda_G =\omega_X$ and $g^{\ast}\lambda_G =\lambda_G$. The dual vector field reads \[
        Z_G =Z_X +X_f =Z_X+\frac{1}{\vert G\vert}\sum_{g\in G} X_{f_g},
    \]
    which is complete as $X_{f_g}$ are zero away from a given compcat subset of $X$.
\end{proof}
Assume that $(X,Z,\omega)$ is of \emph{finite type}, i.e. the zeroes of $Z$ is compact, then in \cite{GPS20} Ganatra-Pardon-Shende constructed an $A_{\infty}$-category $\mathcal{W} (X,\lambda ,Z)$ of $X$, called the \emph{wrapped Fukaya category} of $X$. In \cite{GPS24S}, they further showed that 
\begin{thm}[Viterbo Restriction]\label{thm:Viterbo-restriction}
    Let $X_1, X_2$ be Weinstein sectors of finite type such that $\bar{X}_1$, the compact symplectic submanifold of $X_1$ with corners containing all critical points of $Z_1$ such that $Z_1$ is outward-pointing along $\partial\bar{X}_1\setminus (\bar{X}_1\cap\partial X_1)$, admits a symplectic embedding into $X_2$ mapping $\partial X_1\cap\bar{X}_1$ into $\partial X_2$, so that the complement $\overline{\bar{X}_2\setminus\bar{X}_1}$ is a Weinstein cobordism, then there exists an essential surjective $A_{\infty}$-functor $\mathcal{W} (X_2)\to\mathcal{W} (X_1)$, called the \emph{Viterbo restriction}.
\end{thm}
Using this and the exhaustion $X_0\subseteq X_1\subseteq\dotsb$ for a given general Weinstein sector, we can formulate the following definition of wrapped Fukaya categories for Weinstein sectors:
\begin{dfn}\label{defn:wrapped-Fukaya-cat-inf-sectors}
    Let $(X,Z,\omega,\phi)$ be a Weinstein sector, then its \emph{big wrapped Fukaya category} $\mathcal{W} (X)$ is defined as the homotopy limit \[
        \mathcal{W}^{big} (X)\simeq\varprojlim_{i\to\infty}\mathcal{W} (\widehat{X}_i),
    \]where $\widehat{X}_i$ is the completion of $X_i$ into a Weinstein sector of finite type. Let $L_{ij}\subseteq X_i$ be a cocore of one of the handles of the Weinstein domain $X_i$, which defines an object in $\mathcal{W}^{big} (X)$ since it's the Viterbo restriction image of any Lagrangian $L_{ij}\subseteq X_k$ in any larger Weinstein domains.\par
    We define the \emph{wrapped Fukaya category} $\mathcal{W} (X)$ of $X$ to be the full $A_{\infty}$-subcategory of $\mathcal{W}^{big} (X)$ generated by these objects $L_{ij}$.
\end{dfn}
\begin{rmk}
    This definition is far from satisfactory. One should expect in general that $\mathcal{W} (X)$ is the $A_{\infty}$-full subcategory of compact objects of the large Fukaya category $\mathcal{W}^{big} (X)$. However, this does not automatically algebraically hold, and one needs some geometric inputs to prove this fact. \par
    In the work of \cite{KW25}, the wrapped Fukaya category was defined as $\mathcal{W}^{big} (X)$ in this paper, and it's also interesting to ask why the mirror behaviours in these two situations are drastically different.
\end{rmk}
From this definition, generation results for Fukaya categories of Weinstein sectors are clear:
\begin{lmm}[Generation results for Weinstein Sectors \cite{GPS24S}]\label{lmm:generation-Weinstein-sector}
    $\mathcal{W} (X)$ is generated by homological cocores of the Weinstein structure $\phi$ on $X$.
\end{lmm}
Given a Stein manifold $X$ and a holomorphic map $\phi\colon X\to\mathbb{C}$ whose critical values lie in a compact subset, the pair $(X,\phi)$ is called a \emph{Landau-Ginzburg model}. As in \cite{GPS20}, a Landau-Ginzburg model is equivalent to a Weinstein sector via the following construction:\par
Pick $M\gg 0$ so that all critical values of $\phi$ lie inside $\phi^{-1} (\mathbb{D}_0 (M))$, then the construction \cite[Definition 2.14]{GPS20} provides a Liouville sector $X_{\phi}$ with sectorial boundary $\partial X_{\phi}\cong_{\diff} \phi^{-1} (M)\times\mathbb{R}$. We will not distinguish between Liouville sectors and Landau-Ginzburg models in the sequel.
\begin{exmp}[Mirror to $\Perf_{\mathbb{G}_m} (\mathbb{P}^1)$]\label{exmp:mirror-to-Coh-G-P1}
    Consider the pair $(\mathbb{C}, e^z+e^{-z})$ which is the $\mathbb{Z}$-fold cover of the sector $(\mathbb{C}^{\ast},z+1/z)$ that is mirror to $\mathbb{P}^1$. From \cite{FLTZ12} we know that $\mathcal{W} (\mathbb{C},e^z+e^{-z})$ is quasi-equivalent to the dg category of $\mathbb{G}_m$-equivariant coherent sheaves on $\mathbb{P}^1$: $\Coh_{\mathbb{G}_m} (\mathbb{P}^1)$. Let's compute this category explicitly.\par
    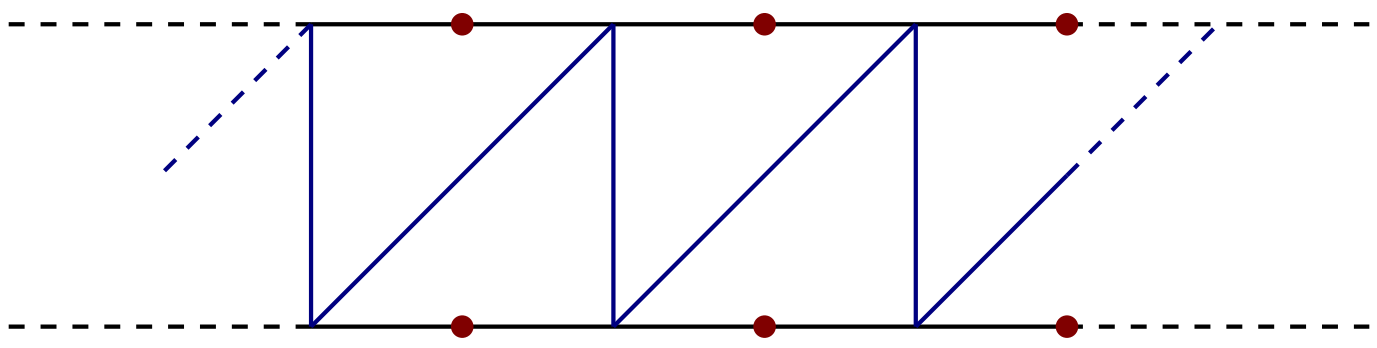
\begin{figure}[ht]
        \centering
        \begin{tikzpicture}
            \draw [thick] (2,2) -- (7,2);
            \draw [thick,dashed] (0,2) -- (2,2);
            \draw [thick, dashed] (7,2) -- (9,2);
            \draw [thick,dashed] (0,0) -- (2,0);
            \draw [thick,dashed] (7,0) -- (9,0);
            \draw [thick] (2,0) -- (7,0);
            \foreach \x in {3,5,7} {
                \node [circle,fill=red!50!black, inner sep=1.5pt] at (\x, 2) {};
                \node [circle,fill=red!50!black,inner sep=1.5pt] at (\x,0) {};
            }
            \draw [thick,blue!50!black] (4,0) -- (4,2);
            \draw [thick,blue!50!black] (4,0) -- (6,2);
            \draw [thick,blue!50!black] (6,2) -- (6,0);
            \draw [thick,blue!50!black] (4,2) -- (2,0);
            \draw [thick,blue!50!black] (6,0) -- (7,1);
            \draw [thick, blue!50!black,dashed] (7,1) -- (8,2);
            \draw [thick, blue!50!black] (2,0) -- (2,2);
            \draw [thick, blue!50!black, dashed] (2,2) -- (1,1);
        \end{tikzpicture}
        \caption{Generating Lagrangians of $\mathcal{W} (\mathbb{C}, e^z+e^{-z})$.}\label{fig:generating-Lag-equiv-P1}
    \end{figure}
    The space $(\mathbb{C},e^z+e^{-z})$ can be pictured as in figure \ref{fig:generating-Lag-equiv-P1}, where the red dots are fibres of $e^z+e^{-z}$ near infinity, and blue lines are Lagrangian submanifolds that are projected diffeomorphically onto $\mathbb{C}^{\ast}$ via the exponential map $z\mapsto e^z$. We claim that these Lagrangian submanifolds generate the wrapped Fukaya category $\mathcal{W} (\mathbb{C},e^z+e^{-z})$.\par
    To prove the claim, first consider the trivial Weinstein sector $X_0\coloneqq(\mathbb{C},z^2) \simeq T^{\ast} [0,1]$, then the Weinstein sector $(\mathbb{C},e^z+e^{-z})$ can be regarded as attaching $\mathbb{Z}$-copies of $1$-handles to the disjoint union \[
        \bigsqcup_{\mathbb{Z}} X_0
    \]in the way depicted in figure \ref{fig:1-handle-attaching-Z}. It's then clear that the wrapped Fukaya category $\mathcal{W} (\mathbb{C},e^z+e^{-z})$ is generated by the union of blue Lagrangians in figure \ref{fig:1-handle-attaching-Z}, as they are generalized cocores to the sector $X_0$ and to the $1$-handles. It's not difficult to check that these Lagrangians are exactly the blue Lagrangians shown in figure \ref{fig:generating-Lag-equiv-P1}.\par
    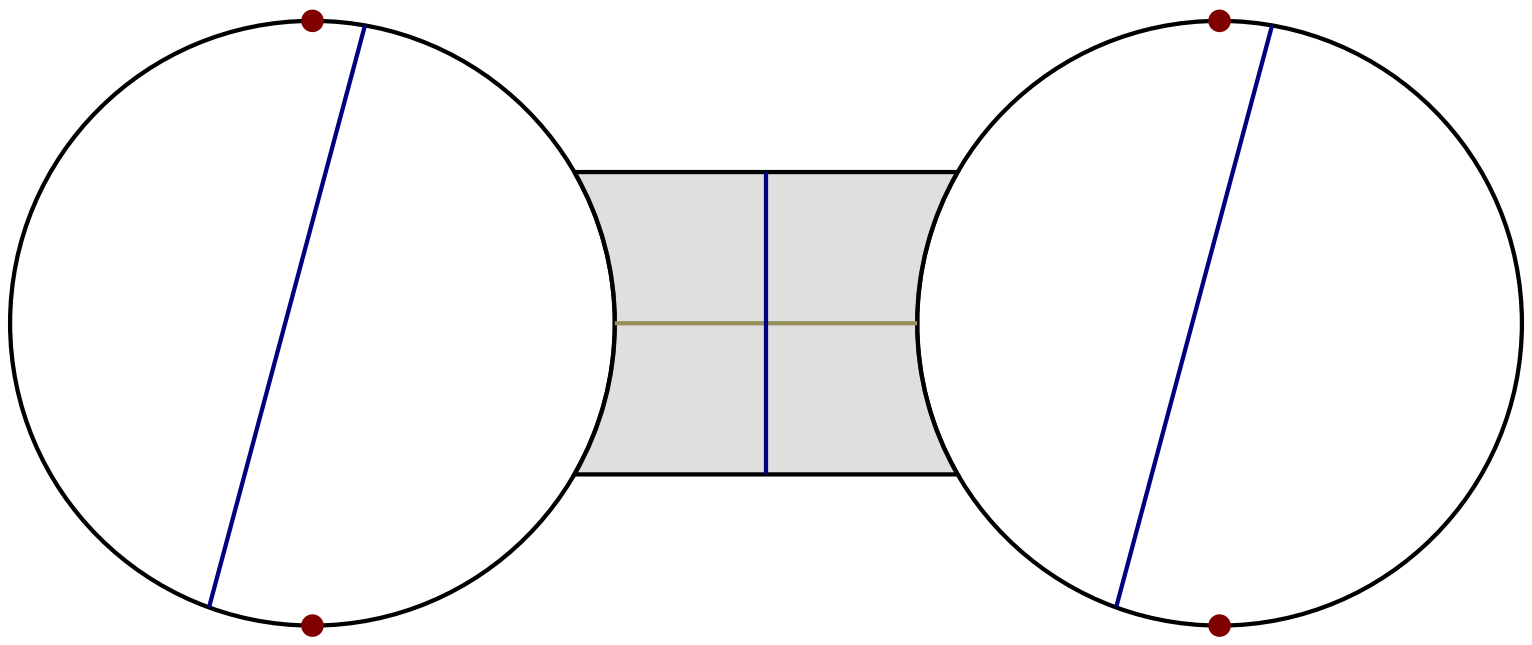
\begin{figure}[ht]
        \centering
        \begin{tikzpicture}
            \begin{scope}
                \draw [thick] (0,0) circle (2);
                \node [circle, fill=red!50!black, inner sep=1.5pt] at (0,2) {};
                \node [circle, fill=red!50!black, inner sep=1.5pt] at (0,-2) {};
            \end{scope}
            \begin{scope}[xshift=6cm]
                \draw [thick] (0,0) circle (2);
                \node [circle, fill=red!50!black, inner sep=1.5pt] at (0,2) {};
                \node [circle, fill=red!50!black, inner sep=1.5pt] at (0,-2) {};
            \end{scope}
            \filldraw[ fill=lightgray, fill opacity=0.5, draw=black, thick ] ({2*cos(30)},{2*sin(30)}) -- ({6+2*cos(150)},{2*sin(150)}) arc [start angle = 150, end angle = 210, radius = 2] -- ({2*cos(-30)},{2*sin(-30)}) arc [start angle = -30, end angle = 30, radius = 2];
            \draw [thick, yellow!50!black] (2,0) -- (4,0);
            \draw [thick, blue!50!black] ({2*cos(250)},{2*sin(250)}) -- ({2*cos(80)},{2*sin(80)});
            \draw [thick,blue!50!black] (3,{-2*sin(30)}) -- (3,{2*sin(30)});
            \draw [thick,blue!50!black] ({6+2*cos(250)},{2*sin(250)}) -- ({6+2*cos(80)},{2*sin(80)});
        \end{tikzpicture}
        \caption{Connecting Two ``Standard Copies'' by a $1$-Handle.}\label{fig:1-handle-attaching-Z}
    \end{figure}
    Before we compute the wrapped Floer cohomology of these Lagrangians, we fix one and label the vertical blue Lagrangians (corresponding to the cocore to the $1$-handle in figure \ref{fig:1-handle-attaching-Z}) in figure \ref{fig:generating-Lag-equiv-P1} by $L_0[k]$, $k\in\mathbb{Z}$, and the slanted Lagrangians (corresponding to the generating Lagrangian of $X_0$ in figure \ref{fig:1-handle-attaching-Z}) in figure \ref{fig:generating-Lag-equiv-P1} by $L_1 [k]$. They correspond under Homological Mirror Symmetry to the objects $\mathcal{O}_{\mathbb{P}^1}[k]$ and $\mathcal{O}_{\mathbb{P}^1} (1) [k]$, where $k$ is the character of the $\mathbb{G}_m$-action.\par
    It's clear from figure \ref{fig:1-handle-attaching-Z} that\[
        CW^{\ast} (L_1[k],L_1[l])=\left\{\begin{array}{ll}
            0,&\text{if}\ l\neq k;\\
            \Bbbk,&\text{if}\ l=k.
        \end{array}\right.\quad\text{and}\quad CW^{\ast} (L_0 [k],L_0[l])=\left\{\begin{array}{ll}
            0,&\text{if}\ l\neq k;\\
            \Bbbk,&\text{if}\ l=k.
        \end{array}\right.
    \]Let's conventionally set $L_1[0]$ to be the Lagrangian in the $1$-handle to the right of $X_0[0]$ where $L_1[0]$ lies in. Under this convention, it is also clear that \[
        CW^{\ast} (L_0[k],L_1 [l])=\left\{\begin{array}{ll}
            0,&\text{if}\ \vert k-l\vert\geq 2;\\ 
            \Bbbk,&\text{otherwise}.
        \end{array}\right.\quad\text{and}\quad CW^{\ast} (L_1[k],L_0[l])=0.
    \]Therefore we have completely compute the wrapped Fukaya category $\mathcal{W} (\mathbb{C}, e^z+e^{-z})$. To summarize, $\mathcal{W} (\mathbb{C},e^z+e^{-z})$ is quasi-equivalent to the category of complexes of perfect modules over path algebras of the infinite quiver 
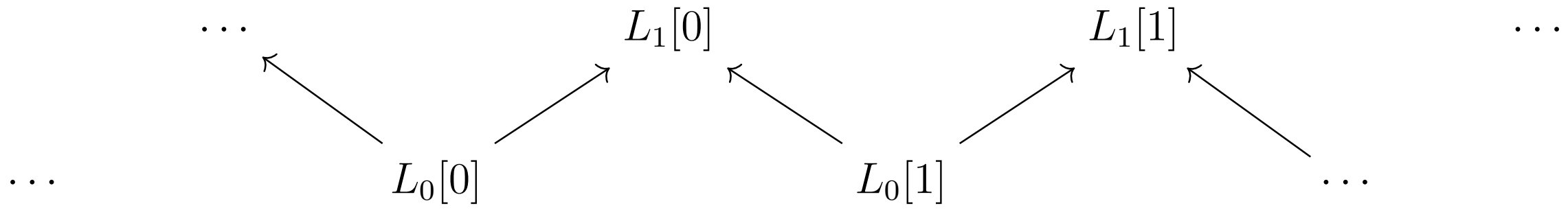
\[
        \begin{tikzcd}
            &\dotsb &&L_1[0] &&L_1 [1] &&\dotsb&\\
            \dotsb&&L_0[0]\arrow[ru]\arrow[lu]&&L_0[1]\arrow[ru]\arrow[lu]&&\dotsb\arrow[lu]
        \end{tikzcd}
    \]
\end{exmp}
\begin{exmp}\label{exmp:wrapped-Fukaya-cat-twice-punctured-equivariant}
    As another example, we will construct a $\mathbb{Z}^2$-cover $\mathcal{X}$ of the Weinstein sector $(\mathbb{C}\setminus\{\pm 1\},2/(1-z^2))$. This Weinstein sector can be regarded as the sectorial gluing of two Weinstein sectors $\left(\mathbb{C}^{\ast},z+z^{-1}\right)$ along one of their boundary stops, and its $\mathbb{Z}^2$-cover can be then realized as the gluing of $\mathbb{Z}$-copies of each sector along $\mathbb{Z}^2$-copies of strips, as shown in figure \ref{fig:Z2-cover-twice-punctured-disk}.\par
    \begin{figure}[ht]
        \centering
        \begin{tikzpicture}       
        \begin{scope}[yshift=6cm]
            \draw [thick] (0,0) -- (0,2);
            \draw [thick] (1,0) -- (1,2);
                \foreach \x in {0.2,0.8,1.4}{
                \filldraw [fill=white, draw=black,thick] (1,\x)
            -- (1.2,\x)
            .. controls (1.25,\x)
                       and (1.35,{0.98*\x+0.02*((\x-0.2)/0.2-4.6)})
            .. (1.4,{0.96*\x+0.04*((\x-0.2)/0.2-4.6)})
            -- (3.7,{0.04*\x+0.96*((\x-0.2)/0.2-4.6)})
            .. controls (3.75,{0.02*\x+0.98*((\x-0.2)/0.2-4.6)})
                       and (3.8,{(\x-0.2)/0.2-4.6})
            .. (4,{(\x-0.2)/0.2-4.6}) -- (4,{(\x-0.2)/0.2 -4.2})
            .. controls (3.7,{(\x-0.2)/0.2-4.2}) and (3.6,{((\x-0.2)/0.2-4.2)-0.2*((\x-0.2)/0.2-4.6-\x)/2.5})
            .. (3.5,{((\x-0.2)/0.2-4.2)-0.3*((\x-0.2)/0.2-4.6-\x)/2.5})
            -- (1.5,{\x+0.4+0.3*((\x-0.2)/0.2-4.6-\x)/2.5})
            .. controls (1.4,{\x+0.4+0.2*((\x-0.2)/0.2-4.6-\x)/2.5}) and (1.3,{\x+0.4})
            ..(1.2,{\x+0.4}) -- (1,{\x+0.4}) -- cycle;
            }
            \foreach \x in {0.4,1,1.6} {
                \node [circle,fill=red,inner sep=1.5pt] at (0,\x) {};
                \node [circle,fill=red,inner sep=1.5pt] at (1,\x) {};
            }
        \end{scope}
        \begin{scope}[yshift=3cm]
            \draw [thick] (0,0) -- (0,2);
            \draw [thick] (1,0) -- (1,2);
            \foreach \x in {0.2,0.8,1.4}{
                \filldraw [fill=white, draw=black, thick] (1,\x)
            -- (1.2,\x)
            .. controls (1.25,\x)
                       and (1.35,{0.98*\x+0.02*((\x-0.2)/0.2-2.2)})
            .. (1.4,{0.96*\x+0.04*((\x-0.2)/0.2-2.2)})
            -- (3.7,{0.04*\x+0.96*((\x-0.2)/0.2-2.2)})
            .. controls (3.75,{0.02*\x+0.98*((\x-0.2)/0.2-2.2)})
                       and (3.8,{(\x-0.2)/0.2-2.2})
            .. (4,{(\x-0.2)/0.2-2.2}) -- (4,{(\x-0.2)/0.2 -1.8})
            .. controls (3.7,{(\x-0.2)/0.2-1.8}) and (3.6,{((\x-0.2)/0.2-1.8)-0.2*((\x-0.2)/0.2-2.2-\x)/2.5})
            .. (3.5,{((\x-0.2)/0.2-1.8)-0.3*((\x-0.2)/0.2-2.2-\x)/2.5})
            -- (1.5,{\x+0.4+0.3*((\x-0.2)/0.2-2.2-\x)/2.5})
            .. controls (1.4,{\x+0.4+0.2*((\x-0.2)/0.2-2.2-\x)/2.5}) and (1.3,{\x+0.4})
            ..(1.2,{\x+0.4}) -- (1,{\x+0.4}) -- cycle;
            }
            \foreach \x in {0.4,1,1.6} {
                \node [circle,fill=red,inner sep=1.5pt] at (0,\x) {};
                \node [circle,fill=red,inner sep=1.5pt] at (1,\x) {};
            }
        \end{scope}
        \begin{scope}
            \draw [thick] (0,0) -- (0,2);
            \draw [thick] (1,0) -- (1,2);
            \foreach \x in {0.2,0.8,1.4}{
                \filldraw [fill=white, draw=black,thick] (1,\x)
            -- (1.2,\x)
            .. controls (1.25,\x)
                       and (1.35,{0.98*\x+0.02*((\x-0.2)/0.2+0.2)})
            .. (1.4,{0.96*\x+0.04*((\x-0.2)/0.2+0.2)})
            -- (3.7,{0.04*\x+0.96*((\x-0.2)/0.2+0.2)})
            .. controls (3.75,{0.02*\x+0.98*((\x-0.2)/0.2+0.2)})
                       and (3.8,{(\x-0.2)/0.2+0.2})
            .. (4,{(\x-0.2)/0.2+0.2}) -- (4,{(\x-0.2)/0.2 +0.6})
            .. controls (3.7,{(\x-0.2)/0.2+0.6}) and (3.6,{((\x-0.2)/0.2+0.6)-0.2*((\x-0.2)/0.2+0.2-\x)/2.5})
            .. (3.5,{((\x-0.2)/0.2+0.6)-0.3*((\x-0.2)/0.2+0.2-\x)/2.5})
            -- (1.5,{\x+0.4+0.3*((\x-0.2)/0.2+0.2-\x)/2.5})
            .. controls (1.4,{\x+0.4+0.2*((\x-0.2)/0.2+0.2-\x)/2.5}) and (1.3,{\x+0.4})
            ..(1.2,{\x+0.4}) -- (1,{\x+0.4}) -- cycle;
            }
            \foreach \x in {0.4,1,1.6} {
                \node [circle,fill=red,inner sep=1.5pt] at (0,\x) {};
                \node [circle,fill=red,inner sep=1.5pt] at (1,\x) {};
            }
        \end{scope}
        \begin{scope}[yshift=3cm, xshift=4cm]
            \draw [thick] (0,0) -- (0,2);
            \draw [thick] (1,0) -- (1,2);
            \foreach \x in {0.4,1,1.6} {
                \node [circle,fill=red,inner sep=1.5pt] at (0,\x) {};
                \node [circle,fill=red,inner sep=1.5pt] at (1,\x) {};
            }
        \end{scope}     
        \begin{scope}[xshift=4cm]
            \draw [thick] (0,0) -- (0,2);
            \draw [thick] (1,0) -- (1,2);
            \foreach \x in {0.4,1,1.6} {
                \node [circle,fill=red,inner sep=1.5pt] at (0,\x) {};
                \node [circle,fill=red,inner sep=1.5pt] at (1,\x) {};
            }
        \end{scope}
        \begin{scope}[yshift=6cm,xshift=4cm]
            \draw [thick] (0,0) -- (0,2);
            \draw [thick] (1,0) -- (1,2);
            \foreach \x in {0.4,1,1.6} {
                \node [circle,fill=red,inner sep=1.5pt] at (0,\x) {};
                \node [circle,fill=red,inner sep=1.5pt] at (1,\x) {};
            }
        \end{scope}
        \begin{scope}[yshift=9cm]
            \node at (0.5,1) {$\vdots$};
        \end{scope}
        \begin{scope}[yshift=-3cm]
            \node at (0.5,1) {$\vdots$};
        \end{scope}
        \begin{scope}[yshift=9cm,xshift=4cm]
            \node at (0.5,1) {$\vdots$};
        \end{scope}
        \begin{scope}[yshift=-3cm,xshift=4cm]
            \node at (0.5,1) {$\vdots$};
        \end{scope}
        \end{tikzpicture}
    \caption{$\mathbb{Z}^2$-cover of the Twice-Punctured Disk with Stops.}\label{fig:Z2-cover-twice-punctured-disk}
    \end{figure}
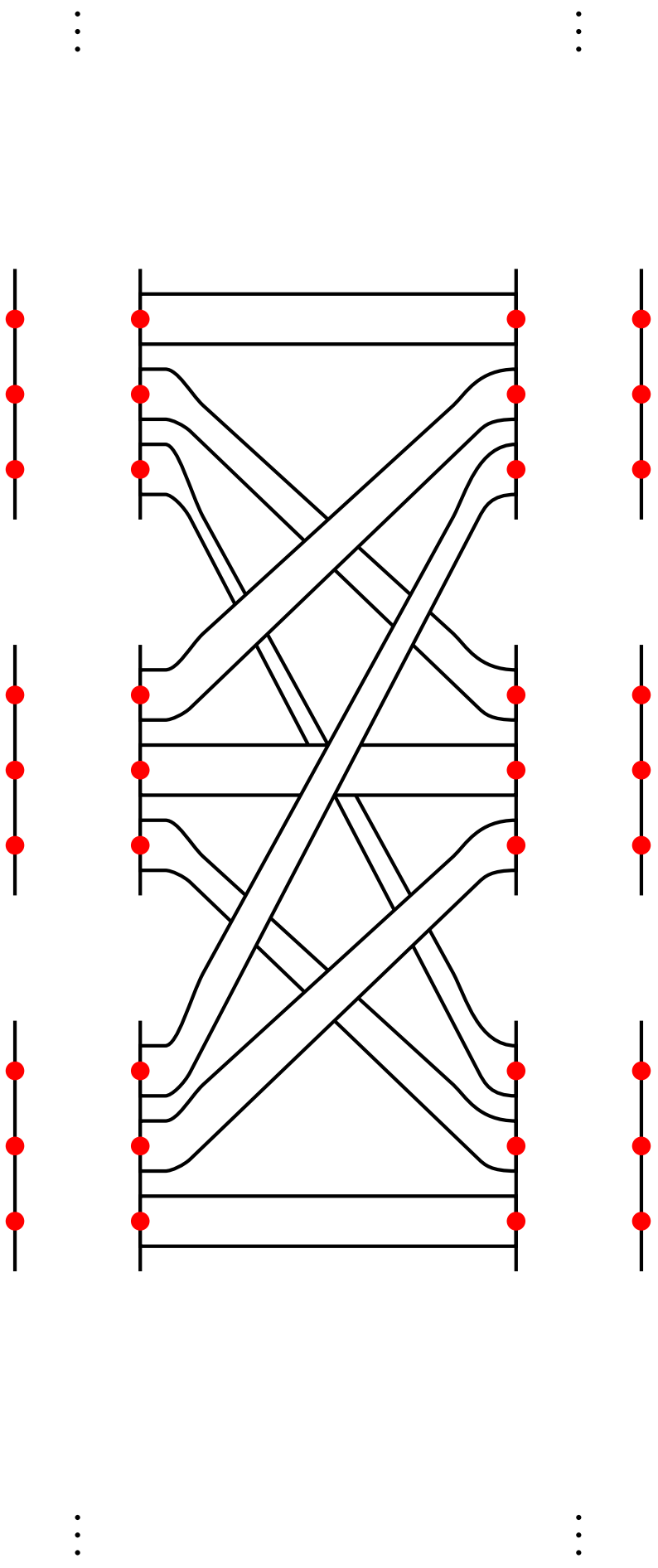
    As we know the generation of each sectors in this gluing, using sectorial descent for the diagram \[
        \begin{tikzcd}
            \mathcal{W} \left(\bigsqcup_{\mathbb{Z}^2}T^{\ast}[0,1]\right)\arrow[r]\arrow[d]&\mathcal{W} \left(\bigsqcup_{\mathbb{Z}}(\mathbb{C},e^{z_1}+e^{-z_1})\right)\\ 
            \mathcal{W}\left(\bigsqcup_{\mathbb{Z}}(\mathbb{C},e^{z_2}+e^{-z_2})\right)&
        \end{tikzcd}
    \]we see that $\mathcal{X}$ is still generated by the union of generators of the Fukaya categories \[\mathcal{W}\left(\bigsqcup_{\mathbb{Z}} (\mathbb{C},e^{z_i}+e^{-z_i})\right),\]which is again invariant under the $\mathbb{Z}^2$-action.\par
    Based on this gluing diagram, we can label generating Lagrangian submanifolds of $\mathcal{W}(\mathcal{X})$ as $L_i [a,b](c)$, where $i\in\{0,1\}$ denotes whether this Larangian is from the left-hand punctured disk or the right, $c\in\{0,1\}$ denotes the ``slope'' of $L_i$ and $(a,b)\in\mathbb{Z}^2$ denotes which connected component and after assigning a canonical $0$ which sheets in the $\mathbb{Z}$-cover it lies in.\par
    To better understand the Floer cohomologies $CW^{\ast}\left(L_i [a,b](c),L_j[c,d](c);\Bbbk\right)$, it's more intuitive to look at the twice-punctured disk model and use imagination for where the chosen Lagrangian gets lifted. See figure \ref{fig:computation-Z2-cover}. Suppose that $i=j$, as the wrapping around the middle circle clockwise shifts the bi-degree by $(-1,1)$, we see that \[
        CW^{\ast}\left(L_i [a_1,b_1](c_1), L_i [a_2,b_2](c_2)\right)\simeq\left\{\begin{array}{ll}
            \Bbbk,&\text{if}\ \vert a_1 -a_2\vert\leq 1,\ b_1= b_2\ \text{and}\ c_2-c_1=1;\\ 
            \Bbbk,&\text{if}\ a_1-b_1=b_2-a_2\ \text{and}\ c_1=c_2;\\ 
            0,&\text{otherwise}.
        \end{array}\right.
    \]
    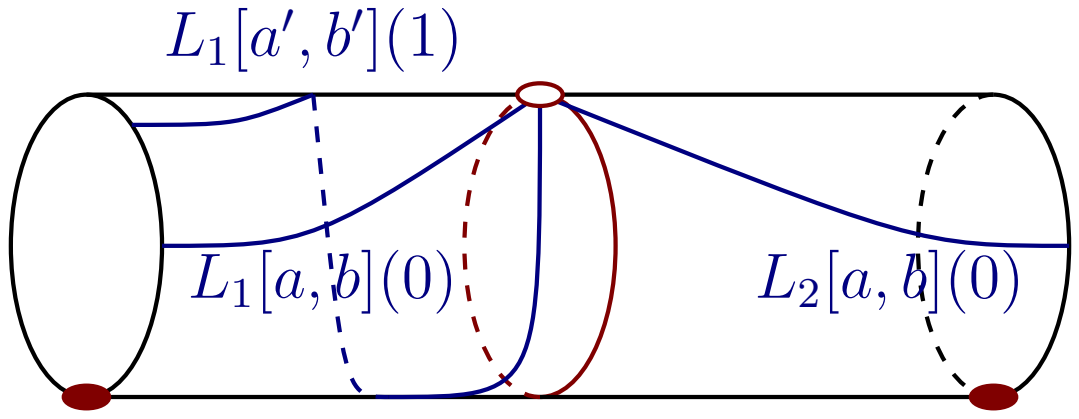
\begin{figure}[ht]
        \centering
        \begin{tikzpicture}
            \begin{scope}
        \draw [thick] (-2,0) arc[start angle=0, end angle=360, x radius=0.5, y radius=1];
        \draw [thick] (-2.5,1) -- (0.5,1);
        \draw [thick] (-2.5,-1) -- (0.5,-1);
        \draw [thick,red!50!black] (0.5,1) arc[start angle=90, end angle=-90, x radius=0.5, y radius=1];
        \draw [thick,red!50!black, dashed] (0.5,-1) arc[start angle=270, end angle=90, x radius=0.5, y radius=1];
        \draw [thick,blue!50!black] (-2,0) .. controls (-1,0) .. (0.5,1) node [pos=0.5,below] {$L_1[a,b](0)$};
        \draw [thick,blue!50!black] (-2.2,{sqrt(1-(0.3/0.5)^2)}) .. controls (-1.5,{sqrt(1-(0.3/0.5)^2)}) .. (-1,1) node [above] {$L_1[a',b'](1)$};
        \draw [thick,dashed,blue!50!black] (-1,1) .. controls (-0.8,-1) .. (-0.5,-1);
        \draw [thick,blue!50!black] (-0.5,-1) .. controls (0.5,-1) .. (0.5,1);
        \end{scope}
        \node [ellipse,
    minimum width=3mm,
    minimum height=1.5mm,
    inner sep=0pt,
    fill=red!50!black,
    draw=red!50!black,
    thick] at (-2.5,-1) {};
    \begin{scope}[xshift=3cm]
        \draw [thick] (-2.5,1) -- (0.5,1);
        \draw [thick] (-2.5,-1) -- (0.5,-1);
        \draw [thick] (0.5,1) arc[start angle=90, end angle=-90, x radius=0.5, y radius=1];
        \draw [thick, dashed] (0.5,-1) arc[start angle=270, end angle=90, x radius=0.5, y radius=1];
        \node [ellipse,
    minimum width=3mm,
    minimum height=1.5mm,
    inner sep=0pt,
    fill=red!50!black,
    draw=red!50!black,
    thick] at (0.5,-1) {};
    \draw [thick,blue!50!black] (1,0) .. controls (0,0) .. (-2.5,1) node [pos=0.5,below] {$L_2[a,b](0)$};
    \end{scope}
    \node [ellipse,
    minimum width=3mm,
    minimum height=1.5mm,
    inner sep=0pt,
    fill=white,
    draw=red!50!black,
    thick] at (0.5,1) {};
        \end{tikzpicture}
        \caption{Computing wrapped Floer cohomology from the base space.\label{fig:computation-Z2-cover}}
    \end{figure}
    Similarly, when $i\neq j$, it can be read from figures \ref{fig:Z2-cover-twice-punctured-disk}, \ref{fig:computation-Z2-cover} that \[
        CW^{\ast}\left(L_i [a_0, b_0](c_i),L_j[a_1,b_1](c_j);\Bbbk\right)=\left\{\begin{array}{ll}
            \Bbbk,&\text{if}\ a_1-b_1 =b_0-a_0;\\ 
            0,&\text{otherwise}.
        \end{array}\right.
    \]This provides a complete characterization of the wrapped Fukaya category. The quiver presentation will consist of a 3D graph, and is omitted here. Note that some parallel description of the same space also appears in \cite{Neumann04} in the study of extended Bloch groups.
\end{exmp}
We end this section with a discussion of a higher-dimension situation.
Recall that for a smooth manifold $M$, with a choice of a Morse function $f\colon M\to\mathbb{R}$ on $M$ and a Riemannian metric $g$ the cotangent bundle $T^{\ast} M$ carries a natural Weinstein structure $(\phi,\lambda_f,\omega_f)$ where $\phi (x,\alpha) =f(x)+\vert\alpha\vert_g^2$, $\lambda_f =df-\lambda_{\can}$ with $\lambda_{\can}$ the canonical Liouville structure on $T^{\ast} M$ given by $\pi^{\ast}\lambda_{\can} =\lambda_{\can}$ where $\pi\colon T^{\ast} M\to M$ is the bundle map, and $\omega_f =d\lambda_f$.\par
A Legendrian embedding $\Lambda\hookrightarrow\partial_{\infty} T^{\ast} M$ provides a symplectic hypersurface $T^{\ast}\Lambda\hookrightarrow T^{\ast} M$, which itself carries a Weinstein structure as above.
Due to \cite{Eliashberg18}, after modifying the Weinstein structure on $T^{\ast} M$ we can get a Weinstein structure making $T^{\ast}\Lambda\subseteq T^{\ast} M$ a Weinstein hypersurface, defining a Weinstein sector due to \cite{GPS20}.\par
One example appears in \cite{FLTZ12}, where $M=\mathbb{T}^n$ and for a given polygon $\Sigma$ of dimension $n$ we can construct a Legendrian link $\Lambda_{\Sigma}\hookrightarrow\partial_{\infty} T^{\ast}\mathbb{T}^n$ to obtain a Weinstein sector $(T^{\ast}\mathbb{T}^n,\Lambda_{\Sigma})$.
The $\mathbb{Z}^n$-cover $(T^{\ast}\mathbb{R}^n,\widetilde{\Lambda}_{\Sigma})$ with the pull-back Weinstein structure is also naturally a Weinstein sector, by our definition.\par
Note that $\pi_1 (T^{\ast}\mathbb{T}^n,\ast)\cong\mathbb{Z}^n$, hence in the $\mathbb{Z}^n$-cover the $1$-handle attachments are altered, while all remaining handles are lifted trivially.
Write $X=(T^{\ast}\mathbb{R}^n,\widetilde{\Lambda}_{\Sigma})$, we can then construct an exhausting sequence of Weinstein cobordisms $X_0\hookrightarrow X_1\hookrightarrow X_2\hookrightarrow\dotsb\hookrightarrow X$ as follows.\par
$X_0$ is a disjoint union of finitely many points, with labels $0$ or $e_i$ where
$$e_i =(0,\dotsb ,0,\underbrace{1}_i,0,\dotsb ,0).$$
$X_1$ is obtained from $X_0$ by attaching $1$-handles to $X_0$ connecting the point $0$ to all the $e_i$, all the remaining handles of higher index lying over these $1$-handles, and half-handles attached to the given handles that contribute to the Legendrian stops near the contact boundary-at-infinity.
$X_2$ is constructed from $X_1$ by adding $0$-handles with labels $e_i+e_j$, $1\leq i,j\leq n$, and repeating the procedure of attaching handles described in the construction of $X_1$.
This procedure can then be continued and we will exhaust all the handles for the lifted Weinstein manifold $(T^{\ast}\mathbb{R}^n,\widetilde{\Lambda}_{\Sigma})$.
It's immediate from the definition that
\begin{lmm}\label{lmm:Fukaya-category-Gmn-equivariant-generation}
    The wrapped Fukaya category $\mathcal{W} (T^{\ast}\mathbb{R}^n,\widehat{\Lambda}_{\Sigma})$ is generated by cotangent fibres $T^{\ast}_q\mathbb{R}^n$ where $q\in\mathbb{Z}^n$ and all linking disks to the (infinite type) Legendrian link $\Lambda_{\Sigma}$ near infinity.
\end{lmm}
On the other hand, we can construct from the data $(T^{\ast}M,\Lambda)$ dg-category of constructible sheaves of categories $\Sh_{\Lambda} (M)$ with singular supports contained in the Whitney stratification $\mathcal{S}_{\Lambda}$ of conic isotropic submanifolds consisting of the zero section together with positive real cones over the Legendrians $\Lambda$ near infinity valued in the dg-category of chain complexes over $\Bbbk$ \cite{GPS24M}.
It's proved in \cite{GPS24M} that
\begin{thm}[Sheaf-Floer correspondence]\label{thm:sheaf-Floer-correspondence}
    Suppose that $M$ is a smooth manifold with $f\colon M\to\mathbb{R}$ a Morse function with finitely many critical points. There is a canonical quasi-equivalence of categories $$
        \Perf\mathcal{W} (T^{\ast} M,\Lambda)^{op}\simeq\Sh_{\Lambda} (M)^{\omega}
    $$where linking disk at a point $p\in\Lambda$ are sent to co-representatives of microstalk functor associated to $p$ and cotangent fibres $T^{\ast}_p M$ with $p$ not in the image of the projection of $\Lambda$ under $\partial_{\infty} T^{\ast} M\to M$ sent to a co-representative of the stalk functor at $p$.
\end{thm}
We will provide a generalization of this theorem to infinite-type sectors in this specific situation.
\begin{prp}\label{prp:inf-type-sheaf-Floer-Rn}
    For $M=\mathbb{R}^n$ and $\widetilde{\Lambda}_{\Sigma}$ the lifted FLTZ skeleton, there is a quasi-equivalence of $A_{\infty}$-categories \[
        \Perf\mathcal{W} (T^{\ast}\mathbb{R}^n,\widetilde{\Lambda}_{\Sigma})^{op}\simeq\Sh_{\widetilde{\Lambda}_{\Sigma}} (\mathbb{R}^n)^{\omega}.
    \]
\end{prp}
\begin{proof}
    Note that the exhaustion $X_0\hookrightarrow X_1\hookrightarrow\dotsb\hookrightarrow X$ by construction are exhaustions by cotangent bundles with Legendrian stops, and therefore for each $X_i=(T^{\ast} M_i,\Lambda_i)$, we have by \cite{GPS24M} a canonical quasi-equivalence of categories \[
        \Perf\mathcal{W} (X_i)^{op}\simeq\Sh_{\Lambda_i} (M_i)^{\omega}.
    \]Under this quasi-equivalence, Viterbo functor is given by the restriction functor \[
    \Sh_{\Lambda_{i+1}} (M_{i+1})\to\Sh_{\Lambda_i} (M_i)
    \]sending a constrctible sheaf to its restriction to an open subset removing some of the strata, which sends compact objects to compact objects using the description \cite[Lemma 4.18]{GPS24M}.
    We can then verify directly using the exhaustion that we can present $\Sh_{\widetilde{\Lambda}_{\Sigma}} (\mathbb{R}^n)$ as a homotopy limit \[
        \Sh_{\widetilde{\Lambda}_{\Sigma}} (\mathbb{R}^n)\simeq\holim\left(\Sh_{\Lambda_0} (M_0)\leftarrow\Sh_{\Lambda_1} (M_1)\leftarrow\dotsb\right).
    \]\cite[Proposition 4.20]{GPS24M} together with the definition of wrapped Fukaya categories then implies the desired quasi-equivalence.
\end{proof}
\subsection{Wrapped Fukaya category of quotients.} Let $L,K\subseteq X/G$ be exact cylindrical spin Lagrangian submanifolds such that $\langle\pi_1 (L),\pi_1 (K)\rangle\cap G= \{1\}$ in $\pi_1 (X/G)$, then we have
\begin{prp}
    $CF^{\ast} (L,K) =\displaystyle\bigoplus_{g\in G} CF^{\ast} (g(\tilde{L}),\tilde{K})$ as graded $\Bbbk$-modules. Moreover, for given positive isotopy $\psi\colon K\leadsto K^+$ with lift $\tilde{\psi}$, the continuation map $c_{\psi}\colon CF^{\ast} (L,K)\to CF^{\ast} (L,K^+)$ coincides with $\displaystyle\bigoplus_{g\in G} c_{\tilde{\psi}}\colon \bigoplus_{g\in G} CF^{\ast} (g(\tilde{L}),\tilde{K})\to\bigoplus_{g\in G} CF^{\ast} (g(\tilde{L}),\tilde{K}^+)$.
\end{prp}
\begin{proof}
    An intersection point $p\in L\cap K$ in $X/G$ corresponds to a point $\tilde{L}\cap g(\tilde{K})$ in $X$ for some $g\in G$, and therefore by definition of Floer cochain complex we have $CF^{\ast} (L,K)=\displaystyle\bigoplus_{g\in G} (g(\tilde{L}), \tilde{K})$. Given a positive isotopy $\psi\colon K\leadsto K^+$ with a given lift $\tilde{K}$, using the fact that $X\to X/G$ is a covering space we can lift the whole isotopy $\psi$ to an isotopy $\tilde{\psi}\colon \tilde{K}\leadsto\tilde{K}^+$. Under the $G$-action we have disjoint copies of isotopies $g(\tilde{\psi})\colon g(\tilde{K})\leadsto g(\tilde{K}^+)$, and therefore provides Floer-theoretic continuation maps \[
        c_{g(\tilde{\psi})}\colon CF^{\ast} (g(\tilde{L}),\tilde{K})\to CF^{\ast} (g(\tilde{L}),\tilde{K}^+)
    \]for each chosen Lagrangian $\tilde{L}$ in $X$. Taking the direct sum provides the desired continuation map.
\end{proof}
Now given a sequence of admissible Lagrangian submanifolds $L_0, L_1,\dotsb ,L_k$ in $X/G$ with intersection points $p_i\in CF^{\ast} (L_{i-1}, L_i)$ for $1\leq i\leq k$ and $p_0\in L_k\cap L_0$, a regular $J$-holomorphic disk $u\in\mathcal{M} (p_0;p_1,\dotsb ,p_k)$ admits a unique natural lift to an element \[
\tilde{u}\in\mathcal{M}(X;\tilde{p}_0;\tilde{p}_1,\dotsb ,\tilde{p}_k)
\]where for fixed chosen lift $\tilde{L}_i$ of $L_i$, the intersection points $\tilde{p}_i\in g_i(\tilde{L}_{i-1})\cap\tilde{L}_i$ lift to intersection points of $g_i (\tilde{L}_{i-1})$ and $\tilde{L}_i$ for some $g_i$ modulo $k$, and hence the lifted holomorphic disk $\tilde{u}$ contributes to a map \[
\mu^k\colon CF^{\ast} (g_{k-1}(\tilde{L}_{k-1}),\tilde{L}_k)\otimes\dotsb\otimes CF^{\ast} (g_1\dotsb g_{k-1}(\tilde{L}_0),g_2\dotsb g_{k-1} (\tilde{L}_1))\to CF^{\ast} (g_1\dotsb g_k (\tilde{L}_0),\tilde{L}_k)[2-k].
\]
Conversely, any regular $J$-holomorphic disk in $u$ with boundary on $\displaystyle\bigsqcup_{g\in G}g(\tilde{L})$ asymptotic to intersections $g(\tilde{L}_{i-1})\cap\tilde{L}_i$ projects to regular $J$-holomorphic disks in $X/G$ with boundaries on $L_i$, and therefore the $A_{\infty}$-structure on $\mathcal{W} (X/G)$ can be described as $A_{\infty}$-structures of $\mathcal{W} (X)$ that is invariant under the $G$-action.\footnote{The same argument was also used in Seidel's work \cite{Seidel11}.} This is enough to conclude the theorem.
\begin{proof}[Proof of Theorem \ref{thm:wrapped-Fukaya-quotient}]
    Note from definition \ref{defn:wrapped-Fukaya-cat-inf-sectors} that the wrapped Fukaya category $\mathcal{W} (X)$ is generated by all the cocores appearing in any of the exhausting sequence of Weinstein domains. In particular, the $G$-invariant Weinstein structure on $X$ provides $G$-copies of critical points associating to a single critical point of any given Weinstein structure on $X/G$. Therefore we conclude that the disjoint union $\displaystyle\bigsqcup_{g\in G} g(\tilde{L})$ of cocores, where $\tilde{L}$ is a lift of cocore of the Weinstein structure on $X/G$, generates the wrapped Fukaya category of $X$. On the chain level we have \[
        CW^{\ast}_{X/G}\left(L,K;\Bbbk\right)=\bigoplus_{g\in G} CW^{\ast}_X (g(\tilde{L}),\tilde{K};\Bbbk),
    \]where the $A_{\infty}$-structures are written as the discussion above, which pairs exactly with the definition of $A_{\infty}$-orbit categories in Definition \ref{defn:Ainf-orbit-category}.
\end{proof}
\subsection{The K\"unneth theorem and deformation invariance.} Recall from \cite{GPS24S} that we have the following K\"unneth theorem for wrapped Fukaya categories:
\begin{thm}[{\cite[Theorem 1.5]{GPS24S}}]\label{thm:Kunneth-theorem-Fukaya-cat}
    For Weinstein sectors $X$ and $Y$ of finite type, there is a quasi-equivalence of $A_{\infty}$-categories \[
        \mathcal{W} (X)\otimes\mathcal{W} (Y)\simeq\mathcal{W} (X\times Y).
    \]
\end{thm}
The result we need to use in our situation is slightly more general: it allows $X$ to be of infinite type, and so $\mathcal{W} (X)$ is defined instead in an ad hoc way. Nonetheless, we can still achieve the following:
\begin{crl}\label{crl:Kunneth-Fukaya-infinite}
    Suppose that $X$ is a Weinstein sector and $Y$ is another Weinstein sector of finite type, then there is a quasi-equivalence of $A_{\infty}$-categories \[
        \mathcal{W} (X)\otimes\mathcal{W} (Y)\simeq\mathcal{W} (X\times Y).
    \]
\end{crl}
\begin{proof}
    Let $\emptyset=X_{-1}\subseteq X_0\subseteq\dotsb\subseteq X$ be an exhaustion of $X$ by Weinstein domains, then for each $n$ we have quasi-equivalences of $A_{\infty}$-categories \[
        \mathcal{W} (X_n)\otimes\mathcal{W} (Y)\xrightarrow{\simeq}\mathcal{W} (X_n\times Y).
    \]Note that Ornaghi's tensor product \cite{Ornaghi25} is functorial, and so the Viterbo restriction $\mathcal{W} (X_{n+1})\twoheadrightarrow\mathcal{W} (X_n)$ gives rise to $A_{\infty}$-functors \[
        \mathcal{W} (X_{n+1})\otimes\mathcal{W} (Y)\twoheadrightarrow\mathcal{W} (X_n)\otimes\mathcal{W} (Y)
    \]which fits into a homotopically commutative diagram of $A_{\infty}$-categories \[
        \begin{tikzcd}
            \mathcal{W} (X_{n+1})\otimes \mathcal{W} (Y)\arrow[r,"\simeq"]&\mathcal{W} (X_{n+1}\times Y)\\ 
            \mathcal{W} (X_n)\otimes\mathcal{W} (Y)\arrow[u,two heads]\arrow[r,"\simeq"]&\mathcal{W} (X_n\times Y)\arrow[u,two heads]
        \end{tikzcd}
    \]Note that $\{X_n\times Y\}_n$ is an exhaustion of $X\times Y$ by Weinstein subdomains, and therefore by taking homotopy limits we get the quasi-equivalence \[
        \mathcal{W}^{big} (X)\otimes\mathcal{W}(Y)\simeq\mathcal{W}^{big} (X\times Y),
    \]and by keeping track of the ad-hoc definition of \ref{defn:wrapped-Fukaya-cat-inf-sectors} we get the desired K\"unneth isomorphism.
\end{proof}
The other essential ingredient is a proof that certain exact symplectic deformation of an infinite-type Weinstein sector does not change the corresponding Fukaya category.
\begin{prp}\label{prp:deformation-invariance-Fukaya}
    Suppose that $X$ is a Weinstein sector with Weinstein structure $\phi$ and $\{\phi_t\}$ a Weinstein homotopy with $\phi_0=\phi$ such that for each $n$, zeroes of $\phi_t\vert_{X_n}$ are contained in a compact subset in the interior of $X_n$. Then $\mathcal{W} (X,\phi_0)\simeq\mathcal{W} (X,\phi_1)$.
\end{prp}
\begin{proof}
    From Lemma 2 of \cite{Jeffs22} we know that $\mathcal{W} (X_n,\phi)\simeq\mathcal{W} (X_n,\phi_1)$ for each $n$, and by definition of the Weinstein homotopy we see that this quasi-isomorphism is compatible with the Viterbo restriction. Therefore by taking the homotopy limit we get the desired deformation invariance.
\end{proof}
\subsection{Kartal's modification.} We end this section with a recall of Kartal's modification \cite{Kartal21} adapted to the current situation. Let $X$ be a Weinstein sector of finite type with Weinstein structure $(\varphi,\lambda)$ and an exact symplectomorphism $\phi\colon X\to X$ that is contact type in a cylindrical end of $X$, we are going to construct a Weinstein structure $\psi$ on $(T^{\ast}\mathbb{R}\times X,e^z+e^{-z})$ such that $\psi$ is invariant under the $\tau^{-1}\times\phi$-action, with $\tau\colon (T^{\ast}\mathbb{R},e^z+e^{-z})\to (T^{\ast}\mathbb{R},e^z+e^{-z})$ the translation-by-$1$ symplectomorphism.\par
The ansatz, from \cite{Kartal21}, is to consider Weinstein structures of the form\[
\pi_X^{\ast}\varphi_z +p^{\ast}\varphi_{(T^{\ast}\mathbb{R},e^z+e^{-z})}
\]where $\{\varphi_z\}$ is a $T^{\ast}\mathbb{R}$-family of Weinstein structures on $X$ so that $\phi^{\ast}\varphi_{z+1}=\varphi_z$ and $p\colon (T^{\ast}\mathbb{R}\times X,e^z+e^{-z})\to (T^{\ast}\mathbb{R},e^z+e^{-z})$ is the projection to the $(T^{\ast}\mathbb{R},e^z+e^{-z})$-factor. We can further assume that $\varphi_{z=x+iy}$ is $y$-invariant and is constant when $x\in\mathbb{Z}$ and so the resulting Weinstein deformation does not change the stop and we can simply write $\{\varphi_s\}_{s\in\mathbb{R}}$. By definition of $\phi$ we are able to find a smooth function $K\colon X\to\mathbb{R}$ where $\phi^{\ast}\lambda_X=\lambda_X+dK$, and the construction of $\{\lambda_s\}$ is given by considering a smooth family of functions $\{K_s\}_{s\in\mathbb{R}}$ that is cylindrical in a cylindrical end of $X$ so that $\lambda_s =\lambda +dK_s$. The condition $\phi^{\ast}\lambda_{s+1}=\lambda_s$ then turns into the condition that \[
    dK_s=d(\phi^{\ast} K_{s+1}+K).
\]We can then find a smooth family of functions $K_s$ such that $\phi^{\ast} K_{s+1} +K-K_s=0$ holds for all $s\in\mathbb{R}$, by exactly the same procedure as described in the proof of Proposition 2.1 in \cite{Kartal21}.\par
To upgrade this Liouville structure to a Weinstein structure, we need to further choose a family of smooth functions $\varphi_s$ together with a constant $C>0$ so that $\varphi_C =p^{\ast}\varphi_{(T^{\ast}\mathbb{R},e^z+e^{-z})} +C^{-1}\varphi_z$ is Lyapunov for the Liouville structure $\lambda_C =p^{\ast}\lambda_{(T^{\ast}\mathbb{R},e^z+e^{-z})} +C^{-1}\lambda_z$. The conclusions drawn in \cite{Kartal21} is then 
\begin{prp}[{Similar to Proposition 2.3 of \cite{Kartal21}}]\label{prp:Kartal-modification-Weinstein}
    There exists a Weinstein structure on $(T^{\ast}\mathbb{R}\times X,e^z+e^{-z})$ of the form \[
        \varphi_{(T^{\ast}\mathbb{R}\times X,e^z+e^{-z})} =p^{\ast}\varphi_{(T^{\ast}\mathbb{R},e^z+e^{-z})}+C^{-1}\varphi_z,\quad \lambda_{(T^{\ast}\mathbb{R}\times X,e^z+e^{-z})}=\pi^{\ast}\lambda_{(T^{\ast}\mathbb{R},e^z+e^{-z})} +C^{-1}\lambda_z
    \]for $C>0$ sufficiently large, that is invariant under the $\tau^{-1}\times\phi$-action on $(T^{\ast}\mathbb{R}\times X,e^z+e^{-z})$.
\end{prp}

\section{$G$-equivariant sheaves\label{sec:G-equiv-sheaves}}
In this section we recall the definition of $G$-equivariant sheaves (from e.g. \cite[\href{https://stacks.math.columbia.edu/tag/03LE}{Tag 03LE}]{stacks-project}) and establish the relation between $G$-equivariant coherent sheaves/vector bundles over a scheme $X$ and $H$-equivariant coherent sheaves/vector bundles over $X$, with $H\subseteq G$ a subgroup.

\subsection{Categories of $G$-equivariant $\mathcal{O}_X$-modules.}
Let $G=(G, m, e)$ a group scheme over $\field$ with multiplication $m : G \times G \rightarrow G$ and identity section $e : \Spec(\field) \rightarrow G$,
and
$X$ a scheme over $\field$
equipped with a $G$-action
$\sigma : G \times X \rightarrow X$. See \cite[\href{https://stacks.math.columbia.edu/tag/022R}{Section 022R}]{stacks-project} and \cite[\href{https://stacks.math.columbia.edu/tag/01II}{Section 01II}]{stacks-project} for some of the basic definitions.

\begin{defn}\label{defn:EquivariantSheaves}
A \emph{$G$-equivariant $\mathcal{O}_X$-module}
$(\mathcal{F}, \theta)$
is an $\mathcal{O}_X$-module $\mathcal{F}$
together with an isomorphism
$$
\theta : \sigma^* \mathcal{F} \xrightarrow{\cong} \left(\pr_X^{G \times X}\right)^* \mathcal{F}
$$
of $\mathcal{O}_{G\times X}$-modules
(called a \emph{$G$-linearization} of $\mathcal{F}$),
where $\pr_X^{G \times X} : G \times X \rightarrow X$ is the canonical projection,
such that the following conditions hold:
\begin{itemize}
\item \emph{(cocycle condition)}
$
(m \times \id_X)^* \theta
=
(\id_G \times \sigma)^* \theta \comp \left(\pr_{G \times X}^{G \times G \times X}\right)^* \theta.
$

\item \emph{(unit condition)}
$(e \times \id_X)^*\theta = \id_{\mathcal{F}}$.

\end{itemize}

A \emph{morphism of $G$-equivariant $\mathcal{O}_X$-modules}
from
$\mathcal{F}=(\mathcal{F}, \theta)$
to
$\mathcal{F}'=(\mathcal{F}', \theta')$
is a morphism
$
f : \mathcal{F} \rightarrow \mathcal{F}'
$
of $\mathcal{O}_X$-modules
that commutes with the $G$-linearizations,
i.e.,
$$
\sigma^* f \comp \theta'
=
\theta \comp \left(\pr_X^{G \times X}\right)^* f
$$
as morphisms of $\mathcal{O}_{G \times X}$-modules.


We denote by
$\Mod_G(X)$
the \emph{category of $G$-equivariant $\mathcal{O}_X$-modules}.
A $G$-equivariant $\mathcal{O}_X$-module $(\mathcal{F},\theta)$
is said to be
\emph{quasi-coherent}, \emph{coherent},
or \emph{locally free of finite rank}
if $\mathcal{F}$
is an $\mathcal{O}_X$-module having the corresponding property.
We denote by
$\QCoh_G(X)$,
$\Coh_G(X)$,
and
$\VB_G(X)$
the full subcategories of $\Mod_G(X)$
whose objects are
$G$-equivariant quasi-coherent, coherent, and finite rank locally free $\mathcal{O}_X$-modules,
respectively,
so that they satisfy
$$
\VB_G(X)\subset\Coh_G(X)\subset\QCoh_G(X)\subset\Mod_G(X).
$$
\end{defn}

\begin{defn}[Internal tensor product of equivariant sheaves]
\label{defn:TensorProductOfEquivariantModules}
Let $(\mathcal{F}, \theta)$ and $(\mathcal{F}', \theta')$ be two $G$-equivariant $\mathcal{O}_X$-modules in $\Mod_G(X)$.
The \emph{(internal) tensor product} of $(\mathcal{F}, \theta)$ and $(\mathcal{F}', \theta')$
$$
(\mathcal{F}, \theta)
\underset{\mathcal{O}_X}{\otimes}
(\mathcal{F}', \theta')
:=
\left(\mathcal{F}\underset{\mathcal{O}_X}{\otimes}\mathcal{F}', \, \theta\otimes\theta'\right)
\in
\Mod_G(X),
$$
is defined as the $\mathcal{O}_X$-module $\mathcal{F}\underset{\mathcal{O}_X}{\otimes}\mathcal{F}'$
together with the $G$-linearization
$$
\theta\otimes\theta'
:
\sigma^*\left(
\mathcal{F}
\underset{\mathcal{O}_X}{\otimes}
\mathcal{F}'
\right)
\xrightarrow{\cong}
\left(\pr_X^{G \times X}\right)^*
\left(\mathcal{F}
\underset{\mathcal{O}_X}{\otimes}
\mathcal{F}'\right)
$$
given by the tensor product $\theta \otimes \theta'$
under the canonical isomorphisms
$$
\sigma^*\left(\mathcal{F}\underset{\mathcal{O}_X}{\otimes}\mathcal{F}'\right)
\cong
\sigma^* \mathcal{F} \underset{\mathcal{O}_{G \times X}}{\otimes} \sigma^* \mathcal{F}'
\quad\text{and}\quad
\left(\pr_X^{G \times X}\right)^* \left(\mathcal{F}\underset{\mathcal{O}_X}{\otimes}\mathcal{F}'\right)
\cong
\left(\pr_X^{G \times X}\right)^*\mathcal{F} \underset{\mathcal{O}_{G \times X}}{\otimes} \left(\pr_X^{G \times X}\right)^*\mathcal{F}'.
$$

The \emph{unit object} for the tensor product is defined by
$(\mathcal{O}_X, \id_{\mathcal{O}_{G\times X}}) \in \Mod_G(X)$,
where
$
\id_{\mathcal{O}_{G\times X}}
:
\sigma^* \mathcal{O}_X
\xrightarrow{\cong}
\left(\pr_X^{G \times X}\right)^* \mathcal{O}_X
$
is induced from the canonical isomorphisms of sheaves
$
\sigma^* \mathcal{O}_X
\cong
\mathcal{O}_{G\times X}
$
and
$
\left(\pr_X^{G \times X}\right)^* \mathcal{O}_X
\cong
\mathcal{O}_{G\times X}.
$

Equipped with this tensor product,
$\Mod_G(X)$
forms a
\emph{symmetric monoidal category}
$
\left(
\Mod_G(X),
\underset{\mathcal{O}_X}{\otimes},
\left(\mathcal{O}_X,\theta(0)\right)
\right).
$
Furthermore, the operation $\underset{\mathcal{O}_X}{\otimes}$ restricts to endow
$\QCoh_G(X)$, $\Coh_G(X)$, and $\VB_G(X)$
with symmetric monoidal category structures as well.

\end{defn}

\begin{defn}[Internal hom of equivariant sheaves]
\label{defn:InternalHomOfEquivariantModules}
Let $(\mathcal{F}, \theta)$ and $(\mathcal{F}', \theta')$ be two $G$-equivariant $\mathcal{O}_X$-modules in $\Mod_G(X)$
such that
$\mathcal{F}$ is finitely presented\footnote{
An $\mathcal{O}_X$-module $\mathcal{F}\in\Mod(X)$
is finitely presented if
(1)
$\mathcal{F}\in\VB(X)$ is finite rank locally free,
or
(2)
$X$ is locally Noetherian and $\mathcal{F}\in\Coh(X)$ is a coherent sheaf.
}.
The \emph{internal hom} of $(\mathcal{F}, \theta)$ and $(\mathcal{F}', \theta')$
$$
\mathcal{H}om_{\mathcal{O}_X}\left( (\mathcal{F}, \theta), (\mathcal{F}', \theta') \right)
:=
\left( \mathcal{H}om_{\mathcal{O}_X}(\mathcal{F}, \mathcal{F}'), \, \Theta \right)
\in
\Mod_G(X),
$$
is defined to be the $\mathcal{O}_X$-module $\mathcal{H}om_{\mathcal{O}_X}(\mathcal{F}, \mathcal{F}')$
together with the $G$-linearization
$$
\Theta
:
\sigma^*\mathcal{H}om_{\mathcal{O}_X}(\mathcal{F}, \mathcal{F}')
\xrightarrow{\cong}
\left(\pr_X^{G \times X}\right)^* \mathcal{H}om_{\mathcal{O}_X}(\mathcal{F}, \mathcal{F}')
$$
given by the conjugation
$
\alpha
\mapsto
\theta^{-1}
\comp
\alpha
\comp
\theta'
$
under the canonical isomorphisms\footnote{
These canonical homomorphisms are isomorphisms
because the morphisms
$\sigma$,
$\pr_X^{G\times X}
:
G \times X \rightarrow X$
are flat
(as any group scheme $G$ is flat over the base)
and $\mathcal{F}$ is finitely presented.
}
\begin{align*}
\sigma^* \mathcal{H}om_{\mathcal{O}_X}(\mathcal{F}, \mathcal{F}')
&\xrightarrow{\cong}
\mathcal{H}om_{\mathcal{O}_{G \times X}}\left(\sigma^* \mathcal{F}, \sigma^* \mathcal{F}'\right)
\quad\text{and}\quad
\\[2mm]
\left(\pr_X^{G \times X}\right)^* \mathcal{H}om_{\mathcal{O}_X}(\mathcal{F}, \mathcal{F}')
&\xrightarrow{\cong}
\mathcal{H}om_{\mathcal{O}_{G \times X}}\left(\left(\pr_X^{G \times X}\right)^*\mathcal{F}, \left(\pr_X^{G \times X}\right)^*\mathcal{F}'\right).
\end{align*}

More precisely, for each open $U \subseteq G \times X$ and a local section
$
\alpha \in \left(\sigma^*\mathcal{H}om_{\mathcal{O}_X}(\mathcal{F}, \mathcal{F}')\right)(U)
$
viewed as an $\mathcal{O}_{G\times X}|_U$-module homomorphism
$
\alpha : \sigma^*\mathcal{F}|_U \rightarrow \sigma^*\mathcal{F}'|_U
$,
the section $\Theta(U)(\alpha)$ is defined by the composition
$$
\left(\theta|_U\right)^{-1}
\comp
\alpha
\comp
\theta'|_U
:
\left(\pr_X^{G \times X}\right)^*\mathcal{F}|_U
\rightarrow
\left(\pr_X^{G \times X}\right)^*\mathcal{F}'|_U.
$$

Equipped with the tensor product and internal hom,
$\VB_G(X)$
(as well as $\Coh_G(X)$ if $X$ is locally Noetherian)
forms a \emph{closed symmetric monoidal category}.

\end{defn}
We end this paragraph by computing sections between $G$-equivariant sheaves.
\begin{prop}\label{prop:InvariantSectionsOfEquivariantInternalHom}
Let
$(\mathcal H,\Theta)
:=
\mathcal Hom_{\mathcal O_X}
\bigl((\mathcal F,\theta),(\mathcal F',\theta')\bigr)$.
For every $G$-invariant open subscheme $U\subseteq X$,
we have
$$
\mathcal H(U)^G
=
\Hom_{\Mod_G(U)}
\left(
(\mathcal F,\theta)|_U,\,
(\mathcal F',\theta')|_U
\right).
$$
\end{prop}

\begin{proof}
A section $f\in\mathcal H(U)$ is an $\mathcal O_U$-module
morphism $f:\mathcal F|_U\rightarrow\mathcal F'|_U$.
By the definition of $G$-invariant sections and of $\Theta$,
the following hold:
\begin{align*}
f\in\mathcal H(U)^G
&\quad\Leftrightarrow\quad
\Theta|_{G\times U}(\sigma_U^*f)
=
\left(\pr_U^{G\times U}\right)^*f
\\
&\quad\Leftrightarrow\quad
\theta_U^{-1}\cdot\sigma_U^*f\cdot\theta'_U
=
\left(\pr_U^{G\times U}\right)^*f
\\
&\quad\Leftrightarrow\quad
\sigma_U^*f\cdot\theta'_U
=
\theta_U\cdot
\left(\pr_U^{G\times U}\right)^*f
\end{align*}
The last equality is exactly the $G$-equivariance condition
for $f$.
\end{proof}

\subsection{The action of character group $\widehat{G}$ on $\Mod_G(X)$.} In the abelian category $\Mod_G (X)$ of $G$-equivariant sheaves on $X$ there is a canonical action of a specific abelian group, $\widehat{G}$, called the character group of $G$, which we will define as follows.
\begin{defn}
A \emph{character} of a group scheme $(G,m,e)$ is a group scheme homomorphism
$
\cha : G \longrightarrow \Gm
$
from $G$ to the multiplicative group scheme
$\Gm := \Spec\left(\field[t, t^{-1}]\right)$.
Equivalently,
a character corresponds to a global section
$
\cha
\in
\Gamma\left(G, \mathcal{O}_G^\times\right)
$
of the trivial $\Gm$-bundle $\mathcal{O}_X^\times$ on $X$
satisfying
$$
m^* \cha
=
\left(\pr_{G,1}^{G \times G}\right)^* \cha
\cdot
\left(\pr_{G,2}^{G \times G}\right)^* \cha
\in
\Gamma\left(G \times G, \mathcal{O}_{G \times G}^\times\right),
$$
where
$\pr_{G,1}^{G \times G}$,
$\pr_{G,2}^{G \times G}
:
G \times G \rightarrow G
$
are the projections onto the first and second factors, respectively,
$
\left(\pr_{G,i}^{G \times G}\right)^* \cha
\in
\Gamma\left(
G \times G,
\left(\pr_{G,i}^{G \times G}\right)^*\mathcal{O}_G^\times
\cong
\mathcal{O}_{G \times G}^\times
\right)
$
is the pullback of $\cha$
by $\pr_{G,i}^{G\times G}$
($i\in\{1,2\}$),
and $\cdot$ denotes the multiplication of elements of the ring
$\mathcal{O}_{G\times G}^\times(G\times G)$.

The set of all characters of $G$ forms an abelian group,
with the group structure induced by the multiplicative group scheme structure of $\Gm$,
denoted by
$$
\widehat{G}
:=
\operatorname{Hom}_{\mathsf{GrpSch}(\field)}\left(G, \Gm\right),
$$
and is called the \emph{character group} of $G$.

\end{defn}

\begin{defn}[Twisting by a character]

For each character
$\cha \in \widehat{G}\subset\Gamma\left(G, \mathcal{O}_G^\times\right)$,
the \emph{trivial line bundle twisted by $\cha$}
is defined by
$$
\mathcal{O}_X(\cha)
:=
(\mathcal{O}_X, \id_{\mathcal{O}_{G\times X}}(\cha))
\in \vb_G(X),
$$
where the underlying $\mathcal{O}_X$-module is $\mathcal{O}_X$
and the $G$-linearization
$$
\id_{\mathcal{O}_{G\times X}}(\cha)
:=
\left(\pr_G^{G \times X}\right)^* \cha \cdot \id_{\mathcal{O}_{G\times X}}
:
\mathcal{O}_{G\times X}
\xrightarrow{\cong}
\mathcal{O}_{G\times X}
$$
(under the canonical identifications
$
\sigma^* \mathcal{O}_X
\cong
\mathcal{O}_{G\times X}
$
and
$
\left(\pr_X^{G \times X}\right)^* \mathcal{O}_X
\cong
\mathcal{O}_{G\times X}
$)
is defined by
the multiplication of
the pullback
$$
\left(\pr_G^{G \times X}\right)^* \cha
\in
\Gamma\left(G \times X, \left(\pr_G^{G\times X}\right)^*\mathcal{O}_G^\times \cong \mathcal{O}_{G \times X}^\times\right)
$$
of $\cha$ by the canonical projection
$\pr_G^{G\times X}:G\times X\rightarrow G$.
More precisely,
we define
$$
\left(\left(\id_{\mathcal{O}_{G\times X}}(\cha)\right)(U)\right)(s)
:=
\left(\left.\left(\pr_G^{G \times X}\right)^* \cha\right|_U\right)
\cdot
s
\in
\mathcal{O}_{G\times X}(U)
$$
for each open
$U\subseteq G\times X$
and
$s\in\mathcal{O}_{G\times X}(U)$,
where
$\cdot$ denotes the multiplication of elements of the ring
$\mathcal{O}_{G\times X}(U)$.

Let $(\mathcal{F}, \theta)$ be a $G$-equivariant $\mathcal{O}_X$-module in $\Mod_G(X)$ and $\cha \in \widehat{G}$.
The \emph{twisting of $(\mathcal{F},\theta)$ by $\cha$}
is defined by
$$
(\mathcal{F},\theta)(\cha)
:=
(\mathcal{F}, \theta(\cha))
\in \Mod_G(X),
$$
where the underlying $\mathcal{O}_X$-module is $\mathcal{F}$
and the $G$-linearization
$$
\theta(\cha)
:=
\left(\pr_G^{G \times X}\right)^* \cha \cdot \theta
:
\sigma^* \mathcal{F}
\xrightarrow{\cong}
\left(\pr_X^{G \times X}\right)^* \mathcal{F}
$$
is defined by the scalar multiplication of
the pullback
$$
\left(\pr_G^{G \times X}\right)^* \cha
\in
\Gamma\left(G \times X, \left(\pr_G^{G\times X}\right)^*\mathcal{O}_G^\times \cong \mathcal{O}_{G \times X}^\times\right)
$$
on the original $G$-linearization $\theta$.
More precisely,
we define
$$
(\theta_{\cha}(U))(s)
:=
\left(\left.\left(\pr_G^{G \times X}\right)^* \cha\right|_U\right)
\cdot
(\theta(U)(s))
\in
\left(\left(\pr_X^{G \times X}\right)^*\mathcal{F}\right)(U)
$$
for each open
$U\subseteq G\times X$
and
$s\in(\sigma^*\mathcal{F})(U)$,
where
$\cdot$ denotes the scalar multiplication of elements of the $\mathcal{O}_{G\times X}(U)$-module
$\left(\left(\pr_X^{G \times X}\right)^*\mathcal{F}\right)(U)$
by elements of the ring $\mathcal{O}_{G\times X}(U)$.

Twisting by characters defines a strict $\widehat{G}$-action on $\Mod_G(X)$
$$
\phi:\widehat{G}\rightarrow\Aut_{\str}(\Mod_G(X)),
\quad
\cha\mapsto \phi_{\cha}
$$
where
$$
\phi_{\cha}:\Mod_G(X)\rightarrow\Mod_G(X),
\quad
(\mathcal{F},\theta)
\mapsto
(\mathcal{F},\theta(\cha)).
$$
This restricts to strict $\widehat{G}$-actions on
$\QCoh_G(X)$, $\Coh_G(X)$, and $\VB_G(X)$.
\end{defn}

\begin{rmk}
There is a natural isomorphism
$$
(\mathcal{F},\theta)(\cha)
\cong
(\mathcal{F},\theta)
\underset{\mathcal{O}_X}{\otimes}
\mathcal{O}_X(\cha)
\in\Mod_G(X).
$$
\end{rmk}

\subsection{Weight space decompositions} Since $\widehat{G}$ is abelian, its action on hom spaces in $\Mod_G (X)$ will induce a decomposition of hom spaces into one-dimensional irreducible $\widehat{G}$-representations, which we call customly \emph{weight space decompositions}.

\begin{defn}
Let $\Lambda$ be an abelian group and $\field$ be a commutative ring.
The \emph{group algebra} $\field[\Lambda]$ is a Hopf algebra with a basis $\{e_\lambda\}_{\lambda \in \Lambda}$,
where the comultiplication $\Delta : \field[\Lambda] \rightarrow \field[\Lambda] \underset{\field}{\otimes} \field[\Lambda]$ and the counit $\varepsilon : \field[\Lambda] \rightarrow \field$ are defined by
$$
\Delta(e_\lambda) = e_\lambda \otimes e_\lambda \quad \text{and} \quad \varepsilon(e_\lambda) = 1
\quad\textnormal{for all $\lambda \in \Lambda$}.
$$

A \emph{right comodule} over $\field[\Lambda]$ is a $\field$-module $V$ equipped with a $\field$-linear map $\rho : V \rightarrow V \underset{\field}{\otimes} \field[\Lambda]$ satisfying the coassociativity and counit axioms:
$$
\rho \comp (\rho \otimes \id_{\field[\Lambda]}) = \rho \comp (\id_V \otimes \Delta)
\quad\text{and}\quad
\rho \comp (\id_V \otimes \varepsilon) = \id_V.
$$
\end{defn}

\begin{prop}[Weight space decomposition of $\field[\Lambda\text{$]$}$-comodules]
\label{prop:WeightSpaceDecompositionForComodule}
A $\field$-module $V$ has a right comodule structure over $\field[\Lambda]$ if and only if $V$ is $\Lambda$-graded, i.e., it decomposes into a direct sum of weight spaces
$$
V = \bigoplus_{\lambda \in \Lambda} V_\lambda,
$$
where $V_\lambda = \{v \in V \mid \rho(v) = v \otimes e_\lambda\}$.
\end{prop}

\begin{proof}
For any $v \in V$, its image under $\rho$ can be uniquely written as a finite sum
$$
\rho(v) = \sum_{\lambda \in \Lambda} P_\lambda(v) \otimes e_\lambda,
$$
where $P_\lambda : V \rightarrow V$ are $\field$-linear maps.
Applying the counit axiom $\rho \comp (\id_V \otimes \varepsilon) = \id_V$ to $v$ yields
$$
v = (\id_V \otimes \varepsilon)(\rho(v)) = (\id_V \otimes \varepsilon)\left(\sum_{\lambda \in \Lambda} P_\lambda(v) \otimes e_\lambda\right) = \sum_{\lambda \in \Lambda} P_\lambda(v) \cdot 1 = \sum_{\lambda \in \Lambda} P_\lambda(v).
$$
This shows that the sum of the images of $P_\lambda$ covers $V$.
Next, applying the coassociativity axiom
$
\rho \comp (\rho \otimes \id_{\field[\Lambda]}) = \rho \comp (\id_V \otimes \Delta)
$
to $v$ gives
$$
\sum_{\lambda, \mu \in \Lambda} P_\mu(P_\lambda(v)) \otimes e_\mu \otimes e_\lambda = \sum_{\lambda \in \Lambda} P_\lambda(v) \otimes e_\lambda \otimes e_\lambda.
$$
Comparing the coefficients of the basis elements $e_\mu \otimes e_\lambda$ on both sides, we obtain $P_\mu(P_\lambda(v)) = \delta_{\lambda, \mu} P_\lambda(v)$.
Using the left-to-right composition convention, this implies $P_\lambda \comp P_\mu = \delta_{\lambda, \mu} P_\lambda$.
This implies that the maps $\{P_\lambda\}_{\lambda \in \Lambda}$ form a family of mutually orthogonal projection operators summing to the identity.
Thus, $V$ completely decomposes into the direct sum of the images $\mathrm{Im}(P_\lambda)$.
It is straightforward to verify that $\mathrm{Im}(P_\lambda) = V_\lambda$, concluding the proof.
\end{proof}

\begin{defn}[Diagonalizable group scheme]
An affine group scheme $G$ over $\field$ is called \emph{diagonalizable} if there exists an abelian group $\Lambda$ such that
$G\cong\Spec(\field[\Lambda])$
as a group scheme. In this case, the group $\cha$ is canonically identified with the character group $\widehat{G}$ of $G$, and its coordinate ring $\mathcal{O}_G(G)$ is isomorphic to the group algebra $\field[\widehat{G}]$ as a Hopf algebra.
\end{defn}

\begin{lemma}[Weight space decomposition of global sections]
\label{lem:GlobalSectionDecomposition}
Let $G$ be a diagonalizable group scheme acting on a scheme $X$, and let $(\mathcal{F}, \theta) \in \Mod_G(X)$ be a $G$-equivariant $\mathcal{O}_X$-module.
Then the space of global sections $V := \mathcal{F}(X)$($=\Gamma(X, \mathcal{F})$) naturally admits a right comodule structure over $\mathcal{O}_G(G)$, and therefore admits a canonical weight space decomposition
$$
\mathcal{F}(X) = \bigoplus_{\cha \in \widehat{G}}
\mathcal{F}(X)_{\cha}.
$$
where
$
\mathcal{F}(X)_{\cha}
:=
\left\{
s \in \mathcal{F}(X)
\left|
(\theta(G\times X))(\sigma^* s)
=
\left(\pr_G^{G \times X}\right)^*\cha
\cdot
\left(\pr_X^{G \times X}\right)^* s
\right.\right\}.
$

\end{lemma}

\begin{proof}
The $G$-linearization $\theta : \sigma^* \mathcal{F} \xrightarrow{\cong} \left(\pr_X^{G \times X}\right)^* \mathcal{F}$ induces a natural map on the global sections.
For any section $s \in V = \mathcal{F}(X)$, its pullback along the action map $\sigma:G\times X \rightarrow X$ yields a section $\sigma^* s \in (\sigma^*\mathcal{F})(G\times X)$.
Applying the isomorphism $\theta$ to this pullback gives
$$
(\theta(G\times X))(\sigma^* s)
\in
\left(\left(\pr_X^{G \times X}\right)^* \mathcal{F}\right)(G\times X)
\cong
\mathcal{F}(X) \underset{\field}{\otimes} \mathcal{O}_G(G)
=
V \underset{\field}{\otimes} \mathcal{O}_G(G).
$$
This defines a $\field$-linear map
$$
\rho
:
V \rightarrow V \underset{\field}{\otimes} \mathcal{O}_G(G),
\quad
s
\mapsto
(\theta(G\times X))(\sigma^* s).
$$

The cocycle condition and the unit condition of $\theta$
translates precisely to the coassociativity axiom
and the counit axiom of $\rho$,
respectively.
Therefore, $\rho$ defines a right comodule structure on $V$ over $\mathcal{O}_G(G) \cong \field[\widehat{G}]$.
By Proposition \ref{prop:WeightSpaceDecompositionForComodule},
this comodule structure completely decomposes $V$ into the direct sum of its weight spaces indexed by the character group $\widehat{G}$.

Under the canonical identification
$
\left(\left(\pr_X^{G \times X}\right)^* \mathcal{F}\right)(G\times X)
\cong
\mathcal{F}(X) \underset{\field}{\otimes} \mathcal{O}_G(G),
$
the pure tensor $s \otimes e_{\cha} \in V \underset{\field}{\otimes} \mathcal{O}_G(G)$ corresponds to
$
\left(\pr_G^{G \times X}\right)^*\cha
\cdot
\left(\pr_X^{G \times X}\right)^* s
\in
\left(\left(\pr_X^{G \times X}\right)^* \mathcal{F}\right)(G\times X),
$
where
$\cdot$ denotes the scalar multiplication of elements of the $\mathcal{O}_{G\times X}(G\times U)$-module
$\left(\left(\pr_X^{G \times X}\right)^*\mathcal{F}\right)(G\times U)$
by elements of the ring $\mathcal{O}_{G\times X}(G\times U)$.
Thus, the condition $\rho(s) = s \otimes e_{\cha}$ translates into the condition
$$
(\theta(G\times X))(\sigma^* s) = \left(\pr_G^{G \times X}\right)^*\cha \cdot \left(\pr_X^{G \times X}\right)^* s,
$$
establishing the formula for the weight space $\mathcal{F}(X)_{\cha}$.
\end{proof}

\begin{prop}[Weight space decomposition of hom spaces]
\label{prop:WeightSpaceDecompositionOfHomSpace}
Let $G$ be a diagonalizable group scheme acting on a scheme $X$.
For any
$(\mathcal{F}, \theta)$,
$(\mathcal{F}', \theta') \in \Mod_G(X)$
such that
$\mathcal{F}$
is finitely presented,
the hom space from $\mathcal{F}$ to $\mathcal{F}'$
in $\Mod(X)$
admits a canonical weight space decomposition
$$
\operatorname{Hom}_{\Mod(X)}(\mathcal{F}, \mathcal{F}') 
= 
\bigoplus_{\cha \in \widehat{G}} \operatorname{Hom}_{\Mod_G(X)}\left( (\mathcal{F}, \theta)(\cha), (\mathcal{F}', \theta') \right).
$$
\end{prop}
\begin{proof}
Since $\mathcal{F}$ is finitely presented,
by Definition \ref{defn:InternalHomOfEquivariantModules}, 
the internal hom $\mathcal{H} := \mathcal{H}om_{\mathcal{O}_X}(\mathcal{F}, \mathcal{F}')$ is naturally equipped with a $G$-linearization $\Theta$ given by $\alpha \mapsto \theta^{-1} \comp \alpha \comp \theta'$, making $(\mathcal{H}, \Theta) \in \Mod_G(X)$.

By Lemma \ref{lem:GlobalSectionDecomposition},
the space of global sections
$$
\operatorname{Hom}_{\Mod(X)}(\mathcal{F}, \mathcal{F}')
=
\mathcal{H}(X)
$$
admits a canonical weight space decomposition
$
\mathcal{H}(X)
=
\bigoplus_{\cha \in \widehat{G}} \mathcal{H}(X)_{\cha},
$
where
the following hold
for $f\in\mathcal{H}(X)$ and $\cha\in\widehat{G}$:
\begin{align*}
f\in \mathcal{H}(X)_{\cha}
\quad
&\Leftrightarrow
\quad
(\Theta(G\times X))(\sigma^* f)
=
\left(\pr_G^{G \times X}\right)^*\cha
\cdot
\left(\pr_X^{G \times X}\right)^* f
\\[2mm]
&\Leftrightarrow
\quad
\theta^{-1} \comp \sigma^* f \comp \theta' = \left(\pr_G^{G \times X}\right)^*\cha \cdot \left(\pr_X^{G \times X}\right)^* f
\\[2mm]
&\Leftrightarrow
\quad
\sigma^* f \comp \theta' = \theta \comp \left( \left(\pr_G^{G \times X}\right)^*\cha \cdot \left(\pr_X^{G \times X}\right)^* f \right)
\\[2mm]
&\Leftrightarrow
\quad
\sigma^* f \comp \theta' = \left( \left(\pr_G^{G \times X}\right)^*\cha \cdot \theta \right) \comp \left(\pr_X^{G \times X}\right)^* f\\[2mm]
&(\because \left(\pr_G^{G \times X}\right)^*\cha \cdot\text{ is a scalar multiplication.})
\\[2mm]
&\Leftrightarrow
\quad
\sigma^* f \comp \theta' = \theta(\cha) \comp \left(\pr_X^{G \times X}\right)^* f
\quad(\because \text{definition of $\theta(\cha)$})
\end{align*}

This equality is precisely the condition for $f$ to be a $G$-equivariant morphism from
$(\mathcal{F}, \theta(\cha))$
to
$(\mathcal{F}', \theta')$
(Definition \ref{defn:EquivariantSheaves}). 
Therefore, $\mathcal{H}(X)_{\cha} = \operatorname{Hom}_{\Mod_G(X)}\left( (\mathcal{F}, \theta)(\cha), (\mathcal{F}', \theta') \right)$, completing the canonical weight space decomposition.
\end{proof}
Finally, as a Corollary to the computation, we get an important intermediate step toward our proof of the main theorem.
\begin{cor}[Embedding of orbit categories]
\label{thm:EmbeddingOfTheOrbitCategory}
Let $G$ be a diagonalizable group scheme acting on a scheme $X$.
Then the forgetful functor $\Mod_G(X) \rightarrow \Mod(X)$ sending objects $(\mathcal{F}, \theta) \mapsto \mathcal{F}$ induces a fully faithful embedding of the orbit categories\footnote{
The same holds true for $\Coh_G(X)$ and $\Coh(X)$ when $X$ is locally Noetherian.
}:
$$
\QCoh_G(X)/\widehat{G}
\hookrightarrow
\QCoh(X).
$$
\end{cor}

\begin{proof}
By definition, the hom space in the orbit category is given by
$$
\operatorname{Hom}_{\QCoh_G(X)/\widehat{G}}\left( (\mathcal{F}, \theta), (\mathcal{F}', \theta') \right)
=
\bigoplus_{\cha \in \widehat{G}} \operatorname{Hom}_{\QCoh_G(X)}\left( (\mathcal{F}, \theta)(\cha), (\mathcal{F}', \theta') \right).
$$
By Proposition \ref{prop:WeightSpaceDecompositionOfHomSpace}, this direct sum exactly coincides with the hom space $\operatorname{Hom}_{\QCoh(X)}(\mathcal{F}, \mathcal{F}')$ in $\QCoh(X)\subseteq\Mod(X)$.
Therefore, the functor $\QCoh_G(X)/\widehat{G} \hookrightarrow \QCoh(X)$ defined by $(\mathcal{F}, \theta) \mapsto \mathcal{F}$ on objects and by the identity map on morphisms is fully faithful.
\end{proof}

\subsection{Restriction to subgroup schemes.}

Let $G$ be a group scheme,
and $H\subseteq G$ a closed subgroup scheme
with the inclusion morphism $i_H : H \hookrightarrow G$.
A $G$-action
$\sigma_G : G \times X \rightarrow X$
on a scheme $X$ naturally induces a $H$-action on $X$ defined by
$$
\sigma_H
:=
(i_H \times \id_X)
\comp
\sigma
:
H\times X\rightarrow X.
$$

\begin{defn}[Restriction of equivariant sheaves]\label{defn:RestrictionOfEquivariantSheaves}
For any $G$-equivariant $\mathcal{O}_X$-module $(\mathcal{F}, \theta) \in \Mod_G(X)$, its \emph{restriction} to $H$ is defined to be the $H$-equivariant $\mathcal{O}_X$-module
$$
(\mathcal{F}, \theta)|_H := (\mathcal{F}, \theta|_H),
$$
where the restricted $H$-linearization
is defined by the pullback of the $G$-linearization
$$
\theta|_H 
:= 
(i_H \times \id_X)^* \theta 
: 
\sigma_H^* \mathcal{F} 
\xrightarrow{\cong} 
\left(\pr_X^{H \times X}\right)^* \mathcal{F}
$$

This operation on objects, together with the identity map on morphisms, naturally defines a functor
$$
\Res_H^G : \Mod_G(X) \rightarrow \Mod_H(X),
$$
called the \emph{restriction functor}.
This also induces restriction functors on the corresponding full subcategories:
$$
\QCoh_G(X) \rightarrow \QCoh_H(X),
\quad
\Coh_G(X) \rightarrow \Coh_H(X),
\quad\text{and}\quad
\VB_G(X) \rightarrow \VB_H(X).
$$
\end{defn}

The following lemma is standard:

\begin{lemma}\cite[Exercise 2.8]{Waterhouse}
\label{lem:DiagonalizableSurjective}
Let $G$ be a diagonalizable group scheme
over a field $\field$,
and let $H \subseteq G$ be a closed subgroup scheme.
Then $H$ is also diagonalizable, and the restriction map on character groups
$
i_H^* : \widehat{G} \rightarrow \widehat{H}
$
is surjective with $\ker(i_H^*) = H^\perp$,
where
$$
H^\perp
:=
\{ \cha \in \widehat{G} \mid \cha(h) = 1 \text{ for all } h \in H \}
\subseteq
\widehat{G}
$$
is the annihilator of $H$ in $\widehat{G}$.
\end{lemma}

\begin{proof}
Since $H$ is a closed subgroup scheme of
a diagonalizable group scheme
$G\cong \Spec(\field[\widehat{G}])$,
there is a Hopf ideal $I$ of $\field[\widehat{G}]$
such that
$H\cong\Spec(\field[\widehat{G}]/I)$
and
$\mathcal{O}_H(H)\cong\field[\widehat{G}]/I$.

We claim that
$$
I=\left<\cha-1|\cha\in H^\perp\right>
\subseteq
\field[\widehat{G}].
$$

For any
$\cha \in H^\perp\subseteq\field[\widehat{G}]\cong\mathcal{O}_G(G)$,
we have
$\cha|_H=1_{\mathcal{O}_H(H)}\in\mathcal{O}_H(H)$,
and hence
$(\cha-1)|_H=0\in \mathcal{O}_H(H)$.
This implies
$\langle \cha - 1 \mid \cha \in H^\perp \rangle \subseteq I$.

Conversely,
suppose
$f = \sum_{\cha \in \widehat{G}} c_{\cha} \cha \in I\subseteq\field[\widehat{G}]\cong\mathcal{O}_G(G)$.
Grouping the terms by their restricted images in $\mathcal{O}_H(H)$, we evaluate $f|_H$:
$$
0 = f|_H = \sum_{\mu} \left( \sum_{\cha|_H = \mu} c_{\cha} \right) \mu \quad \text{in } \mathcal{O}_H(H),
$$
where $\mu$ ranges over the distinct non-zero images of characters in $\mathcal{O}_H(H)$. Since each $\mu$ is a group-like element in the Hopf algebra $\mathcal{O}_H(H)$, these distinct elements $\mu$ are linearly independent over $\field$
\cite[\S2.2]{Waterhouse}.
Consequently, the coefficient of each $\mu$ must vanish; that is, $\sum_{\cha|_H = \mu} c_{\cha} = 0$ for all such $\mu$. 

For each $\mu$ appearing in the sum, fix a representative $\cha_0 \in \widehat{G}$ such that $\cha_0|_H = \mu$. Any other character $\cha$ restricting to $\mu$ satisfies $(\cha \cha_0^{-1})|_H = 1$, meaning $\cha \cha_0^{-1} \in H^\perp$. Thus, we can write $\cha = h\cha_0$ for some $h \in H^\perp$. We can then rewrite the sub-sum corresponding to $\mu$ as follows:
$$
\sum_{h \in H^\perp} c_{h\cha_0} (h\cha_0) = \sum_{h \in H^\perp} c_{h\cha_0} (h - 1)\cha_0 + \left( \sum_{h \in H^\perp} c_{h\cha_0} \right) \cha_0.
$$
Since the sum of the coefficients $\sum c_{h\cha_0}$ is exactly zero, the second term vanishes. This leaves only a linear combination of elements of the form $h-1$ for $h \in H^\perp$. Summing over all such $\mu$, we conclude that $f \in \langle \cha - 1 \mid \cha \in H^\perp \rangle$.

Thus, $I = \langle \cha - 1 \mid \cha \in H^\perp \rangle$. By the first isomorphism theorem for rings, this yields the natural isomorphism $\mathcal{O}_H(H) \cong \field[\widehat{G}]/I \cong \field[\widehat{G}/H^\perp]$,
concluding that $H$ is diagonalizable
with
$\widehat{G}/H^\perp\cong\widehat{H}$.
\end{proof}

\begin{prop}[$H$-invariants of hom spaces]
\label{prop:HInvariantsOfHomSpace}
Let $G$ be a diagonalizable group scheme acting on a scheme $X$, and $H\subseteq G$ a closed subgroup scheme.
For any
$(\mathcal{F}, \theta)$,
$(\mathcal{F}', \theta') \in \Mod_G(X)$
such that
$\mathcal{F}$
is finitely presented,
the hom space
from
$(\mathcal{F}, \theta|_H)$
to
$(\mathcal{F}', \theta'|_H)$
in $\Mod_H(X)$
admits a canonical decomposition
$$
\operatorname{Hom}_{\Mod_H(X)}\left( (\mathcal{F}, \theta|_H), (\mathcal{F}', \theta'|_H) \right)
=
\bigoplus_{\cha \in H^\perp} \operatorname{Hom}_{\Mod_G(X)}\left( (\mathcal{F}, \theta)(\cha), (\mathcal{F}', \theta') \right).
$$
\end{prop}

\begin{proof}
Note that $H$ is diagonalizable
by Lemma \ref{lem:DiagonalizableSurjective}.
Therefore,
by Proposition \ref{prop:WeightSpaceDecompositionOfHomSpace},
we have two decompositions
\begin{align*}
\operatorname{Hom}_{\Mod(X)}(\mathcal{F}, \mathcal{F}') 
&=
\bigoplus_{\cha \in \widehat{G}} \operatorname{Hom}_{\Mod_G(X)}\left( (\mathcal{F}, \theta)(\cha), (\mathcal{F}', \theta') \right)
\\
&=
\bigoplus_{\chaa \in \widehat{H}} \operatorname{Hom}_{\Mod_H(X)}\left( (\mathcal{F}, \theta|_H)(\chaa), (\mathcal{F}', \theta'|_H) \right).
\end{align*}

We claim that
the following holds for each $\cha\in\widehat{G}$:
$$
\operatorname{Hom}_{\Mod_G(X)}\left( (\mathcal{F}, \theta)(\cha), (\mathcal{F}', \theta') \right)
\subseteq
\operatorname{Hom}_{\Mod_H(X)}\left( (\mathcal{F}, \theta|_H)(i_H^*\cha), (\mathcal{F}', \theta'|_H) \right).
$$
Indeed,
$
f\in\operatorname{Hom}_{\Mod_G(X)}\left( (\mathcal{F}, \theta)(\cha), (\mathcal{F}', \theta') \right)
$
satisfies
$$
\sigma_G^* f\comp \theta'
=
\left(\left(\pr_G^{G\times X}\right)^*\cha \cdot \theta\right)
\comp
\left(\pr_X^{G\times X}\right)^* f
$$
Taking $(i_H^*\times \id_X)^*$
on both sides yields
$$
\sigma_H^*f\comp\left(\theta'|_H\right)
=
\left(
\left(\pr_H^{H\times X}\right)^*
i_H^*\cha
\cdot
\theta|_H
\right)
\comp
\left(\pr_X^{H\times X}\right)^*f,
$$
where
$
\theta|_H
(i_H^*\cha)
=
\left(
\left(\pr_H^{H\times X}\right)^*
i_H^*\cha
\cdot
\theta|_H
\right).
$
This implies
$
f\in
\operatorname{Hom}_{\Mod_H(X)}\left( (\mathcal{F}, \theta|_H)(i_H^*\cha), (\mathcal{F}', \theta'|_H) \right),
$
establishing the claim.

Therefore,
we have the following sub-decomposition
for each $\chaa\in\widehat{H}$:
$$
\operatorname{Hom}_{\Mod_H(X)}\left( (\mathcal{F}, \theta|_H)(\chaa), (\mathcal{F}', \theta'|_H) \right)
=
\bigoplus_{\substack{\cha \in \widehat{G} \\[1mm] i_H^*\cha = \chaa}} \operatorname{Hom}_{\Mod_G(X)}\left( (\mathcal{F}, \theta)(\cha), (\mathcal{F}', \theta') \right).
$$
Now the proposition is a special case when $\chaa=1$.
\end{proof}
Similar to corollary \ref{thm:EmbeddingOfTheOrbitCategory}, we have an embedding of orbit categories to restrict $G$-equivariant vector bundles to $H$-equivariant ones.
\begin{crl}\label{crl:orbit-category-GH-restriction}
    Let $G$ be a group scheme and $H\subseteq G$ be a subgroup scheme of $G$. Then there is an embedding of abelian categories \[
        \QCoh_G (X)/H^{\perp}\hookrightarrow\QCoh_H (X).
    \]
\end{crl}
\begin{proof}
By definition, the hom space in the orbit category is given by
$$
\operatorname{Hom}_{\QCoh_G(X)/H^\perp}\left( (\mathcal{F}, \theta), (\mathcal{F}', \theta') \right)
=
\bigoplus_{\cha \in H^\perp} \operatorname{Hom}_{\QCoh_G(X)}\left( (\mathcal{F}, \theta)(\cha), (\mathcal{F}', \theta') \right).
$$
By Proposition \ref{prop:HInvariantsOfHomSpace},
this direct sum exactly coincides with the hom space
\[\operatorname{Hom}_{\QCoh_H(X)}
\left((\mathcal{F}, \theta|_H), (\mathcal{F}', \theta'|_H) \right)\]
in
$\QCoh_H(X)\subseteq\Mod_H(X)$.
Therefore, the functor $\QCoh_G(X)/H^\perp \hookrightarrow \QCoh_H(X)$ defined by $(\mathcal{F}, \theta) \mapsto (\mathcal{F},\theta|_H)$
on objects and by the identity map on morphisms is fully faithful.
\end{proof}

\subsection{Product of group scheme actions.}

Suppose that group schemes $G$ and $H$
act on schemes $X$ and $Y$
with action homomorphisms
$\sigma_X : G \times X \rightarrow X$
and
$\sigma_Y : H \times Y \rightarrow Y$,
respectively.
They naturally induce an action of
the product group scheme $G\times H$
on the product scheme $X\times Y$
given by the action morphism
$$
\sigma_{X\times Y}
:=
\sigma_X\times \sigma_Y
:
(G\times H)\times (X\times Y)
\rightarrow
X\times Y.
$$
\begin{defn}[External tensor product of equivariant sheaves]
\label{defn:ExternalTensorProductOfEquivariantModules}
Let $(\mathcal{F},\theta)\in\Mod_G(X)$ and
$(\mathcal{G},\eta)\in\Mod_H(Y)$.
The \emph{external tensor product} of
$(\mathcal{F},\theta)$ and $(\mathcal{G},\eta)$
$$
(\mathcal{F},\theta)\boxtimes(\mathcal{G},\eta)
:=
\left(
\mathcal{F}\boxtimes\mathcal{G},
\
\theta\boxtimes\eta
\right)
\in\Mod_{G\times H}(X\times Y)
$$
is defined as the
$\mathcal{O}_{X\times Y}$-module
$
\mathcal{F}\boxtimes\mathcal{G}
:=
\pi_X^*\mathcal{F}
\underset{\mathcal{O}_{X\times Y}}{\otimes}
\pi_Y^*\mathcal{G}
$
together with the
$G\times H$-linearization
$$
\theta\boxtimes\eta
:=
\tilde{\pi}_{G\times X}^*\theta\otimes\tilde{\pi}_{H\times Y}^*\eta
:
\sigma_{X\times Y}^*(\mathcal{F}\boxtimes\mathcal{G})
\xrightarrow{\cong}
\left(\pr_{X\times Y}^{(G\times H)\times(X\times Y)}\right)^*
(\mathcal{F}\boxtimes\mathcal{G})
$$
where
\begin{align*}
\tilde{\pi}_{G\times X}^*\theta
&:
\tilde{\pi}_{G\times X}^*(\sigma_X^*\mathcal{F})
\xrightarrow{\cong}
\tilde{\pi}_{G\times X}^*\left(\left(\pr_X^{G\times X}\right)^{\ast}\mathcal{F}\right),
\\[2mm]
\tilde{\pi}_{H\times Y}^*\eta
&:
\tilde{\pi}_{H\times Y}^*(\sigma_Y^*\mathcal{G})
\xrightarrow{\cong}
\tilde{\pi}_{H\times Y}^*\left(\left(\pr_Y^{H\times Y}\right)^*\mathcal{G}\right)
\end{align*}
are pull-backs of $\theta$ and $\eta$
by the natural projections
\begin{align*}
\tilde{\pi}_{G\times X}
&:=
\pr_{G\times X}^{(G\times H)\times(X\times Y)}
:
(G\times H)\times(X\times Y)
\rightarrow
G\times X,
\\[2mm]
\tilde{\pi}_{H\times Y}
&:=
\pr_{H\times Y}^{(G\times H)\times(X\times Y)}
:
(G\times H)\times(X\times Y)
\rightarrow
H\times Y,
\end{align*}
respectively,
and
the tensor product
$
\tilde{\pi}_{G\times X}^*\theta\otimes\tilde{\pi}_{H\times Y}^*\eta
$
is given under the canonical isomorphisms
\begin{align*}
\sigma_{X\times Y}^*(\mathcal{F}\boxtimes\mathcal{G})
&=
\sigma_{X\times Y}^*
\Big(
\pi_X^*\mathcal{F}
\underset{\mathcal{O}_{X\times Y}}{\otimes}
\pi_Y^*\mathcal{G}
\Big)
\\[2mm]
&\cong
(\sigma_{X\times Y}\comp\pi_X)^*\mathcal{F}
\underset{\mathcal{O}_{(G\times H)\times(X\times Y)}}{\otimes}
(\sigma_{X\times Y}\comp\pi_Y)^*\mathcal{G}
\\[2mm]
&=
(\tilde{\pi}_{G\times X}\comp \sigma_X)^*\mathcal{F}
\underset{\mathcal{O}_{(G\times H)\times(X\times Y)}}{\otimes}
(\tilde{\pi}_{H\times Y}\comp \sigma_Y)^*\mathcal{G}
\\[2mm]
&\cong
\tilde{\pi}_{G\times X}^*\left(\sigma_X^*\mathcal{F}\right)
\underset{\mathcal{O}_{(G\times H)\times(X\times Y)}}{\otimes}
\tilde{\pi}_{H\times Y}^*\left(\sigma_Y^*\mathcal{G}\right)
\end{align*}
and
\begin{align*}
\left(\pr_{X\times Y}^{(G\times H)\times(X\times Y)}\right)^*
(\mathcal{F}\boxtimes\mathcal{G})
&=
\left(\pr_{X\times Y}^{(G\times H)\times(X\times Y)}\right)^*
\Big(
\pi_X^*\mathcal{F}
\underset{\mathcal{O}_{X\times Y}}{\otimes}
\pi_Y^*\mathcal{G}
\Big)
\\
&\cong
\left(\pr_{X\times Y}^{(G\times H)\times(X\times Y)}\comp\pi_X\right)^*\mathcal{F}
\underset{\mathcal{O}_{(G\times H)\times(X\times Y)}}{\otimes}
\left(\pr_{X\times Y}^{(G\times H)\times(X\times Y)}\comp\pi_Y\right)^*\mathcal{G}
\\
&\cong
\left(\tilde{\pi}_{G\times X}\comp\pr_X^{G\times X}\right)^*\mathcal{F}
\underset{\mathcal{O}_{(G\times H)\times(X\times Y)}}{\otimes}
\left(\tilde{\pi}_{H\times Y}\comp\pr_Y^{H\times Y}\right)^*\mathcal{G}
\\
&\cong
\tilde{\pi}_{G\times X}^*\left(\left(\pr_X^{G\times X}\right)^*\mathcal{F}\right)
\underset{\mathcal{O}_{(G\times H)\times(X\times Y)}}{\otimes}
\tilde{\pi}_{H\times Y}^*\left(\left(\pr_Y^{H\times Y}\right)^*\mathcal{G}\right).
\end{align*}

For morphisms
$
f
\in
\Hom_{\Mod_G(X)}
((\mathcal{F},\theta),(\mathcal{F}',\theta'))
$
and
$
g
\in
\Hom_{\Mod_H(Y)}
((\mathcal{G},\eta),(\mathcal{G}',\eta')),
$
define
$$
f\boxtimes g
:=
\pi_X^*f
\otimes
\pi_Y^*g
\in
\Hom_{\Mod_{G\times H}(X\times Y)}
\big(
\left(
\mathcal{F}\boxtimes\mathcal{G},
\
\theta\boxtimes\eta
\right),
\left(
\mathcal{F}'\boxtimes\mathcal{G}',
\
\theta'\boxtimes\eta'
\right)
\big).
$$


This defines a bifunctor
$$
\boxtimes:
\Mod_G(X)\times\Mod_H(Y)
\rightarrow
\Mod_{G\times H}(X\times Y),
$$
which restricts to equivariant (quasi-)coherent sheaves and
equivariant vector bundles.
\end{defn}

\begin{prop}\label{prop:BoxTensorFunctorFullyFaithful}
The functor
$$
\boxtimes:
\Mod_G(X)\otimes\Mod_H(Y)
\hookrightarrow
\Mod_{G\times H}(X\times Y)
$$
is fully faithful,
that is,
\begin{align*}
\boxtimes:
\Hom_{\Mod_G(X)}&
\big(
(\mathcal{F},\theta),
(\mathcal{F}',\theta')
\big)
\
\otimes
\
\Hom_{\Mod_H(Y)}
\big(
(\mathcal{G},\eta),
(\mathcal{G}',\eta')
\big)
\\ 
&\xrightarrow{\cong}
\Hom_{\Mod_{G\times H}(X\times Y)}
\big(
(\mathcal{F},\theta)\boxtimes(\mathcal{G},\eta),
(\mathcal{F}',\theta')\boxtimes(\mathcal{G}',\eta')
\big)
\end{align*}
\end{prop}

\section{Twisted Künneth formula\label{sec:kunneth-formula}}
After introducing equivariant sheaves, it remains for us to establish the Künneth formula, which is what we essentially used in the proof of theorem \ref{thm:twisted-tensor-P1-bundle}. We start with establishing the dg enhancement of the triangulated category of quasi-coherent sheaves, which contains dg enhancements of bounded derived categories of coherent sheaves and perfect complexes as dg full subcategories.
\subsection{Derived dg-categories.}
Let $X$ be a scheme over $\field$, and let 
$\mathfrak{U}:=(U_i)_{i\in I}$
be an open cover of $X$.
For any sheaf $\mathcal{H}$
of abelian groups
on $X$,
recall that
the \emph{\v{C}ech cochain complex}
$
\check{C}^{\bullet}(\mathfrak{U},\mathcal{H})
$
of $\mathcal{H}$
with respect to $\mathfrak{U}$ is given by
collection of
abelian groups
$$
\check{C}^p(\mathfrak{U},\mathcal{H})
:=
\prod_{i_0<\cdots<i_p}
\mathcal{H}\left(U_{i_0\cdots i_p}\right)
\quad\text{for }
p\in\mathbb{Z}_{\ge0},
$$
where
$
U_{i_0\cdots i_p}
:=
U_{i_0}\cap\cdots\cap U_{i_p},
$
with a differential
$
\delta:
\check{C}^p(\mathfrak{U},\mathcal{H})
\rightarrow
\check{C}^{p+1}(\mathfrak{U},\mathcal{H})
$
defined by
$$
\left(\delta\alpha\right)_{i_0\cdots i_{p+1}}
:=
\sum_{j=0}^{p+1}
(-1)^j
\left.\alpha_{i_0\cdots\hat{i_j}\cdots i_{p+1}}\right|_{U_{i_0\cdots i_{p+1}}}
\in
\mathcal{H}(U_{i_0\cdots i_{p+1}})
$$
for
$i_0<\cdots<i_{p+1}$,
where we denote by
$
\alpha
=
\left(
\alpha_{i_0\cdots i_p}
\in\mathcal{H}(U_{i_0\cdots i_p})
\right)_{i_0<\cdots<i_p}.
$
\begin{dfn}\label{dfn:dg-enhancement-Coh}
Let $X$ be a scheme admitting an affine open cover
$\mathfrak{U}:=(U_i)_{i\in I}$.
The
\emph{coherent derived dg-category of $X$}
with respect to
$\mathfrak{U}$
is a dg-category
$$
\Coh(X;\mathfrak{U})_{\dg}
$$
defined as follows:
\begin{itemize}
\item
\emph{Objects} are given by bounded cochain complexes
$$
\mathcal{F}^{\bullet}
=
\left(
\bigoplus_{\objectdegree\in\mathbb{Z}}
\mathcal{F}^{\objectdegree},
d_{\mathcal{F}}:=\left(d_{\mathcal{F}}^{\objectdegree}:\mathcal{F}^{\objectdegree}\rightarrow \mathcal{F}^{\objectdegree+1}\right)_{\objectdegree\in\mathbb{Z}}\right)
$$
of coherent $\mathcal{O}_X$-modules
$
\mathcal{F}^{\objectdegree}
\in
\Coh(X).
$
That is,
the differential
$d_{\mathcal{F}}=(d_{\mathcal{F}}^{\objectdegree})_{\objectdegree\in\mathbb{Z}}$,
where
$d_{\mathcal{F}}^{\objectdegree}
\in
\Hom_{\Coh(X)}(\mathcal{F}^{\objectdegree},\mathcal{F}^{\objectdegree+1}),
$
satisfies
$d_{\mathcal{F}}^{\objectdegree}\comp d_{\mathcal{F}}^{\objectdegree+1}=0$
for each $\objectdegree\in\mathbb{Z}$.

\item
Given two objects
$\mathcal{F}$ and $\mathcal{F}'$,
we first define a bounded cochain complex
of coherent $\mathcal{O}_X$-modules
given by
$$
\mathcal{H}om_{\mathcal{O}_X}^q\left(\mathcal{F},\mathcal{F}'\right)
:=
\bigoplus_{\objectdegree\in\mathbb{Z}}
\mathcal{H}om_{\mathcal{O}_X}\left(\mathcal{F}^{\objectdegree},\mathcal{F}'^{\objectdegree+q}\right)
\quad\text{for }q\in\mathbb{Z},
$$
with
a differential
$
d_{\mathcal{H}om}:
\mathcal{H}om_{\mathcal{O}_X}^q(\mathcal{F},\mathcal{F}')
\rightarrow
\mathcal{H}om_{\mathcal{O}_X}^{q+1}(\mathcal{F},\mathcal{F}')
$
defined by
$$
d_{\mathcal{H}om}f:=d_{\mathcal{F}}\comp f - (-1)^q f\comp d_{\mathcal{F}'}.
$$

The \emph{dg-morphism space}
from $\mathcal{F}$ to $\mathcal{F}'$
is a cochain complex of $\field$-modules
whose degree $k$-component is given by
\begin{align*}
\hom_{\Coh(X;\mathfrak{U})_{\dg}}^k
\left(
\mathcal{F},\mathcal{F}'
\right)
:&=
\bigoplus_{p+q=k}
\check{C}^p\big(\mathfrak{U},\mathcal{H}om_{\mathcal{O}_X}^q
(\mathcal{F},\mathcal{F}')
\big),
\\[2mm]
&=
\bigoplus_{p+q=k}
\prod_{i_0<\cdots<i_p}\bigoplus_{\objectdegree\in\mathbb{Z}}
\Hom_{\Coh(U_{i_0\cdots i_p})}\left(\mathcal{F}^{\objectdegree}|_{U_{i_0\cdots i_p}},\mathcal{F}'^{\objectdegree+q}|_{U_{i_0\cdots i_p}}\right)
\end{align*}
and the \emph{differential} is given by
the total differential
$$
D
:=
\delta + (-1)^p d_{\mathcal{H}om}.
$$
%
%

\item
The \emph{composition}
is defined by
\begin{align*}
(f\comp f')_{i_0\cdots i_{p+p'}}^{\objectdegree}
:&=
(-1)^{qp'}
\left.
f_{i_0\cdots i_p}^{\objectdegree}
\right|_{U_{i_0\cdots i_{p+p'}}}
\comp
\left.
{f'}_{i_p\cdots i_{p+p'}}^{\objectdegree+q}
\right|_{U_{i_0\cdots i_{p+p'}}}
\\[3mm]
&\in
\Hom_{\Coh(U_{i_0\cdots i_{p+p'}})}\left(\mathcal{F}^{\objectdegree}|_{U_{i_0\cdots i_{p+p'}}},\mathcal{F}''^{\objectdegree+q+q'}|_{U_{i_0\cdots i_{p+p'}}}\right)
\end{align*}
for
\begin{align*}
f
&:=
\left(
f_{i_0\cdots i_p}^{\objectdegree}
\in
\Hom_{\Coh(U_{i_0\cdots i_p})}\left(\mathcal{F}^{\objectdegree}|_{U_{i_0\cdots i_p}},\mathcal{F}'^{\objectdegree+q}|_{U_{i_0\cdots i_p}}\right)
\right)_{
\substack{
i\in\mathbb{Z}
\\
i_0<\cdots<i_p
}
}
\in
\check{C}^p\big(
\mathfrak{U},\mathcal{H}om_{\mathcal{O}_X}^q(\mathcal{F},\mathcal{F}')
\big),
\\[2mm]
f'
&:=
\left(
{f'}_{i_0\cdots i_{p'}}^{\objectdegree}
\in
\Hom_{\Coh(U_{i_0\cdots i_{p'}})}\left(\mathcal{F}'^{\objectdegree}|_{U_{i_0\cdots i_{p'}}},\mathcal{F}''^{\objectdegree+q'}|_{U_{i_0\cdots i_{p'}}}\right)
\right)_{
\substack{
i\in\mathbb{Z}
\\
i_0<\cdots<i_p
}
}
\in
\check{C}^{p'}\big(
\mathfrak{U},\mathcal{H}om_{\mathcal{O}_X}^{q'}(\mathcal{F}',\mathcal{F}'')
\big).
\end{align*}

\item
The \emph{unit}
$
\id_{\mathcal{F}}
:=
\left(\left(\id_{\mathcal{F}}\right)_{i_0\cdots i_p}^{\objectdegree}\right)_{
\substack{
i\in\mathbb{Z}
\\
i_0<\cdots<i_p
}
}
\in
\check{C}^p\big(
\mathfrak{U},\mathcal{H}om_{\mathcal{O}_X}^q(\mathcal{F},\mathcal{F})
\big)
$
of $\mathcal F$ is given by
$$
\left(\id_{\mathcal{F}}\right)_{i_0}^{\objectdegree}
:=
\id_{\mathcal{F}^{\objectdegree}|_{U_{i_0}}}
\in
\Hom_{\Coh(U_{i})}\left(\mathcal{F}^{\objectdegree}|_{U_{i_0}},\mathcal{F}^{\objectdegree}|_{U_{i_0}}\right)
\quad\text{for each }
i_0\in I,\
\objectdegree\in\mathbb{Z}
$$
when
$p=0$, $q=0$,
and
zero otherwise.
\end{itemize}

These operations satisfy $D^2=0$, associativity,
and the Leibniz rule
$$
D(f\comp f')
=
Df\comp f'+(-1)^{|f|}f\comp Df'
$$
for homogeneous morphisms $f$ and $f'$.
\end{dfn}
We customly write $\Perf (X;\mathfrak{U})_{dg}$ for the full dg subcategory of perfect complexes (locally \emph{isomorphic} to bounded complexes of vector bundles).
\begin{prp}\label{prp:verifying-dg}
    $\Coh (X;\mathfrak{U})_{dg}$ is a dg category.
\end{prp}
\begin{proof}
    It remains for us to check that the composition $\comp$ is associative. This can be checked directly and with a great care of signs.
\end{proof}
\begin{prp}\label{prp:dg-enhancement}
$\Coh(X,\mathfrak{U})_{\dg}$
is a dg-enhancement of $D^b(\Coh(X))$.
\end{prp}
\begin{proof}
    Recall that for flasque $\mathcal{O}_X$-modules $\mathcal{F}$ we have $\check{H}^p (\mathfrak{U},\mathcal{F})=0$ for $p>0$ (\cite[Chapter 3, Proposition 4.3]{Hartshorne77}) and any injective $\mathcal{O}_X$-modules are flasque (\cite[Chapter 3, Lemma 2.4]{Hartshorne77}). Since global sections on affine varieties are exact functors, we see that the \v{C}ech complex $\check{C}^{\bullet} (\mathfrak{U};-)\colon\Coh (\mathcal{O}_X)\to\Ch (\Bbbk)$, regarded as a functor between abelian categories, is exact.\par
    Given a complex of sheaves $\mathcal{F}^{\bullet}\in\Ch (\mathcal{O}_X)$, let $\mathcal{F}^{\bullet}\to\mathcal{I}^{\bullet,\bullet}$ be the double complex of its injective resolution, then from exactness we see that $\check{C}^{\bullet} (\mathfrak{U};\mathcal{F}^{\bullet})\to\check{C}^{\bullet} (\mathfrak{U};\mathcal{I}^{\bullet,\bullet})$ is a quasi-isomorphism of chain complexes. We can therefore replace $\mathcal{F}^{\bullet}$ by the total complex $\Tot^{\bullet} (\mathcal{I}^{\bullet,\bullet})$ of injective resolutions.\par
    The \v{C}ech complex $\check{C}^{\bullet}\left(\Tot^{\bullet} (\mathcal{I}^{\bullet,\bullet})\right)$ can then be written as a double complex. Using the aforementioned fact about injective $\mathcal{O}_X$ modules we see that $\displaystyle\check{H}^{\bullet} (\mathcal{I}^{p,q})$ is supported in degree $0$ for any $(p,q)$. Spectral sequences associated to double complexes \cite[\href{https://stacks.math.columbia.edu/tag/012X}{Section 012X}]{stacks-project} provide us a spectral sequence $E^{p,q}_d$ so that $E_2^{p,q}=\displaystyle\check{H}^q\left(\mathfrak{U};H^p(\Tot^{\bullet}(\mathcal{I}^{\bullet,\bullet}))\right)$, which is supported at $q=0$, and therefore the spectral sequence degenerates at the first page and therefore we have isomorphisms\[
    H^p\left(\Tot^{\bullet}(\check{C}^{\bullet} (\mathfrak{U};\Tot^{\bullet} (\mathcal{I}^{\bullet,\bullet})))\right)\cong H^p (\Tot^{\bullet} (\mathcal{I}^{\bullet,\bullet}))\cong H^p (\mathcal{F}^{\bullet}).
    \]
    Now for any objects $\mathcal{F}^{\bullet}$ and $\mathcal{G}^{\bullet}$, we have isomorphisms of $\field$-modules\[
    H^0\hom^{\bullet}_{\Coh (X)} (\mathcal{F}^{\bullet},\mathcal{G}^{\bullet})_{dg}\simeq\bigoplus_{p\in\mathbb{Z}}\hom_{\Coh (X)} (\mathcal{F}^p,\mathcal{G}^p).
    \]Together with the fact that shifting by $[1]$ in $H\Mod (\Coh(X)_{dg}))$ preserves $\Coh (X;\mathfrak{U})_{dg}$, we finish the proof.
\end{proof}
\subsection[equivariant dg categories]{$G$-equivariant derived dg-categories.} 
Let $G$ be an affine gorup scheme and let
$\mathfrak U:=(U_i)_{i\in I}$
be an affine open cover of $X$
such that each $U_i$ is $G$-invariant.
For a $G$-equivariant $\mathcal O_X$-module $(\mathcal H,\Theta)$,
its \emph{$G$-invariant Čech cochain complex} with respect to
$\mathfrak U$ is
given by a collection of abelian groups
$$
\check C_G^p\big(\mathfrak U,(\mathcal H,\Theta)\big)
:=
\prod_{i_0<\cdots<i_p}
\mathcal H(U_{i_0\cdots i_p})^G
\quad\text{for }p\in\mathbb{Z}_{\ge0},
$$
where
$U_{i_0\cdots i_p}:= U_{i_0}\cap\cdots\cap U_{i_p}$,
with a differential
$
\delta:
\check{C}_G^p\big(\mathfrak{U},(\mathcal{H},\Theta)\big)
\rightarrow
\check{C}_G^{p+1}\big(\mathfrak{U},(\mathcal{H},\Theta)\big)
$
defined by
$$
(\delta s)_{i_0\cdots i_{p+1}}
:=
\sum_{j=0}^{p+1}(-1)^j
\left.
s_{i_0\cdots\widehat{i_j}\cdots i_{p+1}}
\right|_{U_{i_0\cdots i_{p+1}}}
\in
\mathcal{H}(U_{i_0\cdots i_{p+1}})^G
$$
for
$i_0<\cdots<i_{p+1}$,
where we denote by
$s=\left(s_{i_0\cdots i_p}\in\mathcal{H}(U_{i_0\cdots i_p})^G\right)_{i_0<\cdots<i_p}.$
\begin{dfn}\label{dfn:GEquivQCohdg}
The
\emph{$G$-equivariant coherent derived dg-category of $X$}
with respect to
$\mathfrak{U}$
is a dg-category
$$
\Coh_G(X;\mathfrak{U})_{\dg}
$$
defined in a similar way as Definition \ref{dfn:dg-enhancement-Coh}, where we replace $\check{C}^{\bullet} (\mathfrak{U},\mathcal{F}^{\bullet})$ by the $G$-invariant \v{C}ech cochain complex $\check{C}^{\bullet}_G (\mathfrak{U},(\mathcal{F}^{\bullet},\theta))$.
\end{dfn}
A direct verification shows that $\Coh_G (X;\mathfrak{U})_{dg}$ is a dg category, and is the dg enhancement of the derived category $D^b(\Coh_G (X))$ of quasi-coherent sheaves.
\begin{rmk}
    The existence of open affine cover of $X$ by $G$-invariants is a somewhat technical assumption. Due to the work of Sumihiro \cite{Sumihiro74} such cover exists when $X$ is normal and $G$ is a connected linear algebraic group, which is enough for the cases discussed in this paper.
\end{rmk}
Note that the strict $\widehat{G}$-action
$
\phi:
\widehat{G}
\rightarrow
\Aut_{\str}(\Coh_G(X))
$
on $\Coh_G(X)$
naturally induces a strict $\widehat{G}$-action
$
\phi:
\widehat{G}
\rightarrow
\Aut_{\str}(\Coh_G(X;\mathfrak{U})_{\dg})
$
on
$\Coh_G(X;\mathfrak{U})_{\dg}$
(Remark \ref{rmk:EmbeddingOfCochainComplexes}).

For a closed subgroup scheme
$H\subseteq G$,
the restriction functor
$\Res^G_H:\Coh_G(X)\rightarrow\Coh_H(X)$
(Definition \ref{defn:RestrictionOfEquivariantSheaves})
induces a dg-functor
$\Res_G^H:\Coh_G(X;\mathfrak{U})_{dg}\rightarrow\Coh_H(X;\mathfrak{U})_{dg}$,
and
we have a restricted strict $H^\perp$-action
$
\phi:
H^\perp
\rightarrow
\Aut_{\str}(\Coh_G(X)_{\dg}).
$
\begin{defn}[Restriction dg-functor]
For a closed subgroup scheme
$H\subseteq G$,
recall the restriction functor
$\Res^G_H:\Coh_G(X)\rightarrow\Coh_H(X)$
(Definition \ref{defn:RestrictionOfEquivariantSheaves}).
Suppose that there exists an affine open cover
$\mathfrak U:=(U_i)_{i\in I}$
of $X$
such that each $U_i$ is $G$-invariant.
Then each $U_i$ is automatically $H$-invariant.
We define the restriction dg-functor
between dg-categories
$$
\Res_G^H
:
\Coh_G(X;\mathfrak{U})_{\dg}
\rightarrow
\Coh_H(X;\mathfrak{U})_{\dg}
$$
whose mapping on objects and morphisms is given as follows:
\begin{itemize}
\item
Objects:
$
(\mathcal{F}^{\bullet},\theta^{\bullet})
\mapsto
(\mathcal{F}^{\bullet},\theta^{\bullet}|_H)
$
where
$$
(\mathcal{F}^{\bullet},\theta^{\bullet}|_H)
:=
\left(
\bigoplus_{\objectdegree\in\mathbb Z}
(\mathcal F^{\objectdegree},\theta^{\objectdegree}|_H),
d_{\mathcal F^\bullet}
:=
\left(d_{\mathcal{F}^{\bullet}}^i:\mathcal{F}^i\rightarrow\mathcal{F}^{i+1}\right)_{i\in\mathbb{Z}}
\right).
$$
Here,
each
$(\mathcal{F}^i,\theta^i|_H)$
is the image of
$(\mathcal{F}^i,\theta^i)$
under the restriction functor
$\Res^G_H:\Coh_G(X)\rightarrow\Coh_H(X)$,
and
$
d_{\mathcal{F}^{\bullet}}^i
\in
\Hom_{\Coh_H(X)}(\mathcal{F}^i,\mathcal{F}^{i+1})
$
is the image of the morphism
$
d_{\mathcal{F}^{\bullet}}^i
\in
\Hom_{\Coh_G(X)}(\mathcal{F}^i,\mathcal{F}^{i+1})
$
under the functor.
\item
Morphisms:
Recall the definition of the hom space
\[
\adjustbox{width=170mm}{%
  $\displaystyle
  \hom_{\Coh_G(X;\mathfrak U)_{\dg}}^{\homdegree}
  \left(
    (\mathcal F^\bullet,\theta^\bullet),
    (\mathcal F'^\bullet,\theta'^\bullet)
  \right)
  =
  \bigoplus_{p+q=\homdegree}
  \prod_{i_0<\cdots<i_p}
  \bigoplus_{\objectdegree\in\mathbb Z}
  \Hom_{\Coh_G(U_{i_0\cdots i_p})}
  \left(
    (\mathcal F^{\objectdegree},\theta^{\objectdegree})
      |_{U_{i_0\cdots i_p}},
    (\mathcal F'^{\objectdegree+q},\theta'^{\objectdegree+q})
      |_{U_{i_0\cdots i_p}}
  \right).$
}
\]
Each summand has a natural restriction
\[
\adjustbox{width=170mm}{$
\displaystyle
\Hom_{\Coh_G(U_{i_0\cdots i_p})}
\left(
(\mathcal F^{\objectdegree},\theta^{\objectdegree})
  |_{U_{i_0\cdots i_p}},
(\mathcal F'^{\objectdegree+q},\theta'^{\objectdegree+q})
  |_{U_{i_0\cdots i_p}}
\right)
\rightarrow
\Hom_{\Coh_H(U_{i_0\cdots i_p})}
\left(
(\mathcal F^{\objectdegree},\theta^{\objectdegree}|_H)
  |_{U_{i_0\cdots i_p}},
(\mathcal F'^{\objectdegree+q},\theta'^{\objectdegree+q}|_H)
  |_{U_{i_0\cdots i_p}}
\right),
$}
\]
given by the functor
$\Res^G_H:\Coh_G(U_{i_0\cdots i_p})\rightarrow\Coh_H(U_{i_0\cdots i_p})$.
Therefore,
there is a natural restriction
$$
\hom_{\Coh_G(X;\mathfrak U)_{\dg}}^{\homdegree}
\left(
(\mathcal F^\bullet,\theta^\bullet),
(\mathcal F'^\bullet,\theta'^\bullet)
\right)
\rightarrow
\hom_{\Coh_H(X;\mathfrak U)_{\dg}}^{\homdegree}
\left(
(\mathcal F^\bullet,\theta^\bullet),
(\mathcal F'^\bullet,\theta'^\bullet)
\right).
$$
\end{itemize}
\end{defn}

\begin{crl}
\label{crl:RestrictionFunctionEm}
Let $G$ be a diagonalizable group scheme acting on a scheme $X$,
and $H \subseteq G$ a closed subgroup scheme.
Let
$\mathfrak U:=(U_i)_{i\in I}$
be an affine open cover of $X$
such that each $U_i$ is $G$-invariant
(and hence $H$-invariant).\par
Then the restriction dg-functor
$
\Res_H^G:\Coh_G(X;\mathfrak{U})_{\dg}
\rightarrow
\Coh_H(X;\mathfrak{U})_{\dg}
$
induces a quasi-embedding of the orbit category:
$$
\Coh_G(X;\mathfrak{U})_{\dg}/H^\perp
\hookrightarrow
\Coh_H(X;\mathfrak{U})_{\dg}.
$$
\end{crl}

\begin{proof}
Let $(\mathcal{F},\theta)$ and $(\mathcal{F}',\theta')$
be two objects in $\Coh_G(X;\mathfrak{U})_{\dg}$.
$$
\hspace{-80mm}
\hom_{\Coh_H(X)_{\dg}}
\left((\mathcal{F}^{\bullet},\theta^{\bullet}|_H),(\mathcal{F}'^{\bullet},\theta'^{\bullet}|_H)\right)
$$
$$
=
\displaystyle
\bigoplus_{p+q=\homdegree}
\prod_{i_0<\cdots<i_p}
\bigoplus_{\objectdegree\in\mathbb Z}
\Hom_{\Coh_H(U_{i_0\cdots i_p})}
\left(
(\mathcal F^{\objectdegree},\theta^{\objectdegree}|_H)
  |_{U_{i_0\cdots i_p}},
(\mathcal F'^{\objectdegree+q},\theta'^{\objectdegree+q}|_H)
  |_{U_{i_0\cdots i_p}}
\right)
$$
By Proposition \ref{prop:HInvariantsOfHomSpace},
each summand is decomposed as
$$
\hspace{-30mm}
\Hom_{\Coh_H(U_{i_0\cdots i_p})}
\left(
(\mathcal F^{\objectdegree},\theta^{\objectdegree}|_H)
  |_{U_{i_0\cdots i_p}},
(\mathcal F'^{\objectdegree+q},\theta'^{\objectdegree+q}|_H)
  |_{U_{i_0\cdots i_p}}
\right)
$$
$$
=
\bigoplus_{\lambda\in H^{\perp}}
\Hom_{\Coh_G(U_{i_0\cdots i_p})}
\left(
(\mathcal F^{\objectdegree},\theta^{\objectdegree})
  |_{U_{i_0\cdots i_p}}(\lambda),
(\mathcal F'^{\objectdegree+q},\theta'^{\objectdegree+q}|_H)
  |_{U_{i_0\cdots i_p}}
\right)
$$
On the other hand,
by definition of orbit categories of dg-categories,
$$
\hspace{-80mm}
\hom_{\Coh_G(X)_{\dg}/H^{\perp}}
\left((\mathcal{F}^{\bullet},\theta^{\bullet}|_H),(\mathcal{F}'^{\bullet},\theta'^{\bullet}|_H)\right)
$$
$$
\hspace{-58mm}
=
\displaystyle
\bigoplus_{\lambda\in H^{\perp}}
\hom_{\Coh_G(X)_{\dg}}
\left((\mathcal{F}^{\bullet},\theta^{\bullet}|_H)(\lambda),(\mathcal{F}'^{\bullet},\theta'^{\bullet}|_H)\right)
$$
$$
=
\displaystyle
\bigoplus_{\lambda\in H^{\perp}}
\bigoplus_{p+q=\homdegree}
\prod_{i_0<\cdots<i_p}
\bigoplus_{\objectdegree\in\mathbb Z}
\Hom_{\Coh_G(U_{i_0\cdots i_p})}
\left(
(\mathcal F^{\objectdegree},\theta^{\objectdegree})
  |_{U_{i_0\cdots i_p}}(\lambda),
(\mathcal F'^{\objectdegree+q},\theta'^{\objectdegree+q}|_H)
  |_{U_{i_0\cdots i_p}}
\right).
$$
Therefore,
we have
$$
\hom_{\Coh_H(X)_{\dg}}
\left((\mathcal{F}^{\bullet},\theta^{\bullet}|_H),(\mathcal{F}'^{\bullet},\theta'^{\bullet}|_H)\right)
=
\hom_{\Coh_G(X)_{\dg}/H^{\perp}}
\left((\mathcal{F}^{\bullet},\theta^{\bullet}|_H),(\mathcal{F}'^{\bullet},\theta'^{\bullet}|_H)\right)
$$
as cochain complexes of $\field$-vector spaces.
\end{proof}

\subsection{Twisted Künneth formula.}
We end this section with the K\"unneth theorem. Suppose that $X$ and $Y$ are $G$- and $H$-schemes that admits an open cover $\mathfrak{U},\mathfrak{V}$ of $G$-invariant, resp. $H$-invariant open affine subschemes. Recall that $X\times Y$ admits an open cover by affine schemes of the form $U\times V\in\mathfrak{U}\times\mathfrak{V}$, which is clearly $G\times H$-invariant, 
Then we can define the equivariant K\"unneth functor
$$
\Perf_G(X;\mathfrak{U})_{\dg}
\otimes
\Perf_H(Y;\mathfrak{V})_{\dg}
\xrightarrow{\Ku}
\Perf_{G\times H}(X\times Y;\mathfrak{U}\times\mathfrak{V})_{\dg}
$$
on dg enhancements as follows:
\begin{defn}
$\Ku$ is a functor such that
\begin{itemize}
\item
Objects:
It maps an object
$
\big((\mathcal{F}^{\bullet},\theta^{\bullet}),(\mathcal{G}^{\bullet},\eta^{\bullet})\big)
$
to an object
$$
\Big(
\bigoplus_{n\in\mathbb{Z}}
\big(
(\mathcal{F}^{\bullet}\boxtimes\mathcal{G}^{\bullet})^n,
(\theta^{\bullet}\boxtimes\eta^{\bullet})^n
\big),
\
d_{\mathcal{F}^{\bullet}\boxtimes\mathcal{G}^{\bullet}}
:=
\big(
d_{\mathcal{F}^{\bullet}\boxtimes\mathcal{G}^{\bullet}}^n
:
(\mathcal{F}^{\bullet}\boxtimes\mathcal{G}^{\bullet})^n
\rightarrow
(\mathcal{F}^{\bullet}\boxtimes\mathcal{G}^{\bullet})^{n+1}
\big)_{n\in\mathbb{Z}}
\Big)
$$
(see Definition \ref{dfn:GEquivQCohdg}),
where
$$
\big(
(\mathcal{F}^{\bullet}\boxtimes\mathcal{G}^{\bullet})^n,
(\theta^{\bullet}\boxtimes\eta^{\bullet})^n
\big)
:=
\bigoplus_{i+j=n}
(
\mathcal{F}^i\boxtimes\mathcal{G}^j,
\
\theta^i\boxtimes\eta^j
)
\quad\text{for each $n\in\mathbb{Z}$}
$$
(note that
$
(\mathcal{F}^i,\theta^i)
\boxtimes
(\mathcal{G}^j,\eta^j)
=
(
\mathcal{F}^i\boxtimes\mathcal{G}^j,
\
\theta^i\boxtimes\eta^j
)
\in\Mod_{G\times H}(X\times Y)
$
was defined in
Definition \ref{defn:ExternalTensorProductOfEquivariantModules}),
and
the differential is determined on each summand
$\mathcal{F}^i\boxtimes\mathcal{G}^j$
(where $i+j=n$)
by
\[
d_{\mathcal F^\bullet\boxtimes\mathcal G^\bullet}^n
\big|_{\mathcal F^i\boxtimes\mathcal G^j}
:=
d_{\mathcal F}^i\boxtimes\id_{\mathcal G^j}
+
(-1)^i\id_{\mathcal F^i}\boxtimes d_{\mathcal G}^j
:
\mathcal{F}^i\boxtimes\mathcal{G}^j
\rightarrow
(\mathcal{F}^{i+1}\boxtimes\mathcal{G}^j)
\oplus
(\mathcal{F}^i\boxtimes\mathcal{G}^{j+1}).
\]

\item
Morphisms:
For homogeneous morphisms $f\in C_{\otimes}^{p+q}$ and $g\in C_{\otimes}^{r+s}$,
their tensor product $f\otimes g$ is mapped to $T^{\ast} (f\otimes g)\in C_{\times}^{p+q+r+s}$, where $T^{\ast}$ is the map defined in section \ref{subsec:Eilenberg-Zilber-contraction-Cech}.
\item Higher structures:
we define the $n$-th K\"unneth functor 
\begin{align*}
\Ku^n\colon &\hom_{\Perf (X;\mathfrak{U})_{dg}\otimes\Perf(Y;\mathfrak{V})_{dg}}\left(\left(
    (\mathcal{F}_{n-1}^{\bullet},\theta_{n-1}^{\bullet}),(\mathcal{G}_{n-1}^{\bullet},\eta_{n-1}^{\bullet})\right),
    \left((\mathcal{F}_n^{\bullet},\theta_n^{\bullet}),(\mathcal{G}_{n}^{\bullet},\eta_n^{\bullet})\right)\right)
    \otimes\dotsb\\ 
    \otimes
    &\hom_{\Perf (X;\mathfrak{U})_{dg}\otimes\Perf(Y;\mathfrak{V})_{dg}}\left(\left(
    (\mathcal{F}_0^{\bullet},\theta_0^{\bullet}),(\mathcal{G}_0^{\bullet},\eta_0^{\bullet})\right),
    \left((\mathcal{F}_1^{\bullet},\theta_1^{\bullet}),(\mathcal{G}_1^{\bullet},\eta_1^{\bullet})\right)\right)\\
    \to &\hom_{\Perf_{G\times H} (X\times Y;\mathfrak{U}\times\mathfrak{V})_{dg}}
    \left((\mathcal{F}_0^{\bullet},\theta_0^{\bullet})\boxtimes (\mathcal{G}_0^{\bullet},\eta_0^{\bullet}),(\mathcal{F}_n^{\bullet},\theta_n^{\bullet})\boxtimes (\mathcal{G}_n^{\bullet},\eta_n^{\bullet})\right)
\end{align*}
provided by the formula 
\begin{center}
\adjustbox{max width=\textwidth}{$\displaystyle
    \Ku^n (f_{n-1}\otimes g_{n-1},\dotsb, f_0\otimes g_0)=\sum_{\ell=1}^{n-1} (-1)^{(n-\ell -1)(t_{\ell}'+u_{\ell}'-1)}\mathbb{H}^{\ast}\left(
        \Ku^{n-\ell} (f_{n-1}\otimes g_{n-1},\dotsb, f_{\ell}\otimes g_{\ell})\circ\Ku^{\ell} (f_0\otimes g_0,\dotsb, f_{\ell-1}\otimes g_{\ell-1})
        \right),
$}\end{center}
where the operator $\mathbb{H}^{\ast}$ is defined in appendix \ref{sec:Eilenberg-Zilber-constraction}.
\end{itemize}
\end{defn}
We can check directly that
\begin{prp}
    The collection $(\Ku^n)_{n\in\mathbb{Z}_{\geq 0}}$ forms an $A_{\infty}$-functor from $\Perf_G (X;\mathfrak{U})_{dg}\otimes\Perf_H (Y;\mathfrak{V})_{dg}$ to $\Perf (X\times Y;\mathfrak{U}\times\mathfrak{V})_{dg}$.
\end{prp}
\begin{proof}
    The full verification will be lengthy, and we only present the verification for the first non-trivial relation to provide the pattern. From appendix \ref{sec:Eilenberg-Zilber-constraction} we see that $T^{\ast} S^{\ast} =\id$ and $S^{\ast} T^{\ast} -\id =\mathbb{H}^{\ast}d+d\mathbb{H}^{\ast}$, and hence 
    \begin{align*}
        d\Ku^2& (f_1\otimes g_1,f_0\otimes g_0)=d\mathbb{H}^{\ast} (\Ku (f_1\otimes g_1)\circ\Ku (f_0\otimes g_0))\\ 
        &=-\mathbb{H}^{\ast} (\Ku (d(f_1\otimes g_1))\circ \Ku (f_0\otimes g_0)) -(-1)^{\vert f_0\otimes g_0\vert }\mathbb{H}^{\ast} (\Ku(f_1\otimes g_1)\circ\Ku (d(f_0\otimes g_0)))\\
        &+T^{\ast} (S^{\ast}(\Ku (f_1\otimes g_1))\circ S^{\ast} (\Ku (f_0\otimes g_0))) -\Ku (f_1\otimes g_1)\circ\Ku (f_0\otimes g_0)\\ 
        &=-(-1)^{\vert f_0\otimes g_0\vert }\Ku^2 (d(f_1\otimes g_1),f_0\otimes g_0) - \Ku^2(f_1\otimes g_1,d(f_0\otimes g_0))\\ 
        &+\Ku((f_1\otimes g_1)\circ (f_0\otimes g_0))-\Ku (f_1\otimes g_1)\circ\Ku (f_0\otimes g_0),
    \end{align*}
    which is exactly the $A_{\infty}$-relation for $\Ku^2$. 
\end{proof}
Now we are ready to prove the K\"unneth embedding.
\begin{thm}\label{thm:Kunneth-embedding}
The prescribed functor
$$
\Ku\colon\Perf_G(X;\mathfrak{U})_{\dg}
\otimes
\Perf_H(Y;\mathfrak{V})_{\dg}
\hookrightarrow
\Perf_{G\times H}(X\times Y;\mathfrak{U}\times\mathfrak{V})_{\dg}
$$

$$
(\mathcal{F},\mathcal{G})
\quad\mapsto\quad
\mathcal{F}\boxtimes\mathcal{G}
~:=~
\pi_X^*\mathcal{F}
\underset{\mathcal{O}_{X\times Y}}{\otimes}
\pi_Y^*\mathcal{G}
$$
$$
f\otimes g
\quad\mapsto\quad
f\boxtimes g
~:=~
\pi_X^*f
\otimes
\pi_Y^*g
$$
is a quasi-embedding of dg-categories.
\end{thm}
\begin{proof}
    From Appendix \ref{sec:Eilenberg-Zilber-constraction} we know that $\Ku^1 =T^{\ast}$ is a quasi-equivalence of chain complexes with quasi-inverse $S^{\ast}$, and therefore after passing to homology category we know that $H^0 (\Ku)\colon H^0 (\Perf_G (X;\mathfrak{U})_{dg}\otimes\Perf_H (Y;\mathfrak{V})_{dg})\to H^0 (\Perf_{G\times H} (X\times Y;\mathfrak{U}\times\mathfrak{V})_{dg})$ is an embedding of categories. This proves that $\Ku$ is a quasi-embedding.
\end{proof}

\begin{crl}
    $\Ku$ induces a quasi-embedding
$$
\Perf_G(X)_{\dg}
\underset{\widehat{G}}{\otimes}
\Perf_G(Y)_{\dg}
\hookrightarrow
\Perf_G(X\times Y)_{\dg}
$$

$$
(\mathcal{F},\mathcal{G})
\quad\mapsto\quad
\mathcal{F}\boxtimes\mathcal{G}
~:=~
\pi_X^*\mathcal{F}
\underset{\mathcal{O}_{X\times Y}}{\otimes}
\pi_Y^*\mathcal{G}
$$
$$
f\otimes g
\quad\mapsto\quad
f\boxtimes g
~:=~
\pi_X^*f
\otimes
\pi_Y^*g.
$$
\end{crl}
\begin{proof}
From theorem \ref{thm:Kunneth-embedding} we get a quasi-embedding\\ 
$\Perf_G (X)_{dg}\otimes\Perf_G(Y)_{dg}\hookrightarrow\Perf_{G\times G} (X\times Y)_{dg}$. On the other hand, consider the diagonal embedding $\Delta\colon G\hookrightarrow G\times G$, its orthogonal subgroup action $\widehat{G}\cong (\Delta (G))^{\perp}\curvearrowright\Perf_{G\times G} (Y)$ preserves the subcategory $\Perf_G (X)_{dg}\otimes\Perf_G(Y)_{dg}$, and therefore from Corollary \ref{crl:RestrictionFunctionEm} and Proposition \ref{prp:twisted-tensor-equal-quotient} we get a quasi-embedding \begin{align*}
    \Perf_G (X)_{dg}\underset{\widehat{G}}{\otimes}\Perf_G(Y)_{dg}&\simeq\left(\Perf_G(X)\otimes\Perf_G(Y)_{dg}\right)\big/\widehat{G}\\ 
    &\hookrightarrow\Perf_{G\times G} (X\times Y)_{dg}\big/\widehat{G}\hookrightarrow\Perf_G (X\times Y)_{dg}.\qedhere
\end{align*}
\end{proof}
\appendix

\section{Eilenberg–-Zilber contraction for \v{C}ech complexes}\label{sec:Eilenberg-Zilber-constraction}

\subsection{Eilenberg--Zilber contraction for chain complexes}

In this subsection,
we follow
the presentation in \cite[\S 3]{Franz} of
the classical construction of Eilenberg--Mac Lane
\cite{EilenbergMacLane1953,EilenbergMacLane1954}
to construct the \emph{Eilenberg--Zilber maps}
for the chain complexes
associated with ordered sets and establish their relations.

For an ordered set $I$,
the chain complex
$
C_{\bullet}(I)
:=
\left(\bigoplus_{p\in\mathbb{Z}}C_p(I),\partial\right)
$
associated with $I$ is defined
by
$$
C_p(I)
:=
\mathbb{Z}\left<\left\{\left.(i_0,\dots,i_p)\in I^{p+1}\right|i_0<\cdots<i_p\right\}\right>
\quad\text{for }p\ge0,
$$
$C_p(I):=0$
for $p<0$,
and
\begin{align*}
\partial:C_p(I)
&\rightarrow
C_{p-1}(I)
\qquad
(p\ge 1),
\\
(i_0,\dots,i_p)
&\mapsto
\sum_{\nu=0}^p
(-1)^{\nu}
(i_0,\dots,\hat{i_{\nu}},\dots,i_p).
\end{align*}
We denote by
$
(i_{\sigma(0)},\dots,i_{\sigma(p)})
:=
\operatorname{sign}(\sigma)
(i_0,\dots,i_p)
$
for $i_0<\cdots<i_p$ in $I$\ and $\sigma\in S_{p+1}$,
and
$(i_0,\dots,i_p):=0$
if any two of $i_0,\dots,i_p$ are identical.

Given two ordered sets $I$ and $J$,
we consider two natural chain complexes
$$
C_{\bullet}(I\times J)
\quad\text{and}\quad
C_{\bullet}(I)\otimes C_{\bullet}(J).
$$
The first one is the chain complex associated with $I\times J$
equipped the lexicographic order.
The second one is the usual tensor product of chain complexes
$C_{\bullet}(I)$ and $C_{\bullet}(J)$,
whose components are explicitly given by
$$
(C_{\bullet}(I)\otimes C_{\bullet}(J))_t
:=
\bigoplus_{p+r=t}
\left(C_p(I)\otimes C_r(J)\right)
$$
with the boundary operator
\begin{align*}
\partial:
(C_{\bullet}(I)\otimes \partial_{\bullet}(J))_t
&\rightarrow
(C_{\bullet}(I)\otimes \partial_{\bullet}(J))_{t-1}
\\[2mm]
(i_0,\dots,i_p)\otimes(j_0,\dots,j_r)
&\mapsto
\partial((i_0,\dots,i_p))\otimes(j_0,\dots,j_r)
+
(-1)^p
(i_0,\dots,i_p)\otimes\partial((j_0,\dots,j_r)).
\end{align*}

We
claim that the two chain complexes are homotopically equivalent to each other,
by constructing explicit chain maps and homotopy,
called the \emph{Eilenberg--Zilber maps}:
\begin{equation}\label{eqn:EZMaps}
\begin{tikzcd}[arrow style=tikz,>=stealth,row sep=2.5em,column sep=3.5em]
C_{\bullet}(I)\otimes C_{\bullet}(J)
  \arrow[r, shift left, "S"]
&
C_{\bullet}(I\times J)
  \arrow[l, shift left, "T"]
  \arrow[loop right, <-, "\mathbb{H}"]
\end{tikzcd}
\end{equation}

\begin{defn}\label{defn:EZMaps}
(1)
The \emph{Shuffle map}
is defined by
\begin{center}
\resizebox{\textwidth}{!}{$\begin{aligned}
S:
C_{\bullet}(I)\otimes C_{\bullet}(J)
&\rightarrow
C_{\bullet}(I\times J)
\\[2mm]
(i_0,\dots,i_p)\otimes(j_0,\dots,j_r)
&\mapsto
\text{
\adjustbox{scale=0.8}{
$\sum_{
0=l_0\le l_1\le\cdots\le l_{p+1}=r
}
(-1)^{l_1+\cdots+l_p}
\big(
(i_0,j_0),\dots,(i_0,j_{l_0}),
(i_1,j_{l_1}),\dots,(i_1,j_{l_2}),
\dots,
(i_p,j_{l_p}),\dots,(i_p,j_r)
\big),
$}}
\end{aligned}$}
\end{center}
where exactly one subscript of $i$ or $j$ increases by $+1$ in each step.

\noindent
(2)
The \emph{Alexander–Whitney map}
is defined by
\begin{align*}
\begin{split}
T:
C_{\bullet}(I\times J)
&\rightarrow
C_{\bullet}(I)\otimes C_{\bullet}(J)
\\[0mm]
\big(\alpha_0=(i_0,j_0),\dots,\alpha_t=(i_t,j_t)\big)
&\mapsto
\sum_{p=0}^t
(i_0,\dots,i_p)\otimes(j_p,\dots,j_t).
\end{split}
\end{align*}

\end{defn}

It is straightforward to check that
$S$ and $T$ are chain maps
and $S\comp T$ is the identity:
$$
S\comp \partial=\partial \comp S,
\qquad
T\comp \partial=\partial \comp T,
\qquad
S\comp T=\id_{C_{\bullet}(I)\otimes C_{\bullet}(J)}.
$$

\begin{prop}\label{prp:chain-homotopy-TS}
The chain maps
$T\bullet S$
and
$\id_{C_{\bullet}(I\times J)}$
are homotopically equivalent to each other.
That is,
there is a chain homotopy
$$
\mathbb{H}:
C_{\bullet}(I\times J)
\rightarrow
C_{\bullet+1}(I\times J)
$$
between
$T\bullet S$
and
$\id_{C_{\bullet}(I\times J)}$,
satisfying
$$
\mathbb{H}\comp \partial+\partial\comp \mathbb{H}
=
T\comp S
-
\id_{C_{\bullet}(I\times J)}.
$$

Moreover, it additionally satisfies
$
S\comp \mathbb{H}=0.
$
\end{prop}

\subsection{Eilenberg--Zilber contraction for \v{C}ech complexes}\label{subsec:Eilenberg-Zilber-contraction-Cech}
Fix objects
$
(\mathcal{F}^{\bullet},\theta^{\bullet}),
$
$
(\mathcal{F}'^{\bullet},\theta'^{\bullet})
\in\Perf_G(X;\mathfrak{U})_{\dg}
$
and
$
(\mathcal{G}^{\bullet},\eta^{\bullet}),
$
$
(\mathcal{G}'^{\bullet},\eta'^{\bullet})
\in\Perf_H(Y;\mathfrak{V})_{\dg}.
$
Consider two cochain complexes
\begin{align*}
C_{\otimes}^{\bullet}
&:=
\hom_{
\Perf_G(X;\mathfrak{U})_{\dg}
\otimes
\Perf_H(Y;\mathfrak{V})_{\dg}
}^{\bullet}
\Big(
\big(
(\mathcal{F}^{\bullet},\theta^{\bullet}),
(\mathcal{G}^{\bullet},\eta^{\bullet})
\big),
\big(
(\mathcal{F}'^{\bullet},\theta'^{\bullet}),
(\mathcal{G}'^{\bullet},\eta'^{\bullet})
\big)
\Big),
\\[2mm]
C_{\times}^{\bullet}
&:=
\hom_{\Perf_{G\times H}(X\times Y;\mathfrak{U}\times\mathfrak{V})_{\dg}}^{\bullet}
\Big(
(\mathcal{F}^{\bullet},\theta^{\bullet})
\boxtimes
(\mathcal{G}^{\bullet},\eta^{\bullet}),
(\mathcal{F}'^{\bullet},\theta'^{\bullet})
\boxtimes
(\mathcal{G}'^{\bullet},\eta'^{\bullet})
\Big).
\end{align*}
We denote homogeneous elements in
\begin{align*}
\begin{split}
&
\hom_{\Perf_{G}(X;\mathfrak{U})_{\dg}}^l
\big(
(\mathcal{F}^{\bullet},\theta^{\bullet}),
(\mathcal{F}'^{\bullet},\theta'^{\bullet})
\big)
=
\bigoplus_{p+q=l}
\check{C}_{G}^p
\Big(
\mathfrak{U},
\mathcal{H}om_{\mathcal{O}_X}^q
\big(
(\mathcal{F}^{\bullet},\theta^{\bullet}),
(\mathcal{F}'^{\bullet},\theta'^{\bullet})
\big)
\Big),
\\[2mm]
&
\hom_{\Perf_{H}(Y;\mathfrak{V})_{\dg}}^m
\big(
(\mathcal{G}^{\bullet},\eta^{\bullet}),
(\mathcal{G}'^{\bullet},\eta'^{\bullet})
\big)
=
\bigoplus_{r+s=m}
\check{C}_{H}^r
\Big(
\mathfrak{V},
\mathcal{H}om_{\mathcal{O}_Y}^s
\big(
(\mathcal{G}^{\bullet},\eta^{\bullet}),
(\mathcal{G}'^{\bullet},\eta'^{\bullet})
\big)
\Big),
\\[2mm]
&
C_{\times}^n
=
\bigoplus_{t+u=n}
\check{C}_{G\times H}^t
\Big(
\mathfrak{U}\times\mathfrak{V},
\mathcal{H}om_{\mathcal{O}_{X\times Y}}^u
\big(
(\mathcal{F}^{\bullet},\theta^{\bullet})
\boxtimes
(\mathcal{G}^{\bullet},\eta^{\bullet}),
(\mathcal{F}'^{\bullet},\theta'^{\bullet})
\boxtimes
(\mathcal{G}'^{\bullet},\eta'^{\bullet})
\big)
\Big)
\end{split}
\end{align*}
of fixed bi-degree $(p,q)$, $(r,s)$, $(t,u)$
by
$$
f:=
\Big(f_{i_0\cdots i_p}^i\Big)
_{\substack{
i\in\mathbb{Z}
\\
i_0<\cdots<i_p
}},
\quad
g:=
\Big(g_{j_0\cdots j_r}^j\Big)
_{\substack{
j\in\mathbb{Z}
\\
j_0<\cdots<j_r
}},
\quad
h:=
\Big(h_{\alpha_0\cdots \alpha_t}^k\Big)
_{\substack{
k\in\mathbb{Z}
\\
\alpha_0<\cdots<\alpha_t
}},
$$
respectively,
under the identifications
\begin{align*}
\begin{split}
&
\check{C}_{G}^p
\Big(
\mathfrak{U},
\mathcal{H}om_{\mathcal{O}_X}^q
\big(
(\mathcal{F}^{\bullet},\theta^{\bullet}),
(\mathcal{F}'^{\bullet},\theta'^{\bullet})
\big)
\Big)
\cong
\prod_{i_0<\cdots<i_p}
\bigoplus_{i\in\mathbb{Z}}
\Hom_{\VB_{G}(U_{i_0\cdots i_p})}
\big(
(\mathcal{F}^i,\theta^i)|_{U_{i_0\cdots i_p}}),
(\mathcal{F}'^i,\theta'^{i+q}|_{U_{i_0\cdots i_p}})
\big),
\\[2mm]
&
\check{C}_{H}^r
\Big(
\mathfrak{V},
\mathcal{H}om_{\mathcal{O}_Y}^s
\big(
(\mathcal{G}^{\bullet},\eta^{\bullet}),
(\mathcal{G}'^{\bullet},\eta'^{\bullet})
\big)
\Big)
\cong
\prod_{j_0<\cdots<j_r}
\bigoplus_{j\in\mathbb{Z}}
\Hom_{\VB_{H}(V_{j_0\cdots j_p})}
\big(
(\mathcal{G}^j,\eta^j)|_{V_{j_0\cdots j_r}}),
(\mathcal{G}'^j,\eta'^{j+s}|_{V_{j_0\cdots j_r}})
\big),
\\[2mm]
&
\check{C}_{G\times H}^t
\Big(
\mathfrak{U}\times\mathfrak{V},
\mathcal{H}om_{\mathcal{O}_{X\times Y}}^u
\big(
(\mathcal{F}^{\bullet},\theta^{\bullet})
\boxtimes
(\mathcal{G}^{\bullet},\eta^{\bullet}),
(\mathcal{F}'^{\bullet},\theta'^{\bullet})
\boxtimes
(\mathcal{G}'^{\bullet},\eta'^{\bullet})
\big)
\Big)
\\[2mm]
&
\hspace{-5mm}
\cong
\prod_{\alpha_0<\cdots<\alpha_t}
\bigoplus_{k\in\mathbb{Z}}
\bigoplus_{\substack{
i+j=k
\\[1mm]
q+s=u
}}
\Hom_{\VB_{G\times H}(W_{\alpha_0\cdots \alpha_t})}
\Big(
\big(
(\mathcal{F}^i,\theta^i)
\boxtimes
(\mathcal{G}^j,\eta^j)
\big)|_{W_{\alpha_0\cdots \alpha_t}},
\big(
(\mathcal{F}'^{i+q},\theta'^{i+q})
\boxtimes
(\mathcal{G}'^{j+s},\eta'^{j+s})
\big)|_{W_{\alpha_0\cdots \alpha_t}})
\Big)
\end{split}
\end{align*}
We denote a homogeneous element in
\begin{center}
\resizebox{\textwidth}{!}{$\begin{aligned}
C_{\otimes}^{n}
&=
\bigoplus_{l+m=n}
\Big(
\hom_{\Perf_{G}(X;\mathfrak{U})_{\dg}}^l
\big(
(\mathcal{F}^{\bullet},\theta^{\bullet}),
(\mathcal{F}'^{\bullet},\theta'^{\bullet})
\big)
\underset{\field}{\otimes}
\hom_{\Perf_{H}(Y;\mathfrak{V})_{\dg}}^m
\big(
(\mathcal{G}^{\bullet},\eta^{\bullet}),
(\mathcal{G}'^{\bullet},\eta'^{\bullet})
\big)
\Big)
\\[2mm]
&
\hspace{-5mm}
\cong
\bigoplus_{l+m=n}
\Bigg(
\bigg(
\bigoplus_{p+q=l}
\check{C}_{G}^p
\Big(
\mathfrak{U},
\mathcal{H}om_{\mathcal{O}_X}^q
\big(
(\mathcal{F}^{\bullet},\theta^{\bullet}),
(\mathcal{F}'^{\bullet},\theta'^{\bullet})
\big)
\Big)
\bigg)
\underset{\field}{\otimes}
\bigg(
\bigoplus_{r+s=m}
\check{C}_{H}^r
\Big(
\mathfrak{V},
\mathcal{H}om_{\mathcal{O}_Y}^s
\big(
(\mathcal{G}^{\bullet},\eta^{\bullet}),
(\mathcal{G}'^{\bullet},\eta'^{\bullet})
\big)
\Big)
\bigg)
\Bigg)
\\[2mm]
&
\hspace{-5mm}
\cong
\bigoplus_{l+m=n}
\bigoplus_{\substack{
p+q=l
\\[1mm]
r+s=m
}}
\prod_{\substack{
i_0<\cdots<i_p
\\[0mm]
j_0<\cdots<j_r
}}
\bigoplus_{\substack{
i\in\mathbb{Z}
\\[1mm]
j\in\mathbb{Z}
}}
\Big(
\Hom_{\VB_{G}(U_{i_0\cdots i_p})}
\big(
(\mathcal{F}^i,\theta^i)|_{U_{i_0\cdots i_p}},
(\mathcal{F}'^{i+q},\theta'^{i+q})|_{U_{i_0\cdots i_p}}
\big)
\\
&\phantom{
\cong
\bigoplus_{l+m=n}
\bigoplus_{\substack{
p+q=l
\\[1mm]
r+s=m
}}
\prod_{\substack{
i_0<\cdots<i_p
\\[0mm]
j_0<\cdots<j_r
}}
\bigoplus_{\substack{
i\in\mathbb{Z}
\\[1mm]
j\in\mathbb{Z}
}}
\Big(
}
\underset{\field}{\otimes}
\Hom_{\VB_{H}(V_{j_0\cdots j_r})}
\big(
(\mathcal{G}^j,\eta^j)|_{V_{j_0\cdots j_r}},
(\mathcal{G}'^{j+s},\eta'^{j+s})|_{V_{j_0\cdots j_s}}
\big)
\Big)
\end{aligned}$}
\end{center}
of multi-degree
$(p,q,r,s)$
by
$$
f\otimes g
=
\Big(
f_{i_0\cdots i_p}^i
\otimes
g_{j_0\cdots j_r}^j
\Big)
_{\substack{
i,j\in\mathbb{Z}
\\
i_0<\cdots<i_p
\\
j_0<\cdots<j_r
}}.
$$

The Eilenberg--Zilber maps
$S$, $T$ and $H$
in
(\ref{eqn:EZMaps})
induce
cochain maps
$S^*$, $T^*$ and$H^*$
between $C_{\otimes}^{\bullet}$ and $C_{\times}^{\bullet}$:

$$
\begin{tikzcd}[arrow style=tikz,>=stealth,row sep=2.5em,column sep=3.5em]
C_{\otimes}^{\bullet}
  \arrow[r, shift left, <-, "S^*"]
&
C_{\times}^{\bullet}
  \arrow[l, shift left, <-, "T^*"]
  \arrow[loop right, "\mathbb{H}^{\ast}"]
\end{tikzcd}
$$

\begin{dfn}

(1)
The chain map 
$
T:
C_{\bullet}(I\times J)
\rightarrow
C_{\bullet}(I)\otimes C_{\bullet}(J)
$
induces a map between cochain complexes
\begin{align*}
T^*:
C_{\otimes}^n
&\rightarrow
C_{\times}^n
\\
\quad
f\otimes g
=
\Big(
f_{i_0\cdots i_p}^i
\otimes
g_{j_0\cdots j_r}^j
\Big)
_{\substack{
i,j\in\mathbb{Z}
\\
i_0<\cdots<i_p
\\
j_0<\cdots<j_r
}}
&\mapsto
T^*(f\otimes g)
:=
\Big(
T^*(f\otimes g)_{\alpha_0\cdots\alpha_t}^{i,j}
\Big)
_{\substack{
i,j\in\mathbb{Z}
\\
i_0<\cdots<i_p
\\
j_0<\cdots<j_r
}},
\end{align*}
whose components are defined by
$$
T^*(f\otimes g)_{\alpha_0\cdots\alpha_t}^{i,j}
:=
\sum_{\substack{
i_0<\cdots<i_p
\\
j_0<\cdots<j_r
}}
(-1)^{qr+qs+si}
T_{\alpha_0\cdots\alpha_t}^{i_0\cdots i_p;j_0\cdots j_r}
\big(
f_{i_0\cdots i_p}^i
\boxtimes
g_{j_0\cdots j_r}^j
\big)\big|_{W_{\alpha_0\cdots\alpha_t}},
$$
where the coefficients
$
T_{\alpha_0\cdots\alpha_t}^{i_0\cdots i_p;j_0\cdots j_r}
\in\mathbb{Z}
$
are determined by the formula
$$
T\big((\alpha_0,\dots,\alpha_t)\big)
=
\sum_{\substack{
i_0<\cdots<i_p
\\
j_0<\cdots<j_r
}}
T_{\alpha_0\cdots\alpha_t}^{i_0\cdots i_p;j_0\cdots j_r}
(i_0,\dots,i_p)\otimes(j_0,\dots,j_r).
$$

It follows from the definition of $T$
(Definition \ref{defn:EZMaps}.(2))
that
$$
T_{\alpha_0\cdots\alpha_t}^{i_0\cdots i_p;j_0\cdots j_r}
=
\begin{cases}
1
&
\text{if }
(\alpha_0,\dots,\alpha_t)
=
\big((i_0,*),\dots,(i_{p-1},*),(i_p,j_0),(*,j_1),\dots,(*,j_r)\big),
\\
0
&
\text{otherwise},
\end{cases}
$$
and therefore we have
$$
T^*(f\otimes g)_{\alpha_0\cdots\alpha_t}^{i,j}
=
(-1)^{qr+qs+si}
\big(
f_{i_0\cdots i_p}^i
\boxtimes
g_{j_p\cdots j_t}^j
\big)\big|_{W_{\alpha_0\cdots\alpha_t}}
$$
where
$(\alpha_0,\dots,\alpha_t)=\big((i_0,j_0),\dots,(i_t,j_t)\big)$.

(2)
The chain map
$
S:
C_{\bullet}(I)\otimes C_{\bullet}(J)
\rightarrow
C_{\bullet}(I\times J)
$
induces
a map between cochain complexes
\begin{align*}
S^*:
C_{\times}^n
&\rightarrow
C_{\otimes}^n
\\
\quad
h
=
\Big(
h_{\alpha_0\cdots\alpha_t}^{i,j}
\Big)
_{\substack{
i,j\in\mathbb{Z}
\\
\alpha_0<\cdots<\alpha_t
}}
&\mapsto
S^*h
:=
\Big(
(S^*h)_{i_0\cdots i_p;j_0\cdots j_r}^{i,j}
\Big)
_{\substack{
i,j\in\mathbb{Z}\phantom{,\ p+r=t}
\\
p,r\in\mathbb{Z},\ p+r=t
\\
i_0<\cdots<i_p\quad\
\\
j_0<\cdots<j_r\quad\
}},
\end{align*}
whose components are defined by
$$
(S^*h)_{i_0\cdots i_p;j_0\cdots j_r}^{i,j}
:=
\sum_{\substack{
\alpha_0\cdots\alpha_t
}}
(-1)^{qr+qs+si}
S_{i_0\cdots i_p;j_0\cdots j_r}^{\alpha_0\cdots\alpha_t}
\boxtimes^{-1}\big(h_{\alpha_0\cdots\alpha_t}^{i,j}\big),
$$
where the coefficients
$
S_{i_0\cdots i_p;j_0\cdots j_r}^{\alpha_0\cdots\alpha_t}
\in\mathbb{Z}
$
are determined by the formula
$$
S\big((i_0,\dots,i_p)\otimes(j_0,\dots,j_r)\big)
=
\sum_{\substack{
\alpha_0<\cdots<\alpha_t
}}
S_{i_0\cdots i_p;j_0\cdots j_r}^{\alpha_0\cdots\alpha_t}
(\alpha_0,\dots,\alpha_t)
$$
and\\
\adjustbox{max width=\textwidth}{$\displaystyle
\boxtimes:
\Hom_{\VB_G(U_{i_0\cdots i_p})}
\big(
(\mathcal{F}^i,\theta^i)|_{U_{i_0\cdots i_p}},
(\mathcal{F}'^{i+q},\theta'^{i+q})|_{U_{i_0\cdots i_p}}
\big)
\underset{\field}{\otimes}
\Hom_{\VB_H(V_{j_0\cdots j_r})}
\big(
(\mathcal{G}^j,\eta^j)|_{V_{j_0\cdots j_r}},
(\mathcal{G}'^{j+s},\eta'^{j+s})|_{V_{j_0\cdots j_r}}
\big)
$}\\
\adjustbox{max width=\textwidth}{
$\displaystyle
\quad
\xrightarrow{\cong}
\Hom_{\VB_{G\times H}(U_{i_0\cdots i_p}\times V_{j_0\cdots j_r})}
\Big(
\big(
(\mathcal{F}^i,\theta^i)
\boxtimes
(\mathcal{G}^j,\eta^j)
\big)
\big|_{U_{i_0\cdots i_p}\times V_{j_0\cdots j_r}},
\big(
(\mathcal{F}'^{i+q},\theta'^{i+q})
\boxtimes
(\mathcal{G}'^{j+s},\eta'^{j+s}
\big)
\big|_{U_{i_0\cdots i_p}\times V_{j_0\cdots j_r}}
\Big)
$}\\
is the isomorphism given in Proposition \ref{prop:BoxTensorFunctorFullyFaithful}.

It follows from the definition of $S$
(Definition \ref{defn:EZMaps}.(1))
that
$$
(S^*h)_{i_0\cdots i_p;j_0\cdots j_r}^{i,j}
=
(-1)^{qr+qs+si}
\boxtimes^{-1}
\bigg(
\sum_{
0=l_0\le l_1\le\cdots l_{p+1}=r
}
(-1)^{l_1+\cdots+l_p}
h_{
(i_0,j_0),\dots,(i_0,j_{l_0}),
(i_1,j_{l_1}),\dots,(i_1,j_{l_2}),
\dots,
(i_p,j_{l_p}),\dots,(i_p,j_r)
}
^{i,j}
\bigg).
$$

\noindent
(3)
The chain homotopy
$
\mathbb{H}:
C_{\bullet}(I\times J)
\rightarrow
C_{\bullet}(I\times J)
$
induces a map between cochain complexes
$$
\mathbb{H}^{\ast}\colon C_{\times}^{\bullet}\to C_{\times}^{\bullet -1}
$$where for each $f\coloneqq\left(f_{\beta_0\dotsb\beta_{t+1}^k}\right)_
{
    \begin{subarray}{c}
        k\in\mathbb{Z}\\ 
        \beta_0<\dotsb<\beta_{t+1}
    \end{subarray}
    }$, we have$$
(\mathbb{H}^{\ast} f)\coloneqq\left((\mathbb{H}^{\ast} f)^k_{\alpha_0\dotsb\alpha_t}\coloneqq\sum_{\beta_0<\dotsb<\beta_{t+1}}\left. \mathbb{H}_{\alpha_0\dotsb\alpha_t}^{\beta_0\dotsb\beta_{t+1}} f_{\beta_0\dotsb\beta_{t+1}}^k\right|_{W_{\alpha_0\dotsb\alpha_t}}\right)
_{
    \begin{subarray}{c}
        k\in\mathbb{Z}\\ 
        \alpha_0<\dotsb<\alpha_t
    \end{subarray}
}
$$where we write \[
    \mathbb{H}((\alpha_0,\dotsb,\alpha_t))\eqqcolon\sum_{\beta_0<\dotsb<\beta_{t+1}}\mathbb{H}_{\alpha_0\dotsb\alpha_t}^{\beta_0\dotsb\beta_{t+1}} (\beta_0,\dotsb ,\beta_{t+1}).
\]
\end{dfn}
Now we can verify that
\begin{prp}\label{prp:STHd-relation}
    \begin{itemize}
        \item Write $\delta$ for the \v{C}ech differential, then $\delta\mathbb{H}^{\ast} +\mathbb{H}^{\ast} \delta=S^{\ast} T^{\ast} -\id$.
        \item $S^{\ast}$ preserves compositions, i.e. $S^{\ast} (f\comp g)=S^{\ast} f\comp S^{\ast} g$.
        \item Write $d_{\mathcal{H}^{\bullet}}$ to be the differential for the internal hom complex, then $d_{\mathcal{H}^{\bullet}}\mathbb{H}^{\ast} =\mathbb{H}^{\ast}d_{\mathcal{H}^{\bullet}}$.
        \item $dT^{\ast} =T^{\ast} d$.
        \item $T^{\ast} S^{\ast} =\id$.
    \end{itemize}
\end{prp}
\begin{proof}
    These are all direct verifications, using the definition and relations in Proposition \ref{prp:STHd-relation}.
\end{proof}



\bibliographystyle{alphaurl}
\bibliography{main}
\end{document}